\documentclass[11pt]{amsart}
\usepackage{a4wide}
\usepackage{verbatim}
\usepackage[T1]{fontenc}
\usepackage[utf8]{inputenc}
\usepackage{amsfonts}
\usepackage{amsthm}
\usepackage{amsmath}
\usepackage{bm}
\usepackage{amssymb,amsbsy,amsopn}
\usepackage{mathrsfs}
\usepackage{graphicx,color,float,subcaption,caption,algorithm}
\usepackage{mathptmx,bbm}
\usepackage[bookmarks=true,bookmarksnumbered=true,bookmarkstype=toc,colorlinks, linkcolor={blue}, citecolor={blue}, urlcolor={blue}]{hyperref}
\usepackage{dsfont}
\usepackage{leftidx}
\usepackage{empheq}
\usepackage{enumerate}
\usepackage[parfill]{parskip}
\usepackage{setspace}
\usepackage{xfrac}
\usepackage{esvect}
\usepackage{multirow}
\usepackage{hhline}
\usepackage{ulem}
\allowdisplaybreaks

\usepackage{tikz}
\usetikzlibrary{arrows,automata,shapes}
\usetikzlibrary{positioning}
\newtheorem{definition}{Definition}[section] 
\newtheorem{theorem}{Theorem}[section]

\newtheorem{lemma}[theorem]{Lemma}
\newtheorem{remark}[theorem]{Remark}
\newtheorem{proposition}[theorem]{Proposition}

\newcommand{\ie}{\setlength{\parskip}{0cm} \setlength{\itemsep}{0cm}}
\DeclareMathAlphabet{\mathcal}{OMS}{cmsy}{m}{n}

\usepackage{amscd}

\newcommand{\la}{\left(}
\newcommand{\ra}{\right)}
\newcommand{\lb}{\left\langle}
\newcommand{\rb}{\right\rangle}

\def\E{\mathbb E}
\def\P{\mathbb P}
\def\b{\mathbb}
\def\bal{\begin{equation*}\begin{aligned}} 
		\def\eal{\end{aligned}\end{equation*}}
\def\be{\begin{equation}\label}
	\def\ee{\end{equation}}
\def\bd{\begin{definition}\label}
	\def\ed{\end{definition}}
\def\bt{\begin{theorem}\label}
	\def\et{\end{theorem}}
\def\bl{\begin{lemma}\label}
	\def\el{\end{lemma}}

\def\R{\mathbb R}

\def\L{\mathbb L}
\def\H{\mathbb H}
\def\V{\mathbb V}

\date{}
\begin{document}
	
\title{Numerical method and Error estimate for stochastic Landau--Lifshitz--Bloch equation}
\author{Ben Goldys}
\address{School of Mathematics and Statistics\\
	The University of Sydney\\
	Sydney 2006, Australia} 
\email{beniamin.goldys@sydney.edu.au}
\author{Chunxi Jiao}
\address{School of Mathematics and Statistics\\
	The University of Sydney\\
	Sydney 2006, Australia 
	\newline
	School of Mathematics and Statistics \\
	The University of New South Wales \\
	Sydney 2052, Australia} 
\email{chunxi.jiao@unsw.edu.au}
\author{Kim-Ngan Le}
\address{School of Mathematics, Monash University, Clayton, Victoria 3800, Australia}
\email{ngan.le@monash.edu}
\thanks{This work was supported by the Australian Research Council Project DP200101866}
\keywords{}
\subjclass[2020]{}

\begin{abstract}
	In this paper we study numerical methods for solving a system of quasilinear stochastic partial differential equations known as the stochastic Landau-Lifshitz-Bloch (LLB) equation on a bounded domain in $\R^d$ for $d=1,2$. Our main results are estimates of the rate of convergence of the Finite Element Method to the solutions of stochastic LLB. 	
	To overcome the lack of regularity of the solution in the case $d=2$, we propose a Finite Element scheme for a regularised version of the equation. We then obtain error estimates of numerical solutions and for the solution of the regularised equation as well as the rate of convergence of this solution to the solution of the stochastic LLB equation. As a consequence, the convergence in probability of the approximate solutions to the solution of  the stochastic LLB equation is derived. To the best of our knowledge this is the first result on error estimates for a system of stochastic quasilinear partial differential equations. A stronger result is obtained in the case $d=1$ due to a new regularity result for the LLB equation which allows us to avoid regularisation.	
\end{abstract}

\maketitle
\tableofcontents

\section{Introduction}
	The Landau-Lifshitz-Bloch (LLB) equation \eqref{eq: sLLB} is a model overcoming a significant limitation of standard models of magnetisation (such as Landau-Lifshitz equation) which are valid only at temperatures close to the Curie point. Mathematical analysis of this equation was initiated in \cite{Le2016}. 
	A stochastic version of the LLB equation was introduced in physics papers ~\cite{Garanin2004},~\cite{Evan2012}. 
	In those works a deterministic version of equation \eqref{eq: sLLB} was modified in order to incorporate random fluctuations into the dynamics of the magnetisation and to describe noise-induced transitions between equilibrium states of the ferromagnet.
	
	The stochastic LLB equation was studied in~\cite{BGL_sLLB}, where the existence of strong (in the sense of PDEs) martingale solutions was proved for domains in $\R^d$, $d=1,2,3$. Furthermore, uniqueness of pathwise solutions and the existence of invariant measures is obtained in cases $d=1,2$. A similar result on the existence of martingale solution on $3$ dimensional domains is also obtained in~\cite{JiangWang2019}.  
	
	Let $D \subset \R^d$, $d=1,2$, be an open bounded domain and let $ m: [0,T] \times D \to \R^3 $ denote the magnetisation vector. 
	Let us note that contrary to the case of the Landau-Lifschitz-Gilbert equation $m$ need not take values in the unit sphere, hence the control of the size of $u$ is not automatically assured. 
	When the temperature exceeds the Curie temperature, the stochastic LLB equation has the form
	\begin{equation}\label{eq: sLLB}
		dm 
		= \la \kappa_1 \Delta m + \gamma m \times \Delta m - \kappa_2 (1+\mu |m|^2) m \ra dt 
		+ \sum_{k=1}^\infty \la \kappa_1 g_k + \gamma m \times g_k \ra \circ dW_k, 
	\end{equation}
	with initial condition $ m(0) = m_0 \in \H^1 $ and homogeneous Neumann boundary condition. 	
	The constants $ \kappa_1, \kappa_2, \gamma, \mu > 0 $ and $ \{ W_k \}_k $ is a family of independent real-valued Wiener processes. 
	We can view the noise part in \eqref{eq: sLLB} as an $ \H^1 $-valued Wiener process
	\begin{equation}\label{def: W in H1}
		W(t) = \sum_{k=1}^\infty g_k W_k(t)\,.
	\end{equation} 
	For simplicity we let $ \kappa_1 = \kappa_2 = \mu = \gamma = 1 $.
	
	In this paper, we propose a numerical scheme for computing approximate solutions of the stochastic LLB equation ~\eqref{eq: sLLB}. 
	The equation is quasilinear because of the nonlinear term $ m \times \Delta m$ which causes difficulties in the analysis of this equation, and in particular in proving the order of convergence. 
\par
In one-dimensional domains we obtain a global error estimate. In order to obtain such global estimate we prove in Theorem \ref{Lemma: m 1D regularity} a new regularity result for the LLB equation in one dimension, which is of independent interest. 
\par
In the two-dimensional case such regularity results for the stochastic LLB equation are not known. Due to the lack of regularity of the solution, we provide local error estimates and local rates of convergence. Here, we were inspired by similar results obtained in \cite{BreitProhl} for the stochastic Navier-Stokes equation.
	To overcome the lack of regularity in $2$-dimensional case, we introduce a regularised version of ~\eqref{eq: sLLB} 
	\begin{equation}\label{eq: sLLB reg}
		\begin{aligned}
			dm^\epsilon 
			= \la -\epsilon \Delta^2 m^\epsilon + \Delta m^\epsilon  + m^\epsilon \times \Delta m^\epsilon - (1+ |m^\epsilon|^2) m^\epsilon \ra dt 
			+ \sum_{k=1}^\infty \la g_k + m^\epsilon \times g_k \ra \circ dW_k, 
		\end{aligned}
	\end{equation}
	with  $ \epsilon \in (0,1) $,  initial condition $ m^\epsilon(0) = m_0 \in \H^1 $ and homogeneous Neumann boundary condition.  
	We will require that the noise $W$ is a $ \H^3 $-valued Wiener process in the form of \eqref{def: W in H1} that satisfies the condition 
	\begin{equation}\label{A: g bound H2}
		C_g:=\sum_{k=1}^\infty |g_k|^2_{{\b W}^{2,\infty} \cap \H^3} 
		< \infty.
	\end{equation}
	Thanks to the regularising term $-\epsilon \Delta^2 m^\epsilon $, the unique solution of~\eqref{eq: sLLB reg} is smooth enough to deduce error estimates. Since the arguments used in the case $ d=1 $ are similar to, but simpler than in the case $ d=2 $, we first focus on the case $ d=2 $ in Sections \ref{Section: regularisation} and \ref{Section: FEM m-epsilon}. The existence and regularity of solution $m^\epsilon$ is proved in Section~\ref{Section: regularisation} by using the Faedo-Galerkin approximations and the method of compactness. Furthermore, the convergence of $m^\epsilon$ to $m$ is also proved in this Section. 
	In Section~\ref{Section: FEM m-epsilon}, we propose a Finite Element scheme for the regularised equation~\eqref{eq: sLLB reg} and prove stabilities of its approximate solutions as well as error estimates and convergences in probability of the approximate solutions to the solution of~\eqref{eq: sLLB}.
	In Section~\ref{Section: 1D case} we consider the case $ d=1 $. We start with Theorem \ref{Lemma: m 1D regularity} that provides a new result on regularity of the stochastic LLB equation in dimensions one. This theorem allows us to approximate the stochastic LLB equation~\eqref{eq: sLLB} directly, without using regularisation. We note that a similar problem in dimension one is considered in \cite{dunst} for the stochastic Landau-Lifshitz-Gilbert equation. In the Appendix, we quote the existence, uniqueness and regularity results from \cite[Theorems 2.2 and 2.3]{BGL_sLLB}, that are used in proofs of the paper.

\section{Regularised equation when $ d=2 $}\label{Section: regularisation}	
\subsection{Uniform estimates}\label{Section: uniform estimates}
	We consider the Neumann Laplacian $\Delta$ with the domain 
	\[D(\Delta)=\left\{\varphi\in\b H^2;\, \frac{\partial\varphi}{\partial \nu}(x)=0,\,\,x\in\partial D\right\}\,,\]
	where $\nu$ stands for the inner normal at $x\in\partial D$. 
	We will also need the operator $\Delta^2$ with the domain $D\left(\Delta^2\right)=\left\{\varphi\in D(\Delta);\,\Delta\varphi\in D(\Delta)\right\}$. 
	Let $ \{e_i\} $ be an orthonormal basis of $ \L^2 $ consisting of eigenvectors of the Neumann Laplacian $\Delta $ in $ D $. 
	Let $ S_n := \text{span}\{e_1,\ldots,e_n\} $ and $ \Pi_n $ be the orthogonal projection from $ \L^2 $ onto $ S_n $ such that
	\begin{align*}
		\lb \Pi_n u, \phi \rb_{\L^2} = \lb u, \phi \rb_{\L^2}, \quad \forall \phi \in S_n. 
	\end{align*} 
	Fix $ \epsilon \in (0,1) $. 
	In order to prove the existence of solution to \eqref{eq: sLLB reg}, we use the Faedo-Galerkin approximation as in \cite{BGL_sLLB}. 
	Consider the approximate equation:
	\begin{equation}\label{eq: sLLB Galerkin}
		\begin{aligned}
			dm_n^\epsilon 
			&= \la -\epsilon \Delta^2 m_n^\epsilon + \Delta m_n^\epsilon  + \Pi_n \la m_n^\epsilon \times \Delta m_n^\epsilon \ra - \Pi_n \la (1+ |m_n^\epsilon|^2) m_n^\epsilon \ra \ra dt \\
			&\quad + \frac{1}{2}\sum_{k=1}^n \Pi_n \la \la m_n^\epsilon \times g_k \ra \times g_k \ra dt 
			+ \sum_{k=1}^n \Pi_n \la g_k + m_n^\epsilon \times g_k \ra dW_k,
		\end{aligned} 
	\end{equation}
	with $ m_n^\epsilon(0) = \Pi_n m_0 $. 
	The map $ S_n \ni v \mapsto \Delta^2 v \in S_n $ is globally Lipschitz. 
	Using \cite[Lemma 3.1]{BGL_sLLB}, we find that there exists a unique local solution of \eqref{eq: sLLB Galerkin}.

	In the following lemmas, we prove a priori estimates of $m_n^\epsilon $ in various norms. 
	The proofs are postponed to Section \ref{Section: Proofs of uniform estimates}. 
	\begin{lemma}\label{Lemma: mn^e est H1}
		Assume that $ |m_0|_{\H^1} \leq C_1 $.  
		For any $ n \in \mathbb{N} $, $ p \in [1,\infty) $ and $ T \in [0,\infty) $,
		\begin{equation*}
			\E \left[ \sup_{t \in [0,T]} |m_n^\epsilon(t)|_{\H^1}^{2p} + \la \int_0^T |\Delta m_n^\epsilon(t)|^2_{\L^2} \ dt \ra^p + \la \epsilon \int_0^T |\nabla \Delta m_n^\epsilon(t)|_{\L^2}^2 \ dt \ra^p \right]
			\leq c,
		\end{equation*}
		for some constant $ c>0 $ that may depend on $ p,T, C_1 $ and $ C_g $ but not on $m_0$, $ n $ and $ \epsilon $.   
	\end{lemma}
	
	\begin{lemma}\label{Lemma: mn^e est H2}
		Assume that $ |m_0|_{\H^2} \leq C_1 $. 
		For any $ n \in \mathbb{N} $, $ p \in [1,\infty) $ and $ T \in [0,\infty) $,
		\begin{equation*}
			\epsilon^{\frac{5}{3}p} \E \left[ \sup_{t \in [0,T]}|\Delta m_n^\epsilon(t)|^{2p}_{\L^2} 
			+ \epsilon^p \la \int_0^T |\Delta^2 m_n^\epsilon(t)|^2_{\L^2} \ dt \ra^p \right] 
			\leq c,
		\end{equation*}
		for some constant $ c>0 $ that may depend on $ p,T, C_1 $ and $ C_g $ but not on $m_0$, $ n $ and $ \epsilon $.   
	\end{lemma}
	
	\begin{lemma}\label{Lemma: mn^e est H3}  
		Assume that $ |m_0|_{\H^3} \leq C_1 $. 
		For any $ n \in \mathbb{N} $, $ p \in [1,\infty) $ and $ T \in [0,\infty) $,
		\begin{equation*}
			\epsilon^{\frac{8}{3}p} \E \left[ \sup_{t \in [0,T]}|\nabla \Delta m_n^\epsilon(t)|^{2p}_{\L^2} 
			+ \epsilon^p \la \int_0^T |\nabla \Delta^2 m_n^\epsilon(t)|^2_{\L^2} \ dt \ra^p \right] 
			\leq c,
		\end{equation*}
		for some constant $ c>0 $ that may depend on $ p,T, C_1 $ and $ C_g $ but not on $m_0$, $ n $ and $ \epsilon $.   
	\end{lemma}

\subsection{Tightness and convergence}\label{Section: tightness and convergence}
	Let $ |m_0|_{\H^3} \leq C_1 $. By Lemmas \ref{Lemma: mn^e est H1} -- \ref{Lemma: mn^e est H3} and \cite[Lemma 4.1]{BGL_sLLB}, 
	\begin{align*}
		\E \left[ |m_n^\epsilon|_{W^{\alpha,p}(0,T; \L^2) \cap L^p(0,T; \H^3) \cap L^2(0,T; \H^5)} \right] \leq c(\epsilon,p),
	\end{align*}
	for $ p \in [2,\infty) $ and $ \alpha \in (0,\frac{1}{2}) $ with $ \alpha - \frac{1}{p} < \frac{1}{2} $ and $ \alpha p > 1 $. 
	Here, the constant $ c(\epsilon,p) $ depends on $ \epsilon $ and $ p $ but not on $ n $. 
	Let $ p=4 $ and 
	\begin{align*}
		B &:= W^{\alpha,4}(0,T; \L^2) \cap L^4(0,T; \H^3) \cap L^2(0,T;\H^5), \\
		E &:= \mathcal{C}([0,T]; \H^{-1}) \cap L^4(0,T; \H^2) \cap L^2(0,T; \H^4).
	\end{align*}
	As in the proof of \cite[Lemma 4.2]{BGL_sLLB}, we deduce from the compact embedding $ B \hookrightarrow E $ that the laws $ \{ \mathcal{L}(m_n^\epsilon) \}_{n \in {\b N}} $ on $ E $ are tight. 
	The following proposition is a result of Skorokhod theorem. 
	\begin{proposition}\label{Prop: Skorohod}
		There exists a probability space $ (\Omega', \mathcal{F}', \P') $, a sequence of random variables $ \{({m_n^\epsilon}',W_n')\} $ 
		and a random variable $ \{{m^\epsilon}', W'\} $ on $ (\Omega',\mathcal{F}', \P') $ taking values in $ E \times \mathcal{C}([0,T]; \H^3) $, such that 
		\begin{equation}\label{eq: Skorohold same law + conv}
			\begin{aligned}
				&\mathcal{L}(m_n^\epsilon, W) = \mathcal{L}({m_n^\epsilon}', W_n') \quad \text{on } E \times \mathcal{C}([0,T]; \H^3), \ n \in \mathbb{N}, \\
				&({m_n^\epsilon}', W_n') \to ({m^\epsilon}', W') \quad \text{strongly on } E \times \mathcal{C}([0,T]; \H^3), \quad \P'\text{-a.s.}
			\end{aligned}
		\end{equation}
	\end{proposition}

	Since $ {m_n^\epsilon}' $ has the same laws as $ m_n^\epsilon $ on $ E $, and $ \mathcal{C}([0,T]; S_n) $ is continuously embedded in $ E $ (by Kuratowski's theorem), we have
	\begin{align*}
		&\sup_{n \in \mathbb{N}} 
		\E' \left[ \sup_{t \in [0,T]} |{m_n^\epsilon}'(t)|_{\H^1}^{2p} 
		+ |\Delta {m_n^\epsilon}'|_{L^2(0,T;\L^2)}^{2p} \right]
		< \infty, \\
		&\sup_{n \in \mathbb{N}}
		\E'\left[ \epsilon^{\frac{5}{3}p} \sup_{t \in [0,T]} |\Delta {m_n^\epsilon}'(t)|_{\L^2}^{2p} 
		+ \epsilon^p |\nabla \Delta {m_n^\epsilon}'|_{L^2(0,T;\L^2)}^{2p} 
		\right]
		< \infty, \\
		&\sup_{n \in \mathbb{N}}
		\E'\left[ \epsilon^{\frac{8}{3}p} \sup_{t \in [0,T]} |\nabla \Delta {m_n^\epsilon}'(t)|_{\L^2}^{2p} 
		+ \epsilon^{\frac{8}{3}p} |\Delta^2 {m_n^\epsilon}'|_{L^2(0,T;\L^2)}^{2p} 
		+ \epsilon^{\frac{11}{3}p} |\nabla \Delta^2 {m_n^\epsilon}'|_{L^2(0,T;\L^2)}^{2p} \right]
		< \infty, 
	\end{align*}
	for $ p \in [1,\infty) $.
	This implies that for fixed $ \epsilon $,
	\begin{align*}
		&\sup_{n \in \mathbb{N}} 
		\E'\left[ \la \int_0^T |(1+|{m_n^\epsilon}'(s)|^2) {m_n^\epsilon}'(s)|_{\L^2}^2 \ ds \ra^p \right] < \infty, \\
		&\sup_{n \in \mathbb{N}}
		\E'\left[ \la \int_0^T |{m_n^\epsilon}'(s) \times \Delta {m_n^\epsilon}'(s)|_{\L^2}^2 \ ds \ra^p \right] < \infty \\
		&\sup_{n \in \mathbb{N}}
		\E'\left[ \la \int_0^T |\Delta^2 {m_n^\epsilon}'(s)|_{\L^2}^2 \ ds \ra^p \right] < \infty, 
	\end{align*}
	and we can deduce the weak convergence $ {m_n^\epsilon}' \rightharpoonup {m^\epsilon}' $ in $ L^2(\Omega; L^\infty(0,T;\H^3) \cap L^2(0,T; \H^5)) $ as $ n \to \infty $ up to subsequences. The details are similar to \cite[Section 5]{BGL_sLLB} and omitted here. 
	We arrive at the following existence, uniqueness and regularity result for the regularised equation \eqref{eq: sLLB reg}, where a martingale solution is defined in the sense of \cite[Definition 2.1]{BGL_sLLB}. 
	\begin{theorem}\label{Theorem: m^e}
		Assume that $ |m_0|_{\H^3} \leq C_1 $. For every $ \epsilon \in (0,1) $, there exists a pathwise unique martingale solution $ (\Omega,\mathcal{F}, \mathbb{F}, \P, W, m^\epsilon) $ of the regularised equation \eqref{eq: sLLB reg} such that 
		\begin{equation}\label{eq: m^e Linf-t est}
			\E\left[ \sup_{t \in [0,T]} \la |m^\epsilon(t)|_{\H^1}^{2p} + \epsilon^{\frac{5}{3}p} |\Delta m^\epsilon(t)|_{\L^2}^{2p} + \epsilon^{\frac{8}{3}p} |\nabla \Delta m^\epsilon(t)|_{\L^2}^{2p} \ra \right]
			\leq c, 
		\end{equation}
		\begin{equation}\label{eq: m^e L2-t est}
			\E \left[ |\Delta m^\epsilon|_{L^2(0,T;\L^2)}^{2p} 
			+ \epsilon^p |\nabla \Delta m^\epsilon|_{L^2(0,T;\L^2)}^{2p} 
			+ \epsilon^{\frac{8}{3}p} |\Delta^2 m^\epsilon|_{L^2(0,T;\L^2)}^{2p} 
			+ \epsilon^{\frac{11}{3}p} |\nabla \Delta^2 m^\epsilon|_{\L^2(0,T;\L^2)}^{2p} \right]
			\leq c, 
		\end{equation}
		and for $ \alpha \in (0,\frac{1}{2}) $, 
		\begin{equation}\label{eq: m^e Calpha-H1 est}
			\epsilon^{\frac{4}{3}p}\,\E \left[ |m^\epsilon|_{\mathcal{C}^\alpha([0,T];\H^1)}^p \right] \leq c , 
		\end{equation}
		for $ p \in [1,\infty) $ and some constant $ c $ independent of $ \epsilon $. 
	\end{theorem}
	\begin{proof}
		The estimates \eqref{eq: m^e Linf-t est} and \eqref{eq: m^e L2-t est} follow directly from the weak convergences and Lemma \ref{Lemma: mn^e est H1} -- \ref{Lemma: mn^e est H3}. 
		We focus on \eqref{eq: m^e Calpha-H1 est}. 
		Let $ \beta \in (0,\frac{1}{2}) $, $ p \in [2,\infty) $, $ \beta - \frac{1}{p} < \frac{1}{2} $ and $ \alpha \in (0, \beta -\frac{1}{p}) $. 
		The following embeddings are continuous:
		\begin{align*}
			W^{1,2}(0,T;\H^1) \hookrightarrow W^{\beta,p}(0,T;\H^1) \hookrightarrow \mathcal{C}^\alpha([0,T];\H^1).
		\end{align*}
		The integral form of the regularised equation \eqref{eq: sLLB reg} reads
		\begin{align*}
			m^\epsilon(t) 
			&= m^\epsilon(0) + \int_0^t \la -\epsilon \Delta^2 m^\epsilon + \Delta m^\epsilon  + m^\epsilon \times \Delta m^\epsilon - (1+ |m^\epsilon|^2) m^\epsilon \ra dt \\
			&\quad + \int_0^t \sum_{k=1}^\infty \la m^\epsilon \times g_k \ra \times g_k \ dt + \int_0^t \sum_{k=1}^\infty \la g_k + m^\epsilon \times g_k \ra dW_k,
		\end{align*}
		Recall that $ \E |m^\epsilon|_{L^2(0,T;\H^1)}^{2p} \leq c $.   
		For a $ W^{1,2}(0,T;\H^1) $-estimate of the drift part, 
		\begin{align*}
			\E \left[ |\epsilon \Delta^2 m^\epsilon|_{L^2(0,T;\H^1)}^{2p} \right]
			&\leq c\E \left[ |\epsilon \Delta^2 m^\epsilon|_{L^2(0,T;\L^2)}^{2p} \right] + c \E \left[ |\epsilon \nabla \Delta^2 m^\epsilon|_{L^2(0,T;\L^2)}^{2p} \right] \\
			&\leq c \la \epsilon^{-\frac{2}{3}p} + \epsilon^{-\frac{5}{3}p} \ra, \\
			\E \left[ |\Delta m^\epsilon|_{L^2(0,T;\H^1)}^{2p} \right]
			&\leq c\E \left[ |\Delta m^\epsilon|_{L^2(0,T;\L^2)}^{2p} \right] + c \E \left[ |\nabla \Delta m^\epsilon|_{L^2(0,T;\L^2)}^{2p} \right] \\
			&\leq c \la 1+\epsilon^{-p} \ra, \\
			\E \left[ |m^\epsilon \times \Delta m^\epsilon|_{L^2(0,T;\H^1)}^{2p} \right]
			&\leq 
			c\E \left[ |m^\epsilon \times \Delta m^\epsilon|_{L^2(0,T;\L^2)}^{2p} +  |\nabla m^\epsilon \times \Delta m^\epsilon|_{L^2(0,T;\L^2)}^{2p} \right] \\
			&\quad + c \E \left[ |m^\epsilon \times \nabla \Delta m^\epsilon|_{L^2(0,T;\L^2)}^{2p} \right] \\
			&\leq 
			c \E \left[ |m^\epsilon|_{L^\infty(0,T;\L^4)}^{2p} |\Delta m^\epsilon|_{L^2(0,T;\L^4)}^{2p} + |\nabla m^\epsilon|_{L^\infty(0,T;\L^4)}^{2p} |\Delta m^\epsilon|_{L^2(0,T;\L^4)}^{2p} \right] \\
			&\quad + c\E \left[ |m^\epsilon|_{L^\infty(0,T;\L^\infty)}^{2p} |\nabla \Delta m^\epsilon|_{L^2(0,T;\L^2)}^{2p} \right] \\
			&\leq 
			c \E \left[ |m^\epsilon|_{L^\infty(0,T;\H^1)}^{4p} \right]^\frac{1}{2} \E \left[ |\Delta m^\epsilon|_{L^2(0,T;\H^1)}^{4p} \right]^\frac{1}{2} \\
			&\quad + c \E \left[ |\nabla m^\epsilon|_{L^\infty(0,T;\H^1)}^{4p} \right]^\frac{1}{2} \E \left[|\Delta m^\epsilon|_{L^2(0,T;\H^1)}^{4p} \right]^\frac{1}{2} \\
			&\quad + c \E \left[ |m^\epsilon|_{L^\infty(0,T;\H^2)}^{4p} \right]^\frac{1}{2} \E \left[ |\nabla \Delta m^\epsilon|_{L^2(0,T;\L^2)}^{4p} \right]^\frac{1}{2} \\
			&\leq 
			c \la \epsilon^{-p} + \epsilon^{-\frac{5}{3}p} \epsilon^{-p}\ra, \\
			\E \left[ | (1+ |m^\epsilon|^2) m^\epsilon|_{L^2(0,T;\H^1)}^{2p}  \right]
			&\leq c \E \left[ |m^\epsilon|_{L^2(0,T;\H^1)}^{2p} + |m^\epsilon|_{L^\infty(0,T;\L^6)}^{6p} + ||m^\epsilon|^2 \nabla m^\epsilon|_{L^2(0,T;\L^2)}^{2p} \right] \\
			&\leq c \la 1 + \E \left[|m^\epsilon|_{L^\infty(0,T;\H^1)}^{6p} + |m^\epsilon|_{L^\infty(0,T;\L^8)}^{4p} |\nabla m^\epsilon|_{L^2(0,T;\L^4)}^{2p} \right] \ra \\
			&\leq c \la 1 + \E \left[|m^\epsilon|_{L^\infty(0,T;\H^1)}^{6p} + |m^\epsilon|_{L^\infty(0,T;\H^1)}^{8p} \right] + \E \left[|\nabla m^\epsilon|_{L^2(0,T;\H^1)}^{2p} \right] \ra \\
			&\leq c, 
		\end{align*}
		Similarly, 
		\begin{align*}
			&\E \left[ \la \int_0^T \sum_{k=1}^\infty |(m^\epsilon \times g_k) \times g_k|_{\H^1} \ dt \ra^{2p} \right] \\
			&\leq c \E \left[ \la \int_0^T \sum_{k=1}^\infty |(m^\epsilon \times g_k) \times g_k|_{\L^2} \ dt \ra^{2p} \right] 
			+ c \E \left[ \la \int_0^T \sum_{k=1}^\infty |(\nabla m^\epsilon \times g_k) \times g_k|_{\L^2} \ dt \ra^{2p} \right] \\
			&\quad + c \E \left[ \la \int_0^T \sum_{k=1}^\infty \la \int_D |m^\epsilon|^2 |\nabla g_k|^2 |g_k|^2 \ dx \ra^\frac{1}{2} \ dt  \ra^{2p} \right] \\
			&\leq c \E \left[ |m^\epsilon|_{L^\infty(0,T;\L^2)}^{2p} + |\nabla m^\epsilon|_{L^\infty(0,T;\L^2)}^{2p} \right] \la \sum_{k=1}^\infty |g_k|_{\L^\infty}^2 \ra^{2p} 
			+ c \E \left[ |m^\epsilon|_{L^\infty(0,T;\L^2)}^{2p} \right] \la \sum_{k=1}^\infty |\nabla g_k|_{\L^\infty} |g_k|_{\L^\infty} \ra^{2p} \\
			&\leq c.
		\end{align*}	
		For an $ W^{\beta,p}(0,T;\H^1) $-estimate of the diffusion part,
		\begin{align*}
			&\E \left[ \left| \sum_{k=1}^\infty \int_s^t \la g_k + m^\epsilon \times g_k \ra \ dW_k \right|^{p}_{W^{\beta,p}(0,T;\H^1)} \right] \\
			&\leq c \E \left[ \int_0^T \la \sum_{k=1}^\infty |g_k + m^\epsilon \times g_k|_{\H^1}^2 \ra^{\frac{p}{2}} dt \right] \\
			&\leq \E \left[ \la \int_0^T \sum_{k=1}^\infty \la |g_k|_{\H^1}^2 + | m^\epsilon|_{\H^1}^2 |g_k|_{\L^\infty}^2 + |m^\epsilon|_{\L^2}^2 |\nabla g_k|_{\L^2}^2 \ra dt \ra^{\frac{p}{2}} \right] \\
			&\leq c \E \left[ | m^\epsilon|_{L^\infty(0,T;\H^1)}^{p} + 1 \right] \la \sum_{k=1}^\infty \la |g_k|_{\L^\infty}^2 + |g_k|_{\H^1}^2 \ra \ra^{\frac{p}{2}} \\
			&\leq c.
		\end{align*}	
		Hence, 
		\begin{align*}
			\E \left[ |m^\epsilon|_{\mathcal{C}^\alpha([0,T]; \H^1)}^{p} \right]
			&\leq c \E \left[ |m^\epsilon|_{W^{\beta,p}(0,T;\H^1)}^{p} \right]
			\leq c \epsilon^{-\frac{4}{3}p}.
		\end{align*}
		The case $ p=1 $ is follows immediately from Jensen's inequality. 
	\end{proof}

\subsection{$ \epsilon $-Convergence}
	\begin{theorem}\label{Theorem: epsi conver}
	For any $\epsilon\in(0,1)$, let $m^\epsilon$ and $m$ be the solution of the regularised equation~\eqref{eq: sLLB reg} and the stochastic LLB equation~\eqref{eq: sLLB}, respectively. Then
	\begin{align}\label{eq: u epsi conver}
		\lim_{\epsilon \to 0} \E \left[ \sup_{t \in [0,T]} |m^\epsilon(t)-m(t)|^2_{\L^2} + \frac{1}{2}\int_0^T |\nabla (m^\epsilon-m)(s)|^2_{\L^2} \ ds \right] = 0. 
	\end{align}
	Furthermore, for every $ \gamma >0 $, there exists a constant $c=c(\gamma)$ independent of $\epsilon$ such that for any $ q,\beta \in (0,1) $,
	\begin{equation}\label{eq: P(u^e)}
		\P \la \sup_{t \in [0,T]} |m^\epsilon(t)-m(t)|_{\L^2}^2 + \frac{1}{2} \int_0^T |\nabla (m^\epsilon-m)(s)|_{\L^2}^2 \ ds > \gamma \epsilon^{1-\beta} \ra
		\leq c\epsilon^{\beta(1-q)} + c K(\epsilon)^{-1},
	\end{equation}
	where $ K(\epsilon) := - \frac{1}{2}\beta q (\ln \epsilon) $.
	\end{theorem}
	\begin{proof}	
		Let $ u^\epsilon = m^\epsilon-m $. 
		We deduce the equation of $u^\epsilon$ :
		\begin{align*}
			du^\epsilon 
			&= -\epsilon \Delta^2 m^\epsilon dt 
			+ \la \Delta u^\epsilon + u^\epsilon \times \Delta m^\epsilon + m \times \Delta u^\epsilon - u^\epsilon - \lb m+m^\epsilon, u^\epsilon \rb m^\epsilon - |m|^2 u^\epsilon \ra dt \\
			&\quad + \frac{1}{2} \sum_{k=1}^\infty (u^\epsilon \times g_k) \times g_k \ dt 
			+ \sum_{k=1}^\infty u^\epsilon \times g_k \ dW_k,
		\end{align*}
		with initial condition $ u^\epsilon(0) =0 $. 
		Applying It{\^o}'s lemma, we have
		\begin{equation}\label{eq: |u|^2}
			\begin{aligned}
				&\frac{1}{2} |u^\epsilon(t)|^2_{\L^2} + \int_0^t |u^\epsilon|^2_{\H^1} \ ds + \int_0^t \int_D |m|^2 |u^\epsilon|^2 \ dx \ ds + \frac{1}{2}\int_0^t \sum_{k=1}^\infty |u^\epsilon \times g_k|^2_{\L^2} \ ds \\
				&= \int_0^t \la -\lb \epsilon \Delta^2 m^\epsilon, u^\epsilon \rb_{\L^2}(s) 
				+ \lb m \times \Delta u^\epsilon, u^\epsilon \rb_{\L^2}(s) 
				- \lb \langle m+m^\epsilon, u^\epsilon \rangle m^\epsilon, u^\epsilon \rb_{\L^2}(s) \ra ds \\ 
				&= \int_0^t \la V_1(s) + V_2(s) + V_3(s) \ra ds. 
			\end{aligned}
		\end{equation}
		We first estimate each $ V_i $ for $ i = 1, \ldots, 3 $. Let $ \delta \in (0,1) $. 
		\begin{align*}
			V_1 
			&= -\epsilon \lb \Delta^2 m^\epsilon, u^\epsilon \rb_{\L^2} \\
			&= -\epsilon \lb \Delta m^\epsilon, \Delta u^\epsilon \rb_{\L^2} \\
			&= -\epsilon |\Delta u^\epsilon|_{\L^2}^2 - \epsilon \lb \Delta m, \Delta u^\epsilon \rb_{\L^2} \\
			&\leq -\frac{1}{2}\epsilon |\Delta u^\epsilon|_{\L^2}^2 + \frac{1}{2} \epsilon |\Delta m|_{\L^2}^2, \\
			V_2 
			&= \lb \nabla u^\epsilon, \nabla m \times u^\epsilon \rb_{\L^2} \\
			&\leq |\nabla u^\epsilon|_{\L^2} |\nabla m|_{\L^4} |u^\epsilon|_{\L^4} \\			
			&\leq c |\nabla u^\epsilon|_{\L^2} \la |u^\epsilon|_{\L^2}^\frac{1}{2} |\nabla u^\epsilon|_{\L^2}^\frac{1}{2} + |u^\epsilon|_{\L^2} \ra |\nabla m|_{\L^4} \\
			&\leq c |\nabla u^\epsilon|_{\L^2}^\frac{3}{2} |u^\epsilon|_{\L^2}^\frac{1}{2} |\nabla m|_{\L^4} 
			+ c |\nabla u^\epsilon|_{\L^2} |u^\epsilon|_{\L^2} |\nabla m|_{\L^4} \\
			&\leq \frac{3}{4}\delta^\frac{4}{3} |\nabla u^\epsilon|_{\L^2}^2 + c \delta^{-4} |u^\epsilon|_{\L^2}^2 |\nabla m|_{\L^4}^4 
			+ \frac{1}{2}\delta |\nabla u^\epsilon|_{\L^2}^2 + c \delta^{-1} |u^\epsilon|_{\L^2}^2 |\nabla m|_{\H^1}^2 \\
			&\leq \frac{5}{4} \delta |\nabla u^\epsilon|_{\L^2}^2 + c \delta^{-4} |u^\epsilon|_{\L^2}^2 \la |\nabla m|_{\L^2}^2 |\Delta m|_{\L^2}^2 + |\nabla m|_{\L^2}^4 + |\nabla m|_{\L^2}^2 + |\Delta m|_{\L^2}^2 \ra, \\
			V_3 
			&= - \int_D \lb m^\epsilon, u^\epsilon \rb^2 \ dx - \int_D  \langle m, u^\epsilon \rangle \langle m^\epsilon, u^\epsilon \rangle \ dx \\
			&\leq - \frac{1}{2}\int_D \lb m^\epsilon, u^\epsilon \rb^2 \ dx + \frac{1}{2} \int_D |m|^2 |u^\epsilon|^2 \ dx
		\end{align*}
		Taking $ \delta = \frac{1}{5} $, we have
		\begin{align*}
			&\sup_{t \in [0,T]} |u^\epsilon(t)|^2_{\L^2} + \frac{1}{2}\int_0^T |\nabla u^\epsilon(s)|^2_{\L^2} \ ds + \int_0^T \int_D |m(s)|^2 |u^\epsilon(s)|^2 \ dx \ dx \\
			&\leq \epsilon \int_0^T |\Delta m(s)|^2_{\L^2} \ ds + \int_0^T \phi(s) \sup_{r \leq s}|u^\epsilon(r)|^2_{\L^2} \ ds,
		\end{align*}
		where
		\begin{align*}
			\phi := c \delta^{-4} \la |\nabla m|_{\L^2}^2 |\Delta m|_{\L^2}^2 + |\nabla m|_{\L^2}^4 + |\nabla m|_{\L^2}^2 + |\Delta m|_{\L^2}^2 \ra.
		\end{align*}
		By Gronwall's Lemma, 
		\begin{align*}
			\sup_{t \in [0,T]} |u^\epsilon(t)|^2_{\L^2}
			&\leq \epsilon \la \int_0^T |\Delta m (s)|^2_{\L^2} \ ds \ra e^{\int_0^T \phi(s) \ ds}, \\
			\frac{1}{2}\int_0^T |\nabla u^\epsilon(s)|^2_{\L^2} \ ds
			&\leq \epsilon \la \int_0^T |\Delta m (s)|^2_{\L^2} \ ds \ra \la 1+ e^{\int_0^T \phi(s) \ ds} \int_0^T \phi(s) \ ds \ra.
		\end{align*}
		Since there exists a constant $ \tilde{c}>0 $ independent of $ \epsilon $ such that
		\begin{align*}
			&\E \left[ |\Delta m|^2_{L^2(0,T;\L^2)} \right] 
			+ \E \left[ \int_0^T \phi(s) \ ds \right] \\
			&\leq c \delta^{-4} \E \left[ (1+T)|\nabla m|^4_{L^\infty(0,T;\L^2)} + |\Delta m|^4_{L^2(0,T;\L^2)} + |\nabla m|^2_{L^2(0,T;\L^2)} + |\Delta m|^2_{L^2(0,T;\L^2)} \right]
			\leq \tilde{c},
		\end{align*}
		the random variables $ \int_0^T |\Delta m|^2_{\L^2} ds $, $ \int_0^T \phi(s) ds $ and thus $ e^{\int_0^T \phi(s) ds} $ are finite, $ \P $-a.s. 
		Hence, we have the almost-sure convergence:
		\begin{align*}
			\lim_{\epsilon \to 0} \la \sup_{t \in [0,T]} |u^\epsilon(t)|^2_{\L^2} + \frac{1}{2}\int_0^T |\nabla u^\epsilon(s)|^2_{\L^2} \ ds \ra = 0, \quad \P \text{-a.s.}
		\end{align*}
		Given the uniform integrability of the family $\left\{u^\epsilon;\,\epsilon\in(0,1)\right\} $ in $ L^2(\Omega; L^\infty(0,T;\L^2) \cap L^2(0,T;\H^1)) $, 
		the convergence~\eqref{eq: u epsi conver} can be deduced by the Vitali convergence theorem.
		
		To prove~\eqref{eq: P(u^e)}, we define a stopping time 
		\begin{align*}
			\tau_K := \inf \ \left\{ t \in [0,\infty) : \int_0^t \phi(s) \ ds \geq K \right\} \wedge T. 	
		\end{align*}
		With the uniform estimates of $ m $, we have from Markov's inequality that
		\begin{align*}
			\P \la \int_0^T \phi(s) \ ds > K \ra
			\leq \tilde{c} K^{-1},
		\end{align*}
		where $ \tau_K \leq T $. 
		We now can estimate the error as follows:
		\begin{align*}
			&\P \la \sup_{t \in [0,T]} |u^\epsilon(t)|_{\L^2}^2 + \frac{1}{2} \int_0^T |\nabla u^\epsilon(s)|_{\L^2}^2 \ ds > \gamma \epsilon^{1-\beta} \ra \\
			&\leq \P \la \sup_{t \in [0,\tau_K]} |u^\epsilon(t)|_{\L^2}^2 + \frac{1}{2} \int_0^{\tau_K} |\nabla u^\epsilon(t)|_{\L^2}^2 \ dt > \gamma \epsilon^{1-\beta} \ra 
			+ \P(\tau_K < T) \\
			&\leq \gamma^{-1} \epsilon^{\beta} \E \left[ \la 1+ e^{\int_0^{\tau_K} \phi(s) \ ds}\la 1+\int_0^{\tau_K} \phi(s) \ ds \ra \ra \int_0^{\tau_K} |\Delta m(s)|_{\L^2}^2 \ ds \right] + \tilde{c} K^{-1} \\
			&\leq \gamma^{-1} \epsilon^{\beta} \la 1+e^{K}(1+K) \ra \E \left[ \int_0^T |\Delta m(s)|_{\L^2}^2 \ ds \right] + \tilde{c} K^{-1} \\
			&\leq c \la \epsilon^{\beta} e^{2K} + K^{-1} \ra, 
		\end{align*}
		where $ c $ may depend on $ \gamma $ and $ \tilde{c} $ but not on $ K $ and $ \epsilon $. 
		Clearly, $ K = K(\epsilon) \to \infty $ as $ \epsilon \to 0 $. Then 
		\begin{equation*}
			\P \la \sup_{t \in [0,T]} |u^\epsilon(t)|_{\L^2}^2 + \frac{1}{2} \int_0^T |\nabla u^\epsilon(s)|_{\L^2}^2 \ ds > \epsilon^{1-\beta} \ra
			\leq c \la \epsilon^{\beta(1-q)} + K(\epsilon)^{-1} \ra
			\to 0, 
		\end{equation*}
		as $ \epsilon \to 0 $.
	\end{proof}

\section{Finite element approximations of $ m^\epsilon $}\label{Section: FEM m-epsilon}
	In this section, we first define an increasing sequence of $ \mathbb{F} $-stopping times $\{\tau^{R,\epsilon}\}_{R\in \mathbb{N}}$ for $ m^\epsilon $ and then propose a numerical scheme to approximate the regularised equation up to the stopping time $\tau^{R,\epsilon}$ for any $R\in \mathbb{N}$. 

\subsection{Stopping time for $ m^\epsilon $}
	For $ R \geq 1 $, let $ (\tau^{R,\epsilon})_{R \in \mathbb{N}} $ be an increasing sequence of $ \mathbb{F} $-stopping times:
	\begin{align*}
		\tau^{R,\epsilon} := \inf \ \{ t \in [0,\infty): |m^\epsilon(t)|_{\H^1} \geq R \} \wedge T. 
	\end{align*}
	By Theorem \ref{Theorem: m^e} and Markov's inequality, 
	\begin{equation}\label{eq: P(H1exit)}
		\P\la \sup_{t \in [0,T]} |m^\epsilon(t)|_{\H^1} > R \ra
		\leq R^{-1} \E \left[ \sup_{t \in [0,T]} |m^\epsilon(t)|_{\H^1} \right]
		\leq c R^{-1}. 
	\end{equation}
	Hence, $ \tau^{R,\epsilon} \to T $ as $ R \to \infty $, $ \P $-a.s. for all $ \epsilon \in (0,1) $. 

\subsection{Approximate solution $ m_h^{(n,\epsilon,R)} $}
	For space discretisation, let $ \mathcal{T}_h $ be a regular triangulation of $ D \subset \R^2 $ into triangles of maximal mesh size $ h $. 
	Let $ \V_h $ denote the space of all continuous piecewise quadratic functions on $ \mathcal{T}_h $. 
	
	For time discretisation, let $ N \in \mathbb{N} $. 
	For every $ n \geq 1 $, we define
	\begin{align*}
		\Delta t := \frac{T}{N}, \quad t_n := n \Delta t, \quad \Delta^{(n)} W := W(t_n) - W(t_{n-1}). 
	\end{align*}
	Similarly, we define the discrete stopping time and increments. 
	The construction follows from \cite{BreitProhl}.
	\begin{align*}
		\tau_n^{\epsilon,R} := \max_{1 \leq l \leq n} \{ t_l : t_l \leq \tau^{\epsilon,R} \}, \quad 
		\Delta \tau_n^{\epsilon,R} := \tau_n^{\epsilon,R} - \tau_{n-1}^{\epsilon,R}, \quad
		\Delta^{(n,\epsilon,R)} W = W(\tau_n^{\epsilon,R}) - W(\tau_{n-1}^{\epsilon,R}). 
	\end{align*}	
	Here, $ \Delta \tau_n^{\epsilon,R} \leq \Delta t $. 
	Moreover, we deduce from the following cases that $ \Delta \tau_n^{\epsilon,R} \leq \Delta \tau_{n-1}^{\epsilon,R} $:
	\begin{enumerate}[(i)]
		\ie
		\item 
		if $ \Delta \tau_n^{\epsilon,R} = 0 $, then naturally $ 0 \leq \Delta \tau_{n-1}^{\epsilon,R} $, 
		
		\item 
		if $ \Delta \tau_n^{\epsilon,R} = \Delta t $, then $ \tau_n^{\epsilon,R} = t_n $ and $ \tau_{n-1}^{\epsilon,R} = t_{n-1} $, implying $ \Delta \tau_n^{\epsilon, R} = \Delta \tau_{n-1}^{\epsilon,R} $. 
	\end{enumerate}
	The case $ 0< \Delta \tau_n^{\epsilon,R} < \Delta t $ is not possible since it implies $ t_{n-1} = \tau_{n-1}^{\epsilon,R} < \tau_n^{\epsilon,R} < t_n $ and thus $ \tau_n^{\epsilon,R} = t_{n-1} $, giving a contradiction. 
	
	We aim to find a sequence of $ \V_h $-valued random variables $ \{\smash{m_h^{(n,\epsilon,R)}}\}_n $ such that
	\begin{equation}\label{eq: tR-FEM} 
		\begin{aligned}
			&\lb m_h^{(n,\epsilon,R)} - m_h^{(n-1,\epsilon,R)}, \phi \rb_{\L^2} \\
			&= - \epsilon \lb \Delta m_h^{(n,\epsilon,R)}, \Delta \phi \rb_{\L^2} \Delta \tau_n^{\epsilon,R} \\ 
			&\quad -\lb \nabla m_h^{(n,\epsilon,R)} + m_h^{(n-1,\epsilon,R)} \times \nabla m_h^{(n,\epsilon,R)}, \nabla \phi \rb_{\L^2} \Delta \tau_n^{\epsilon,R} \\
			&\quad -\lb (1+|m_h^{(n-1,\epsilon,R)}|^2) m_h^{(n,\epsilon,R)}, \phi \rb_{\L^2} \Delta \tau_n^{\epsilon,R} \\
			&\quad + \frac{1}{2} \sum_{k=1}^\infty \lb \la m_h^{(n-1,\epsilon,R)} \times g_k \ra \times g_k, \phi \rb_{\L^2} \Delta \tau_n^{\epsilon, R} \\
			&\quad + \sum_{k=1}^\infty \lb \la g_k + m_h^{(n-1,\epsilon,R)} \times g_k \ra , \phi \rb_{\L^2} \Delta^{(n)} W_k \frac{\Delta \tau_n^{\epsilon,R}}{\Delta t}, 
			\quad \forall \phi \in \V_h, \quad \P\text{-a.s.}
		\end{aligned} 
	\end{equation}
	Fix $ \omega \in \Omega $, given $ m_h^{(n-1,\epsilon,R)}(\omega) $ and $ W^{(n)}(\omega) $, there exists a unique solution $ m_h^{(n,\epsilon,R)}(\omega) \in \V_h $ of \eqref{eq: tR-FEM} by Lax-Milgram theorem. 
	Note that the solution $ \smash{m_h^{(n,\epsilon,R)}} $ is a random variable. 
	
	For simplicity, we denote $ m_h^{(j,\epsilon,R)} $ by $ m_h^{(j)} $ for $ 0\leq j \leq N $. A priori estimates of $ m_h^{(j)} $ are proved in the following lemma. 
	\begin{lemma}\label{Lemma: mhn L2}
		Assume that $ |m_h^{(0)}|^2_{\L^2} \leq C_1 $ for all $ h $. There exists a positive constant $ c $ independent of $ \epsilon, R, h, \Delta t $ such that for $ p \in [1,\infty) $,
		\begin{align*}
			\E \left[ \max_{l \leq n} |m_h^{(l)}|^{2p}_{\L^2} \right] 
			+ \E \left[ \la \sum_{j=1}^n |m_h^{(j)} - m_h^{(j-1)}|^2_{\L^2} \ra^p 
			+ \la \sum_{j=1}^n \la |m_h^{(j)}|^2_{\H^1} + \epsilon |\Delta m_h^{(j)}|_{\L^2}^2 \ra \Delta \tau_j^{\epsilon,R} \ra^p \right] \leq c.
		\end{align*}
	\end{lemma}
	\begin{proof}
		Let $ \phi = m_h^{(n)} $ in \eqref{eq: tR-FEM}. 
		Then, 
		\begin{equation}\label{eq: 1/2|mhn|^2}
			\begin{aligned}
				&\frac{1}{2}|m_h^{(n)}|^2_{\L^2} + \frac{1}{2} |m_h^{(n)} - m_h^{(n-1)}|^2_{\L^2} - \frac{1}{2} |m_h^{(n-1)}|^2_{\L^2} \\
				&= -\la \epsilon |\Delta m_h^{(n)}|_{\L^2}^2 + |m_h^{(n)}|^2_{\L^2} +|\nabla m_h^{(n)}|^2_{\L^2} \ra \Delta \tau_n^{\epsilon,R} 
				- \int_D |m_h^{(n-1)}|^2 |m_h^{(n)}|^2 \ dx \ \Delta \tau_n^{\epsilon,R} \\
				&\quad + \frac{1}{2} \sum_{k=1}^\infty \lb m_h^{(n-1)} \times g_k, -m_h^{(n)} \times g_k \rb_{\L^2} \Delta \tau_n^{\epsilon,R} \\
				&\quad + \sum_{k=1}^\infty \lb g_k + m_h^{(n-1)} \times g_k, m_h^{(n)} \rb_{\L^2} \Delta^{(n)} W_k \frac{\Delta \tau_n^{\epsilon,R}}{\Delta t},
			\end{aligned}
		\end{equation}
		For the Stratonovich term, 
		\begin{align*}
			\frac{1}{2} \sum_{k=1}^\infty \lb m_h^{(n-1)} \times g_k, -m_h^{(n)} \times g_k \rb_{\L^2} \Delta \tau_n^{\epsilon,R} 
			&\leq \frac{1}{4} \sum_{k=1}^\infty |g_k|_{\L^\infty}^2 \la |m_h^{(n-1)}|_{\L^2}^2 + |m_h^{(n)}|_{\L^2}^2 \ra \Delta \tau_n^{\epsilon,R} \\
			&\leq \frac{1}{4} C_g \la |m_h^{(n-1)}|_{\L^2}^2 + |m_h^{(n)}|_{\L^2}^2 \ra \Delta \tau_n^{\epsilon,R}.
		\end{align*}
		For the diffusion term, the ratio $ \frac{\Delta \tau_n^{\epsilon,R}}{\Delta t} = \mathbbm{1}(t_n \leq \tau^{\epsilon,R}) $ and
		\begin{align*}
			&\sum_{k=1}^\infty \lb g_k + m_h^{(n-1)} \times g_k, m_h^{(n)} \rb_{\L^2} \Delta^{(n)} W_k \\
			&= \lb m_h^{(n)}, \Delta^{(n)} W + m_h^{(n-1)} \times \Delta^{(n)} W \rb_{\L^2} \\
			&= \lb m_h^{(n)} - m_h^{(n-1)}, \Delta^{(n)} W + m_h^{(n-1)} \times \Delta^{(n)} W \rb_{\L^2} + \lb m_h^{(n-1)}, \Delta^{(n)} W \rb_{\L^2} \\
			&\leq \frac{1}{4} |m_h^{(n)} - m_h^{(n-1)}|_{\L^2}^2 + 4 |\Delta^{(n)} W + m_h^{(n-1)} \times \Delta^{(n)} W |_{\L^2}^2 + \lb m_h^{(n-1)}, \Delta^{(n)} W \rb_{\L^2},
		\end{align*}
		where the first term on the right-hand side can be moved to the left of \eqref{eq: 1/2|mhn|^2}. 
		Then we have
		\begin{align*}
			&\E \left[ \max_{l \leq n} \la \sum_{j=1}^l|\Delta^{(j)} W + m_h^{(j-1)} \times \Delta^{(j)} W|^2_{\L^2} \mathbbm{1}(t_j \leq \tau^{\epsilon,R}) \ra^p \right] \\
			&\leq n^{p-1} \E \left[ \sum_{j=1}^n |\Delta^{(j)} W + m_h^{(j-1)} \times \Delta^{(j)} W|^{2p}_{\L^2} \ \mathbbm{1}(t_j \leq \tau^{\epsilon,R}) \right] \\
			&\leq c n^{p-1} \E \left[ \sum_{j=1}^n (1+|m_h^{(j-1)}|^{2p}_{\L^2}) \ \mathbbm{1}(t_j \leq \tau^{\epsilon,R}) (\Delta t)^p \right] \\
			&\leq c T^{p-1} \E \left[ \sum_{j=1}^n (1+|m_h^{(j-1)}|^{2p}_{\L^2}) \ \Delta \tau_j^{\epsilon,R} \right],
		\end{align*}
		and by Burkholder-Davis-Gundy inequality, 
		\begin{align*}
			\E\left[ \max_{l \leq n} \left| \sum_{j=1}^l \lb m_h^{(j-1)} \mathbbm{1}(t_j \leq \tau^{\epsilon,R}), \Delta^{(j)} W \rb_{\L^2} \right|^p \right] 
			&\leq b_p \E \left[ \la \sum_{j=1}^n |m_h^{(j-1)}|^2_{\L^2} \mathbbm{1}(t_j \leq \tau^{\epsilon,R}) \Delta t \ra^{\frac{p}{2}} \right] \\
			&\leq b_p^2 \E \left[ \la \sum_{j=1}^n |m_h^{(j-1)}|^2_{\L^2} \Delta \tau_j^{\epsilon,R} \ra^p \right] + 1 \\
			&\leq c b_p^2 \E \left[ \sum_{j=1}^n |m_h^{(j-1)}|^{2p}_{\L^2} \ \Delta \tau_j^{\epsilon,R} \right] + 1.
		\end{align*}
		Taking the $ p $th power, the sum over $ j=1,\ldots,l $, the max over $ l $ and the expectation, we deduce from \eqref{eq: 1/2|mhn|^2} and the fact $ \Delta \tau_j^{\epsilon,R} \leq \Delta t $ for all $ j $ that
		\begin{align}
			&\E \left[ \max_{l \leq n} |m_h^{(l)}|^{2p}_{\L^2} 
			+ \la \sum_{j=1}^n |m_h^{(j)} - m_h^{(j-1)}|^2_{\L^2} \ra^p 
			+ \la \sum_{j=1}^n \la \epsilon |\Delta m_h^{(j)}|_{\L^2}^2 + |m_h^{(j)}|^2_{\H^1} \ra \Delta \tau_j^{\epsilon,R} \ra^p \right] \nonumber \\
			&\leq 
			c\E \left[ |m_h^{(0)}|^{2p}_{\L^2} \right] 
			+ c\E \left[ \sum_{j=1}^n \la |m_h^{(j-1)}|_{\L^2}^{2p} + |m_h^{(j)}|_{\L^2}^{2p} \ra \Delta \tau_j^{\epsilon,R} + 1 \right]  \nonumber \\ 
			&\leq 
			c\la 1+ \E \left[ \sum_{j=1}^n \max_{l \leq j} |m_h^{(l)}|_{\L^2}^{2p} \ \Delta t \right] \ra, \label{eq: 1/2|mhn|^2 v2}
		\end{align}
		where the constant $ c $ may depend on $ p, T, C_1, C_g, b_p $ but not on $ \epsilon, R, h $ or $ \Delta t $. 
		Thus, by the (discrete) Gronwall's lemma, 
		\begin{align*}
			\E\left[ \max_{l \leq n} |m_h^{(n)}|^{2p}_{\L^2} \right] \leq c e^{cT}.
		\end{align*}
		The result follows by taking the other expectation terms in the inequality \eqref{eq: 1/2|mhn|^2 v2} as the subject. 
	\end{proof}
	
	\begin{remark}
		Assume that $ |m_h^{(0)}|_{\H^2} \leq C_1 $ for all $ h $. Then Lemma \ref{Lemma: mhn L2} holds with the sum starting from $ j=0 $ and with the term $ m_h^{(-1)} $ set to $ m_h^{(0)} $. 
	\end{remark}
	
\subsection{Error analysis $ e_h^{(n,\epsilon,R)} $}
	Let $ e_h^{(n,\epsilon,R)} := m^\epsilon(\tau_n^{\epsilon,R}) - m_h^{n,R} $. 
	Note that $ \Delta^{(n)} W = \int_{t_{n-1}}^{t_n} dW(s) $. 
	We first derive an equation for $ e_h^{(n,\epsilon,R)}$: for $ \phi \in \V_h $, 
	\begin{equation}\label{eq: <ehn,phi>}
		\begin{aligned}
			\lb e_h^{(n,\epsilon,R)}, \phi \rb_{\L^2}
			&= \lb m^\epsilon(\tau_{n-1}^{\epsilon,R})-m_h^{(n-1)}, \phi \rb_{\L^2} \\
			&\quad - \epsilon \int_{\tau_{n-1}^{\epsilon,R}}^{\tau_n^{\epsilon,R}} \lb \Delta m^\epsilon(s) - \Delta m_h^{(n)}, \Delta \phi \rb_{\L^2} \ ds \\
			&\quad - \int_{\tau_{n-1}^{\epsilon,R}}^{\tau_n^{\epsilon,R}} \lb \nabla m^\epsilon(s) - \nabla m_h^{(n)} + m^\epsilon(s) \times \nabla m^\epsilon(s) - m_h^{(n-1)} \times \nabla m_h^{(n)}, \nabla \phi \rb_{\L^2} \ ds \\
			&\quad - \int_{\tau_{n-1}^{\epsilon,R}}^{\tau_n^{\epsilon,R}} \lb (1 + |m^\epsilon(s)|^2) m^\epsilon(s) - (1+|m_h^{(n-1)}|^2) m_h^{(n)}, \phi \rb_{\L^2} \ ds \\
			&\quad + \frac{1}{2} \sum_{k=1}^\infty \int_{\tau_{n-1}^{\epsilon,R}}^{\tau_n^{\epsilon,R}} \lb \la m^\epsilon(s) - m_h^{(n-1)} \ra \times g_k, g_k \times \phi \rb_{\L^2} \ ds \\
			&\quad + \int_{t_{n-1}}^{t_n} \lb \phi, \la m^\epsilon(s)- m_h^{(n-1)} \ra \times dW(s) \rb_{\L^2} \frac{\Delta \tau_n^{\epsilon,R}}{\Delta t}.
		\end{aligned}
	\end{equation}
	Similarly, for simplicity, we write $ e_h^{(n)} $ instead of $ e_h^{(n,\epsilon,R)} $. 
	
	We now introduce the orthogonal projection $ P_h $ from $ \L^2 $ onto $ \V_h $ and its properties before proving convergence results in Theorem~\ref{Theorem: eh conve}. 
	Let $ m_h^{(0)} = P_h m_0 $. 
	Since $ m_h^{(n)} $ takes values in $ \V_h $,
	\begin{align*}
		P_h e_h^{(n)} - e_h^{(n)} 
		&= P_h m^\epsilon(\tau_n^{\epsilon,R}) - m_h^{(n)} - \la m^\epsilon(\tau_n^{\epsilon,R}) - m_h^{(n)} \ra 
		= P_h m^\epsilon(\tau_n^{\epsilon,R}) - m^\epsilon(\tau_n^{\epsilon,R}).
	\end{align*}
	For the approximation properties of $ P_h $, let $ k=0, 1, 2 $ with $ \H^0 := \L^2 $,
	\begin{equation}\label{eq: Ph approx}
		|P_h e_h^{(n)} - e_h^{(n)}|_{\H^k} \leq c h |m^\epsilon(\tau_n^{\epsilon,R})|_{\H^{k+1}}, \quad
		|P_h e_h^{(n)}|_{\H^k} \leq |e_h^{(n)}|_{\H^k},
	\end{equation}

	\begin{theorem}\label{Theorem: eh conve}
	Assume that $|m_0|_{\H^3} \leq C_1$ and $ Nh \leq C_2 $. There exists a constant $c$ independent of $\epsilon$ and $R$ such that for any $\epsilon\in(0,1)$ and $R\geq 1$,
	\begin{equation}\label{eq: eh (h,Delta t)-est}
		\E \left[ \max_{n \leq N} |e_h^{(n)}|^2_{\L^2} + \sum_{n=1}^N |\nabla e_h^{(n)}|^2_{\L^2} \Delta \tau_n^{\epsilon,R} \right]
		\leq
		c \la h+(\Delta t)^{\frac{\alpha}{2}} \ra \epsilon^{-\frac{22}{3}} R^8 e^{c \epsilon^{-3} R^8}.
	\end{equation}
	\end{theorem}
	We postpone the proof to Section \ref{Section: Proofs of error estimates}.

\subsection{Convergence in probability}	
	We write $ c^* $, $ c^\dagger $ and $ c' $ instead of $ c $ for the constant in \eqref{eq: eh (h,Delta t)-est}, \eqref{eq: P(H1exit)} and \eqref{eq: P(u^e)}, respectively.  
	For $ q,\beta \in (0,1) $, define
	\begin{equation}\label{eq: epsilon and R in (h,Delta t)}
	\begin{aligned}
		R(h,\Delta t) &:= \la -\frac{q\beta}{2c^*} \ln (h + (\Delta t)^\frac{\alpha}{2}) \ra^\frac{1}{9}, \\
		\epsilon_1(h,\Delta t) &:= (R(h,\Delta t))^{-\frac{1}{3}}, \quad 
		\epsilon (h,\Delta t) := \max \left\{ \epsilon_1(h,\Delta t), \la (\Delta t)^{\frac{\alpha q}{2}(\beta +3)} \ra^{\frac{3}{8}} \right\}, 
	\end{aligned}
	\end{equation}
	where $ R(h,\Delta t) $, $ \epsilon_1 (h,\Delta t) $ and $ \epsilon_2 (h,\Delta t) $ converge to $ 0 $ as $ h,\Delta t \to 0 $.
	
	Using Theorems \ref{Theorem: epsi conver} and \ref{Theorem: eh conve}, we prove the convergence of the solution $ m_h^{(n)} $ of \eqref{eq: tR-FEM} to the solution $ m $ of \eqref{eq: sLLB} below. 
	\begin{theorem}\label{Theorem: conver 2D}
		For every $ \gamma > 0 $, there exists a constant $c = c(\gamma)$ independent of $ \epsilon $, $h$ and $\Delta t$ such that for any $q,\beta \in (0,1)$,
		\begin{align*}
			\P \la \max_{n \leq N} |m(t_n) - m_h^{(n)}|^2_{\L^2} + \sum_{n=1}^N |\nabla m(t_n) - \nabla m_h^{(n)}|^2_{\L^2} \Delta t > \gamma \epsilon(h,\Delta t)^{1-\beta} + \gamma (h + (\Delta t)^\frac{\alpha}{2})^{1-\beta} \ra 
			&\to 0,
		\end{align*}
		as $ h,\Delta t \to 0 $, where $ \epsilon(h,\Delta t) $ is given in \eqref{eq: epsilon and R in (h,Delta t)}. 
	\end{theorem}
	\begin{proof}
		We choose $ R = R(h, \Delta t) $ and $ \epsilon = \epsilon(h,\Delta t) \geq \epsilon_1(h,\Delta t) $ as defined in \eqref{eq: epsilon and R in (h,Delta t)}. 
		Then 
		\begin{align*}
			\epsilon^{-\frac{22}{3}} e^{c^* \epsilon^{-3} R^8} \leq e^{2c^* R^9} = \la h+(\Delta t)^{\frac{\alpha}{2}} \ra^{-q\beta}.
		\end{align*}
		By Theorem \ref{Theorem: eh conve}, 
		\begin{align*}
			&\P \la \max_{n \leq N} |m^\epsilon(\tau^{\epsilon,R}_n) - m_h^{(n)}|^2_{\L^2} + \sum_{n=1}^N |\nabla m^\epsilon(\tau^{\epsilon,R}_n) 
			- \nabla m_h^{(n)}|^2_{\L^2} \Delta \tau_n^{\epsilon,R} > \gamma (h + (\Delta t)^{\frac{\alpha}{2}})^{1-\beta} \ra \\
			&\leq 
			\E\left[\max_{n \leq N} |e_h^{(n)}|^2_{\L^2} + \sum_{n=1}^N |\nabla e_h^{(n)}|^2_{\L^2} \Delta \tau_n^{\epsilon,R} \right] \gamma^{-1} (h + (\Delta t)^{\frac{\alpha}{2}})^{\beta-1} \\
			&\leq c^* \la h+(\Delta t)^{\frac{\alpha}{2}} \ra^{1-q\beta} \gamma^{-1} (h + (\Delta t)^{\frac{\alpha}{2}})^{\beta-1} \\
			&\leq c^* \gamma^{-1} \la h+(\Delta t)^{\frac{\alpha}{2}} \ra^{\beta(1- q)}.
		\end{align*}
		Since $ t_n = \tau_n^{\epsilon,R} $ when $ \tau^{\epsilon,R} \geq T $, we have from \eqref{eq: P(H1exit)} and the inequality above that
		\begin{align*}
			&\P \la \max_{n \leq N} |m^\epsilon(t_n) - m_h^{(n)}|^2_{\L^2} 
			+ \sum_{n=1}^N |\nabla m^\epsilon(t_n) - \nabla m_h^{(n)}|^2_{\L^2} \Delta t > \gamma (h + (\Delta t)^{\frac{\alpha}{2}})^{1-\beta} \ra \\
			&\leq \P \la \max_{n \leq N} |m^\epsilon(\tau^{\epsilon,R}_n) - m_h^{(n)}|^2_{\L^2} + \sum_{n=1}^N |\nabla m^\epsilon(\tau^{\epsilon,R}_n) - \nabla m_h^{(n)}|^2_{\L^2} \Delta \tau^{\epsilon,R} > \gamma (h + (\Delta t)^{\frac{\alpha}{2}})^{1-\beta} \ra \\
			&\quad + \P(\tau^{\epsilon,R} < T) \\
			&\leq c^* \gamma^{-1}  \la h+(\Delta t)^{\frac{\alpha}{2}} \ra^{\beta(1- q)} + c^\dagger (R(h,\Delta t))^{-1},
		\end{align*}
		where the right-hand side converges to $ 0 $ as $ h, \Delta t \to 0 $. 		
		Let 
		\begin{align*}
			f(\epsilon,h,\Delta t)
			&:= \sup_{t \in [0,T]} |m(t)-m^\epsilon(t)|^2_{\L^2} + \max_{n \leq N} |m^\epsilon(t_n) - m_h^{(n)}|^2_{\L^2} \\
			&\quad + \int_0^T |\nabla m(t) - \nabla m^\epsilon(t)|^2_{\L^2} \ dt 
			+  \sum_{n=1}^N |\nabla m^\epsilon(t_n) - \nabla m_h^{(n)}|^2_{\L^2} \Delta t \\
			&\quad + \sum_{n=1}^N \int_{t_{n-1}}^{t_n} |\nabla m^\epsilon(t) - \nabla m^\epsilon(t_n)|_{\L^2}^2 \ dt. 
		\end{align*}
		We observe that
		\begin{align*}
			&\P \la \sum_{n=1}^N \int_{t_{n-1}}^{t_n} |\nabla m^\epsilon(t) - \nabla m^\epsilon(t_n)|_{\L^2}^2 \ dt > \gamma (\Delta t)^{\frac{\alpha}{2}(1-\beta)} \ra \\
			&\leq \gamma^{-1}(\Delta t)^{\frac{\alpha}{2}(\beta-1)} \E \left[ \sum_{n=1}^N \int_{t_{n-1}}^{t_n} |\nabla m^\epsilon(t) - \nabla m^\epsilon(t_n)|_{\L^2}^2 \ dt \right] \\
			&\leq \gamma^{-1}(\Delta t)^{\frac{\alpha}{2}(\beta-1)} \E \left[ |m^\epsilon|_{\mathcal{C}^\alpha([0,T];\H^1)}^2 (\Delta t)^{2\alpha} T \right] \\
			&\leq c_1 \gamma^{-1} T (\Delta t)^{\frac{\alpha}{2}(\beta+3)} \epsilon^{-\frac{8}{3}}.
		\end{align*}
		Since $ \epsilon = \epsilon(h,\Delta t) \geq ((\Delta t)^{\frac{\alpha q}{2}(\beta +3 )})^{\frac{3}{8}} $, we have
		\begin{align*}
			\P \la \sum_{n=1}^N \int_{t_{n-1}}^{t_n} |\nabla m^\epsilon(t) - \nabla m^\epsilon(t_n)|_{\L^2}^2 \ dt > \gamma (\Delta t)^{\frac{\alpha}{2}(1-\beta)} \ra 
			\leq c_1 \gamma^{-1} T (\Delta t)^{\frac{\alpha}{2}(\beta +3)(1-q)}
			\to 0, 
		\end{align*}
		as $ \Delta t \to 0 $.
		
		Combining with \eqref{eq: P(u^e)}, we obtain the convergence in probability of the approximate solutions $m_h^{(n)}$ of the regularised equation to the solution $m$ of~\eqref{eq: sLLB}: 
		\begin{align*}
			&\P \la \frac{1}{2} \max_{n \leq N} |m(t_n) - m_h^{(n)}|_{\L^2}^2 
			+ \frac{1}{2} \sum_{n=1}^N \int_{t_{n-1}}^{t_n} |\nabla m(t) - \nabla m_h^{(n)}|_{\L^2}^2 \ ds > \frac{1}{2}\gamma \epsilon^{1-\beta} + \frac{1}{2}\gamma \la h+(\Delta t)^\frac{\alpha}{2} \ra^{1-\beta} \ra \\
			&\leq \P \la f(\epsilon,h,\Delta t) > \frac{1}{2}\gamma \epsilon^{1-\beta} + \frac{1}{2}\gamma \la h+\Delta t^\frac{\alpha}{2} \ra^{1-\beta} \ra \\
			&\leq 
			\P\la \sup_{t \in [0,T]} |m(t)-m^\epsilon(t)|^2_{\L^2} + \int_0^T |\nabla m(t) - \nabla m^\epsilon(t)|^2_{\L^2} \ dt > \frac{1}{4}\gamma \epsilon^{1-\beta} \ra \\
			&\quad + \P\la \max_{n \leq N} |m^\epsilon(t_n) - m_h^{(n)}|^2_{\L^2} + \sum_{n=1}^N |\nabla m^\epsilon(t_n) - \nabla m_h^{(n)}|^2_{\L^2} \Delta t > \frac{1}{8}\gamma (h + (\Delta t)^{\frac{\alpha}{2}})^{1-\beta} \ra \\
			&\quad + \P \la \sum_{n=1}^N \int_{t_{n-1}}^{t_n} |\nabla m^\epsilon(t) - \nabla m^\epsilon(t_n)|_{\L^2}^2 \ dt > \frac{1}{8}\gamma (\Delta t)^{\frac{\alpha}{2}(1-\beta)} \ra \\
			&\leq 
			c' \la \epsilon^{\beta(1-q)} -2(\beta q)^{-1} (\ln \epsilon)^{-1} \ra 
			+ 8c^* \gamma^{-1}  \la h+(\Delta t)^{\frac{\alpha}{2}} \ra^{\beta(1- q)} + c^\dagger \la R(h,\Delta t) \ra^{-1} \\
			&\quad + 8c_1 \gamma^{-1} T (\Delta t)^{\frac{\alpha}{2}(\beta +3)(1-q)},
		\end{align*}
		where $ \epsilon $ and $ R $ are given by \eqref{eq: epsilon and R in (h,Delta t)}, and then the right-hand side converges to $ 0 $ as $ h,\Delta t \to 0 $.	
	\end{proof}

\section{Error analysis when $ d=1 $}\label{Section: 1D case}

\subsection{Additional regularity of $ m $ and $ \H^1 $-stopping time}
	\begin{theorem}\label{Lemma: m 1D regularity}
		Assume that $ |m_0|_{\H^2} \leq C_1 $. 
		Then for the pathwise solution $ m $ of the stochastic LLB equation \eqref{eq: sLLB} in Theorem \ref{Theorem: sLLB existence}, there exists a positive constant $ c $ such that
		\begin{align*}
			\E \left[ \sup_{t \in [0,T]} |\Delta m(t)|^{2p}_{\L^2} + \la \int_0^{T} |\nabla \Delta m(s)|^2_{\L^2} \ ds \ra^{2p} \right] 
			&\leq c, \quad
			\E \left[ |m|_{\mathcal{C}^\alpha([0,T];\H^1)}^{p} \right]
			\leq c, 
		\end{align*}
		for $ p \in [1,\infty) $ and $ \alpha \in (0,\frac{1}{2}) $. 
	\end{theorem}	
	We again postpone the proof of the estimates in Theorem \ref{Lemma: m 1D regularity} to Section \ref{Section: Proofs of uniform estimates}. 

	For $ R \geq 1 $, let $ (\tau^R)_{R \in \mathbb{N}} $ be an increasing sequence of $ \mathbb{F} $-stopping times with
	\begin{align*}
		\tau^R := \inf \ \{ t \in [0,\infty): |m(t)|_{\H^1} \geq R \} \wedge T. 
	\end{align*}
	By Theorem \ref{Theorem: sLLB existence} and Markov inequality, 
	\begin{equation}\label{eq: P(H1exit) 1D}
		\P\la \sup_{t \in [0,T]} |m(t)|_{\H^1} > R \ra
		\leq R^{-1} \E \left[ \sup_{t \in [0,T]} |m(t)|_{\H^1} \right]
		\leq c R^{-1} \to 0, 
		\quad 
		\text{as $ R \to \infty $}.
	\end{equation}

\subsection{Error analysis $ e_h^{(n,0,R)} $}
	Let $ \mathcal{T}_h $ be a regular triangulation of $ D \subset \R $ into sub-intervals. 
	As in Section \ref{Section: FEM m-epsilon}, we define the discrete stopping time and increments
	\begin{align*}
		\tau_n^R := \max_{1 \leq l \leq n} \{ t_l : t_l \leq \tau^R \}, \quad 
		\Delta \tau_n^R := \tau_n^R - \tau_{n-1}^R.
	\end{align*}
	For the one-dimensional problem, we study the finite-element scheme \eqref{eq: tR-FEM} with $ \epsilon = 0 $, $ \V_h $ the space of all continuous piecewise linear functions on $ \mathcal{T}_h $, and $ (m, \tau_n^R, \Delta \tau_n^R) $ in place of $ (m^\epsilon, \tau_n^{\epsilon,R}, \Delta \tau_n^R) $. 
	
	For simplicity, we denote $ m_h^{(j,0,R)} $ and $ e_h^{(j,0,R)} $, by $ m_h^{(j)} $ and $ e_h^{(j)} $, respectively, for $ 0\leq j \leq N $. 
	
	\begin{lemma}\label{Lemma: mhn L2 1D}
		Assume that $ |m_h^{(0)}|^2_{\L^2} \leq C_1 $ for all $ h $. There exists a positive constant $ c $ independent of $ R, h, \Delta t $ such that for $ p \in [1,\infty) $, 
		\begin{align*}
			\E \left[ \max_{l \leq n} |m_h^{(l)}|^{2p}_{\L^2} \right] + \E \left[ \la \sum_{j=1}^n |m_h^{(j)} - m_h^{(j-1)}|^2_{\L^2} \ra^p + \la \sum_{j=1}^n |m_h^{(j)}|^2_{\H^1} \Delta \tau_j^R \ra^p \right] \leq c.
		\end{align*}
	\end{lemma}
	\begin{proof}
		The estimates can be obtained in a similar manner to Lemma \ref{Lemma: mhn L2}.
	\end{proof}

	The order of convergence is proved in the following theorem and the proof is given in Section \ref{Section: Proofs of error estimates}.
	\begin{theorem}\label{Theorem: conver 1D}
		For every $ \gamma > 0 $, there exists a constant $ c=c(\gamma) $ independent of $h$ and $\Delta t$ such that for any $q,\beta \in (0,1)$,
		\begin{align*}
			&\P \la \max_{n \leq N} |m(t_n) - m_h^{(n)}|^2_{\L^2} + \sum_{n=1}^N |\nabla m(t_n) - \nabla m_h^{(n)}|^2_{\L^2} \Delta t > \gamma (h + (\Delta t)^\alpha)^{1-\beta} \ra \\ 
			&\leq c \gamma^{-1} \la h+(\Delta t)^{\alpha} \ra^{\beta(1- q)} + c \la -\frac{q\beta}{c} \ln (h + (\Delta t)^\alpha) \ra^{-\frac{1}{4}} .
		\end{align*}
	\end{theorem}

\section{Proof of uniform estimates}\label{Section: Proofs of uniform estimates}

\subsection{Proof of Lemma \ref{Lemma: mn^e est H1}}
	The result follows from \cite[Lemmas 3.2 and 3.3]{BGL_sLLB}. We only need to address the additional $ \epsilon \Delta^2 m_n^\epsilon $ term in the estimates. 	
	For $ |m_n^\epsilon|^2_{\L^2} $, we observe that 
	\begin{align*}
		\lb m_n^\epsilon, -\epsilon \Delta^2 m_n^\epsilon \rb_{\L^2} 
		&= -\epsilon |\Delta m_n^\epsilon|^2_{\L^2} \leq 0,
	\end{align*}
	which does not change the conclusion of \cite[Lemma 3.2]{BGL_sLLB}: 
	\begin{equation}\label{eq: |mn^e| L2 est}
		\E \left[ \sup_{t \in [0,T]} |m_n^\epsilon(t)|^{2p}_{\L^2} + \la \int_0^T |\nabla m_n^\epsilon(t)|^2_{\L^2} \ dt \ra^p + \la \int_0^T  \int_D (1+|m_n^\epsilon(t)|^2) |m_n^\epsilon(t)|^2 \ dx \ dt \ra^p \right] \leq c,
	\end{equation}
	for $ p \geq 1 $ and $ T \in [0,\infty) $. 
	Similarly, for $ |\nabla m_n^\epsilon|^2_{\L^2} $, we observe that
	\begin{equation}\label{eq: epsilon part in |D|^2}
		\lb -\Delta m_n^\epsilon, -\epsilon \Delta^2 m_n^\epsilon \rb_{\L^2} 
		= -\epsilon |\nabla \Delta m_n^\epsilon|_{\L^2}^2.
	\end{equation}
	Then we can gather from \cite[Lemma 3.3]{BGL_sLLB}:
	\begin{equation}\label{eq: |Dmn^e| L2 est}
		\E \left[ \sup_{t \in [0,T]} |\nabla m_n^\epsilon(t)|^{2p}_{\L^2} + \la \int_0^T |\Delta m_n^\epsilon(t)|^2_{\L^2} \ dt \ra^p + \la \epsilon \int_0^T |\nabla \Delta m_n^\epsilon(t)|_{\L^2}^2 \ dt \ra^p \right] \leq c,
	\end{equation}
	for $ p \geq 1 $ and $ T \in [0,\infty) $, as desired.

\subsection{Proof of Lemma \ref{Lemma: mn^e est H2}}
	Applying It{\^o}'s lemma to $ \frac{1}{2} |\Delta m_n^\epsilon|^2_{\L^2} $:
	\begin{equation}\label{pf: 1/2|Delta m_n^e|^2}
		\begin{aligned}
			\frac{1}{2} |\Delta m_n^\epsilon(t)|^2_{\L^2}
			&= 
			\frac{1}{2} |\Delta m_n^\epsilon(0)|^2_{\L^2}
			-\epsilon \int_0^t |\Delta^2 m_n^\epsilon(s)|^2_{\L^2} \ ds
			-\int_0^t |\nabla \Delta m_n^\epsilon(s)|^2_{\L^2} \ ds \\
			&\quad + \int_0^t \lb \Delta^2 m_n^\epsilon(s), m_n^\epsilon(s) \times \Delta m_n^\epsilon(s) \rb_{\L^2} \ ds \\
			&\quad - \int_0^t \lb \Delta^2 m_n^\epsilon(s), (1+ |m_n^\epsilon(s)|^2) m_n^\epsilon(s) \rb_{\L^2} \ ds \\
			&\quad + \frac{1}{2} \sum_{k=1}^n \int_0^t \lb \Delta^2 m_n^\epsilon(s), \la m_n^\epsilon(s) \times g_k \ra \times g_k \rb_{\L^2} \ ds \\
			&\quad + \frac{1}{2} \sum_{k=1}^n \int_0^t |\Delta (g_k + m_n^\epsilon(s) \times g_k) |^2_{\L^2} \ ds \\
			&\quad + \sum_{k=1}^n \int_0^t \lb \Delta^2 m_n^\epsilon(s), m_n^\epsilon(s) \times g_k + g_k \rb_{\L^2} \ dW_k(s) \\
			&= \frac{1}{2} |\Delta m_n^\epsilon(0)|^2_{\L^2} - \epsilon \int_0^t |\Delta^2 m_n^\epsilon(s)|^2_{\L^2} \ ds - \int_0^t |\nabla \Delta m_n^\epsilon(s)|^2_{\L^2} \ ds \\
			&\quad + \sum_{i=1}^4 \int_0^t H_i(s) \ ds + H_5(t).
		\end{aligned}
	\end{equation}
	We deduce from the continuous embedding $ \H^1 \hookrightarrow \L^4 $, Gagliardo-Nirenberg inequality and Young's inequality that
	\begin{align*}
		H_1 
		&= \lb \Delta^2 m_n^\epsilon, m_n^\epsilon \times \Delta m_n^\epsilon \rb_{\L^2} \\
		&\leq |\Delta^2 m_n^\epsilon|_{\L^2} |m_n^\epsilon|_{\L^4} |\Delta m_n^\epsilon|_{\L^4} \\ 
		&\leq c |\Delta^2 m_n^\epsilon|_{\L^2} |m_n^\epsilon|_{\H^1} \la |\Delta^2 m_n^\epsilon|_{\L^2}^\frac{1}{4} |\Delta m_n^\epsilon|_{\L^2}^\frac{3}{4} + |\Delta m_n^\epsilon|_{\L^2} \ra \\ 
		&\leq \frac{5}{8}(\delta \epsilon)^{\frac{8}{5} \beta} |\Delta^2 m_n^\epsilon|_{\L^2}^2 
		+ c (\delta \epsilon)^{-\frac{8}{3}\beta} |m_n^\epsilon|_{\H^1}^\frac{8}{3} |\Delta m_n^\epsilon|_{\L^2}^2 \\
		&\quad + \frac{1}{2} (\delta \epsilon)^{\frac{8}{5} \beta} |\Delta^2 m_n^\epsilon|_{\L^2}^2 + c(\delta \epsilon)^{-\frac{8}{5}\beta} |m_n^\epsilon|_{\H^1}^2 |\Delta m_n^\epsilon|_{\L^2}^2 \\
		&\leq 
		\frac{9}{8}(\delta \epsilon)^{\frac{8}{5} \beta} |\Delta^2 m_n^\epsilon|_{\L^2}^2 
		+ c (\delta \epsilon)^{-\frac{8}{3}\beta} \la |m_n^\epsilon|_{\H^1}^\frac{8}{3} + |m_n^\epsilon|_{\H^1}^2 \ra |\Delta m_n^\epsilon|_{\L^2}^2,
	\end{align*}
	for $ \delta \in (0,1) $ and $ \beta >0 $. 
	Similarly, since $ \H^1 \hookrightarrow \L^8 $, we have 
	\begin{align*}
		H_2 
		&= - \lb \Delta^2 m_n^\epsilon, (1+ |m_n^\epsilon|^2) m_n^\epsilon \rb_{\L^2} \\
		&= - |\Delta m_n^\epsilon|^2_{\L^2} + \lb \nabla \Delta m_n^\epsilon, 2 \langle m_n^\epsilon, \nabla m_n^\epsilon \rangle m_n^\epsilon + |m_n^\epsilon|^2 \nabla m_n^\epsilon \rb_{\L^2} \\
		&\leq - |\Delta m_n^\epsilon|^2_{\L^2} + 3|\nabla \Delta m_n^\epsilon|_{\L^2} |m_n^\epsilon|_{\L^8}^2 |\nabla m_n^\epsilon|_{\L^4} \\
		&\leq - |\Delta m_n^\epsilon|^2_{\L^2} + c |\nabla \Delta m_n^\epsilon|_{\L^2} |m_n^\epsilon|_{\H^1}^2 |\nabla m_n^\epsilon|_{\H^1} \\
		&\leq - |\Delta m_n^\epsilon|^2_{\L^2} + \frac{1}{2}\delta |\nabla \Delta m_n^\epsilon|^2_{\L^2} 
		+ c \delta^{-1} \la |m_n^\epsilon|_{\H^1}^6 + |m_n^\epsilon|_{\H^1}^4 |\Delta m_n^\epsilon|_{\L^2}^2 \ra.
	\end{align*}	
	For the Stratonovich correction, 
	\begin{align*}
		H_3 
		&= \frac{1}{2} \sum_{k=1}^n \lb \Delta^2 m_n^\epsilon, \la m_n^\epsilon \times g_k \ra \times g_k \rb_{\L^2} \\
		&\leq -\frac{1}{2} \sum_{k=1}^n \lb \nabla \Delta m_n^\epsilon, \la \nabla m_n^\epsilon \times g_k + m_n^\epsilon \times \nabla g_k \ra \times g_k + (m_n^\epsilon \times g_k) \times \nabla g_k \rb_{\L^2} \\
		&\leq \frac{1}{4}\delta |\nabla \Delta m_n^\epsilon|^2_{\L^2} + c\delta^{-1} |m_n^\epsilon|^2_{\H^1}.
	\end{align*}
	For the It{\^o} correction, 
	\begin{align*}
		H_4
		&= \frac{1}{2} \sum_{k=1}^n |\Delta (g_k + m_n^\epsilon \times g_k) |^2_{\L^2} 
		\leq c (1+|\Delta m_n^\epsilon|^2_{\L^2} + |m_n^\epsilon|^2_{\H^1}). 
	\end{align*}	
	For the diffusion term, we first observe that for $ k \geq 1 $, 
	\begin{align*}
		\lb \Delta^2 m_n^\epsilon, g_k + m_n^\epsilon \times g_k \rb_{\L^2}
		&= -\lb \nabla \Delta m_n^\epsilon, \nabla g_k + \nabla m_n^\epsilon \times g_k + m_n^\epsilon \times \nabla g_k \rb_{\L^2} \\
		&\leq c|\nabla \Delta m_n^\epsilon|_{\L^2} \la |\nabla g_k|_{\L^\infty} (1+|m_n^\epsilon|_{\L^2}) + |g_k|_{\L^\infty} |\nabla m_n^\epsilon|_{\L^2} \ra,
	\end{align*}
	and
	\begin{align*}
		\lb \Delta^2 m_n^\epsilon, g_k + m_n^\epsilon \times g_k \rb_{\L^2}^2
		&\leq c |g_k|^2_{{\b W}^{1,\infty}} |\nabla \Delta m_n^\epsilon|_{\L^2}^2 \la 1+|m_n^\epsilon|^2_{\H^1}\ra.
	\end{align*}
	Then by Burkholder-Davis-Gundy inequality, for $ p \in [1,\infty) $, 
	\begin{align*}
		\E \left[ \sup_{t \in [0,T]} |H_5|^p \right]
		&\leq b_p \E \left[ \la \sum_{k=1}^n \int_0^T \lb \Delta^2 m_n^\epsilon(t), g_k + m_n^\epsilon(t) \times g_k \rb_{\L^2}^2 \ dt \ra^\frac{p}{2} \right] \\
		&\leq c b_p \E \left[ \la \int_0^T \sum_{k=1}^n |g_k|^2_{{\b W}^{1,\infty}} |\nabla \Delta m_n^\epsilon(t)|^2_{\L^2} \la 1 + |m_n^\epsilon(t)|^2_{\H^1} \ra \ dt \ra^\frac{p}{2} \right] \\
		&\leq c b_p \E \left[ \sup_{t \in [0,T]} 2^{\frac{p}{2}} \delta^{-\frac{p}{2}} \la 1+|m_n^\epsilon(t)|^2_{\H^1}\ra^{\frac{p}{2}} \la \int_0^T \frac{\delta}{2} |\nabla \Delta m_n^\epsilon(t)|^2_{\L^2} \ dt \ra^{\frac{p}{2}} \right] \\
		&\leq 2^{p-1} c b_p^2 \delta^{-p} \E \left[ \sup_{t \in [0,T]} \la 1+|m_n^\epsilon(t)|^2_{\H^1} \ra^p \right] \\
		&\quad + 2^{-(p+1)} \delta^p \E \left[\la \int_0^T |\nabla \Delta m_n^\epsilon(t)|^2_{\L^2} \ dt \ra^p \right].
	\end{align*}	
	
	Let $ \beta \geq \frac{5}{8} $. Then for $ 0< \delta < \frac{4}{9} $,
	\begin{align*}
		-\epsilon + \frac{9}{8}(\delta \epsilon)^{\frac{8}{5}\beta} < -\frac{1}{2}\epsilon, \quad 
		-1 + \frac{1}{2}\delta + \frac{1}{4}\delta < -\frac{1}{2}.
	\end{align*}
	Taking the $ p $th power, the supremum over $ t $ and the expectation on both sides of \eqref{pf: 1/2|Delta m_n^e|^2},
	\begin{align*}
		&\E \left[ \sup_{t \in [0,T]}|\Delta m_n^\epsilon(t)|^{2p}_{\L^2} 
		+ \frac{1}{2} \la \int_0^T |\nabla \Delta m_n^\epsilon(t)|^2_{\L^2} \ dt \ra^p + \epsilon^p \la \int_0^T |\Delta^2 m_n^\epsilon(t)|^2_{\L^2} \ dt \ra^p \right] \\
		&\leq 
		c\E \left[|\Delta m_n^\epsilon(0)|^{2p}_{\L^2} \right] 
		+ c \E \left[ \sup_{t \in [0,T]} \la 1+ |m_n^\epsilon(t)|^2_{\H^1} + |m_n^\epsilon(t)|^6_{\H^1}\ra^p  \right] 
		+ c\E \left[ \la \int_0^T |\Delta m_n^\epsilon(t)|^2_{\L^2} \ dt \ra^p \right] \\
		&\quad + c \epsilon^{-\frac{8}{3}\beta p} \E \left[ \sup_{t \in [0,T]} \la |m_n^\epsilon(t)|_{\H^1}^{\frac{8}{3}p} + |m_n^\epsilon(t)|_{\H^1}^{2p} + |m_n^\epsilon(t)|_{\H^1}^{4p} \ra \la \int_0^T |\Delta m_n^\epsilon(t)|_{\L^2}^2 \ dt \ra^p \right] \\
		&\leq c \la 1+ \epsilon^{-\frac{8}{3}\beta p} \ra,
	\end{align*}
	for some constant $ c>0 $ that may depend on $ C_1, C_g, \delta, p, T $ but not $ n, \epsilon $. 
	Here, the last inequality holds by the estimate of $ m_n^\epsilon $ in $ L^{2p}(\Omega;L^\infty(0,T; \H^1) \cap L^2(0,T; \H^2)) $ in Lemma \ref{Lemma: mn^e est H1}. 
	Taking $ \beta = \frac{5}{8} $ and multiplying both sides by $ \epsilon^{\frac{5}{3}p} $ yields the desired inequality.

\subsection{Proof of Lemma \ref{Lemma: mn^e est H3}}
	Applying It{\^o}'s lemma to $ \frac{1}{2}|\nabla \Delta m_n^\epsilon|_{\L^2}^2 $,
	\begin{align*}
		\frac{1}{2}|\nabla \Delta m_n^\epsilon(t)|^2_{\L^2} 
		&= 
		\frac{1}{2} |\nabla \Delta m_n^\epsilon(0)|^2_{\L^2} - \epsilon \int_0^t |\nabla \Delta^2 m_n^\epsilon(s)|_{\L^2}^2 \ ds - \int_0^t |\Delta^2 m_n^\epsilon(s)|_{\L^2}^2 \ ds \\
		&\quad -\int_0^{t} \lb \Delta^3 m_n^\epsilon(s), m_n^\epsilon(s) \times \Delta m_n^\epsilon(s) \rb_{\L^2} \ ds \\
		&\quad + \int_0^{t} \lb \Delta^3 m_n^\epsilon(s), (1+|m_n^\epsilon(s)|^2) m_n^\epsilon(s) \rb_{\L^2} \ ds \\
		&\quad - \frac{1}{2} \sum_{k=1}^n \int_0^{t} \lb \Delta^3 m_n^\epsilon(s), (m_n^\epsilon(s) \times g_k) \times g_k \rb_{\L^2} \ ds \\
		&\quad + \frac{1}{2} \sum_{k=1}^n \int_0^{t} |\nabla \Delta \la m_n^\epsilon(s) \times g_k + g_k \ra |^2_{\L^2} \ ds \\
		&\quad - \sum_{k=1}^n \int_0^{t} \lb \Delta^3 m_n^\epsilon(s), m_n^\epsilon(s) \times g_k + g_k \rb_{\L^2} dW_k(s) \\
		&= \frac{1}{2} |\nabla \Delta m_n^\epsilon(0)|^2_{\L^2} - \epsilon \int_0^t |\nabla \Delta^2 m_n^\epsilon(s)|_{\L^2}^2 \ ds - \int_0^t |\Delta^2 m_n^\epsilon(s)|_{\L^2}^2 \ ds \\
		&\quad + \sum_{i=1}^4 \int_0^{t} T_i(s) \ ds + T_5(t). 
	\end{align*}
	For $ \delta \in (0,1) $, we have
	\begin{align*}
		T_1 
		&= -\lb \Delta^3 m_n^\epsilon, m_n^\epsilon \times \Delta m_n^\epsilon \rb_{\L^2} \\
		&= - 2 \lb \Delta^2 m_n^\epsilon, \nabla m_n^\epsilon \times \nabla \Delta m_n^\epsilon \rb_{\L^2} \\
		&\leq 2 |\Delta^2 m_n^\epsilon|_{\L^2} |\nabla m_n^\epsilon|_{\L^4} |\nabla \Delta m_n^\epsilon|_{\L^4} \\
		&\leq c |\Delta^2 m_n^\epsilon|_{\L^2} |\nabla m_n^\epsilon|_{\L^4} \la |\nabla \Delta m_n^\epsilon|_{\L^2}^\frac{1}{2} |\Delta^2 m_n^\epsilon|_{\L^2}^\frac{1}{2} + |\nabla \Delta m_n^\epsilon|_{\L^2} \ra \\
		&\leq \la \frac{3}{4}\delta^{\frac{4}{3}} + \frac{1}{2}\delta \ra |\Delta^2 m_n^\epsilon|^2_{\L^2} + c \la \delta^{-4} |\nabla m_n^\epsilon|_{\L^4}^4 + \delta^{-1} |\nabla m_n^\epsilon|_{\L^4}^2 \ra |\nabla \Delta m_n^\epsilon|_{\L^2}^2 \\
		&\leq \frac{5}{4}\delta |\Delta^2 m_n^\epsilon|^2_{\L^2} + c \delta^{-4} \la |\nabla m_n^\epsilon|_{\L^2}^2 |\Delta m_n^\epsilon|_{\L^2}^2 + |\nabla m_n^\epsilon|_{\L^2}^4 + |\nabla m_n^\epsilon|_{\H^1}^2 \ra |\nabla \Delta m_n^\epsilon|_{\L^2}^2 \\
		&= \frac{5}{4}\delta |\Delta^2 m_n^\epsilon|^2_{\L^2} + c (\epsilon^{\frac{8}{3}}\delta^4)^{-1} \la |\nabla m_n^\epsilon|_{\L^2}^2 +1 \ra \la \epsilon^{\frac{5}{3}}|\Delta m_n^\epsilon|_{\L^2}^2 \ra \la \epsilon |\nabla \Delta m_n^\epsilon|_{\L^2}^2 \ra \\
		&\quad + c (\epsilon \delta^4)^{-1} \la |\nabla m_n^\epsilon|_{\L^2}^2 + |\nabla m_n^\epsilon|_{\L^2}^4 \ra \la \epsilon |\nabla \Delta m_n^\epsilon|_{\L^2}^2 \ra.
	\end{align*}
	Similarly, for $ T_2 $, we first observe that
	\begin{align*}
		\Delta \la |m_n^\epsilon|^2 m_n^\epsilon \ra 
		&= |m_n^\epsilon|^2 \Delta m_n^\epsilon + 2 |\nabla m_n^\epsilon|^2 m_n^\epsilon + 2 \langle m_n^\epsilon, \Delta m_n^\epsilon \rangle m_n^\epsilon + 4\langle m_n^\epsilon, \nabla m_n^\epsilon \rangle \nabla m_n^\epsilon.
	\end{align*}
	Then,
	\begin{align*}
		T_2
		&= \lb \Delta^3 m_n^\epsilon, (1+|m_n^\epsilon|^2) m_n^\epsilon \rb_{\L^2} \\
		&= -|\nabla \Delta m_n^\epsilon|^2_{\L^2} + \lb \Delta^2 m_n^\epsilon, \Delta \la |m_n^\epsilon|^2 m_n^\epsilon \ra \rb_{\L^2} \\
		&\leq 
		-|\nabla \Delta m_n^\epsilon|^2_{\L^2} 
		+ c|\Delta^2 m_n^\epsilon|_{\L^2} \la |\nabla m_n^\epsilon|^2_{\L^8} |m_n^\epsilon|_{\L^4} + |\Delta m_n^\epsilon|_{\L^4} |m_n^\epsilon|^2_{\L^8} \ra \\
		&\leq 
		-|\nabla \Delta m_n^\epsilon|^2_{\L^2} 
		+ c|\Delta^2 m_n^\epsilon|_{\L^2} \la |\nabla \Delta m_n^\epsilon|_{\L^2}^\frac{3}{4} |\nabla m_n^\epsilon|_{\L^2}^\frac{5}{4} + |\nabla m_n^\epsilon|_{\L^2}^2 \ra |m_n^\epsilon|_{\H^1} \\
		&\quad + c|\Delta^2 m_n^\epsilon|_{\L^2} \la |\Delta^2 m_n^\epsilon|_{\L^2}^\frac{1}{2} |\nabla m_n^\epsilon|_{\L^2}^\frac{1}{2} + |\nabla m_n^\epsilon|_{\L^2} \ra |m_n^\epsilon|^2_{\H^1} \\
		&\leq 
		-|\nabla \Delta m_n^\epsilon|^2_{\L^2} 
		+ \frac{1}{2}\delta |\Delta^2 m_n^\epsilon|^2_{\L^2} 
		+ c \delta^{-1} \la |\nabla \Delta m_n^\epsilon|_{\L^2}^\frac{3}{2} |\nabla m_n^\epsilon|_{\L^2}^\frac{5}{2} |m_n^\epsilon|^2_{\H^1} + |m_n^\epsilon|^6_{\H^1} \ra \\
		&\quad + \frac{3}{4} \delta^\frac{4}{3} |\Delta^2 m_n^\epsilon|_{\L^2}^2 + c \delta^{-4} |\nabla m_n^\epsilon|_{\L^2}^2 |m_n^\epsilon|^8_{\H^1} 
		+ \frac{1}{2} \delta |\Delta^2 m_n^\epsilon|_{\L^2}^2 + c \delta^{-1} |m_n^\epsilon|_{\H^1}^6 \\
		&\leq 
		\la -1 + c \delta^{-1} \ra |\nabla \Delta m_n^\epsilon|^2_{\L^2} 
		+ \frac{5}{4}\delta |\Delta^2 m_n^\epsilon|^2_{\L^2} 
		+ c \delta^{-4} \la |m_n^\epsilon|_{\H^1}^6 + |m_n^\epsilon|_{\H^1}^{10} + |m_n^\epsilon|_{\H^1}^{18} \ra.
	\end{align*}
	For the Stratonovich correction, we have 
	\begin{align*}
		T_3 
		&= -\sum_{k=1}^n \lb \Delta^3 m_n^\epsilon, (m_n^\epsilon \times g_k) \times g_k \rb_{\L^2} \\
		&= -\sum_{k=1}^n \lb \Delta^2 m_n^\epsilon, (\Delta m_n^\epsilon \times g_k) \times g_k + (m_n^\epsilon \times \Delta g_k) \times g_k + 2 (\nabla m_n^\epsilon \times \nabla g_k) \times g_k \rb_{\L^2} \\
		&\quad -\sum_{k=1}^n \lb \Delta^2 m_n^\epsilon, (m_n^\epsilon \times g_k) \times \Delta g_k + 2 \la \nabla m_n^\epsilon \times g_k + m_n^\epsilon \times \nabla g_k \ra \times \nabla g_k \rb_{\L^2} \\
		&\leq 
		\sum_{k=1}^n |\Delta^2 m_n^\epsilon|_{\L^2} 
		\la |\Delta m_n^\epsilon|_{\L^4} |g_k|^2_{\L^8} 
		+ 2 |m_n^\epsilon|_{\L^\infty} |g_k|_{\L^\infty} |\Delta g_k|_{\L^2} \ra \\
		&\quad + \sum_{k=1}^n |\Delta^2 m_n^\epsilon|_{\L^2} \la 4 |\nabla m_n^\epsilon|_{\L^2} |g_k|_{\L^\infty} |\nabla g_k|_{\L^\infty} + |m_n^\epsilon|_{\L^2} |\nabla g_k|_{\L^\infty}^2 \ra \\
		&\leq 
		c \la |\Delta^2 m_n^\epsilon|^\frac{3}{2}_{\L^2} |\nabla m_n^\epsilon|_{\L^2}^\frac{1}{2} + |\Delta^2 m_n^\epsilon|_{\L^2} |\nabla m_n^\epsilon|_{\L^2} \ra \sum_{k=1}^n |g_k|^2_{\H^1} \\
		&\quad + c \la |\Delta^2 m_n^\epsilon|^\frac{5}{4}_{\L^2} |m_n^\epsilon|_{\L^2}^\frac{3}{4} + |\Delta^2 m_n^\epsilon|_{\L^2} |m_n^\epsilon|_{\L^2} \ra \sum_{k=1}^n |g_k|_{\L^\infty} |\Delta g_k|_{\L^2} \\
		&\quad + c |\Delta^2 m_n^\epsilon|_{\L^2} |m_n^\epsilon|_{\H^1} \sum_{k=1}^n |g_k|_{{\b W}^{1,\infty}}^2 \\
		&\leq 
		\frac{3}{4}\delta^\frac{4}{3} |\Delta^2 m_n^\epsilon|^2_{\L^2} + c \delta^{-4} |\nabla m_n^\epsilon|_{\L^2}^2 
		+ \frac{1}{2} \delta |\Delta^2 m_n^\epsilon|_{\L^2}^2 + c \delta^{-1} |\nabla m_n^\epsilon|_{\L^2}^2 \\
		&\quad + \frac{5}{8} \delta^\frac{8}{5} |\Delta^2 m_n^\epsilon|^2_{\L^2} + c \delta^{-\frac{8}{3}} |m_n^\epsilon|_{\L^2}^2 
		+ \frac{1}{2} \delta |\Delta^2 m_n^\epsilon|^2_{\L^2} + \frac{1}{2} c\delta^{-1} |m_n^\epsilon|_{\L^2}^2  \\
		&\quad + \frac{1}{2} \delta |\Delta^2 m_n^\epsilon|^2_{\L^2} + c \delta^{-1} |m_n^\epsilon|^2_{\H^1} \\
		&\leq 
		\frac{23}{8} \delta |\Delta^2 m_n^\epsilon|^2_{\L^2} 
		+ c \delta^{-4} |m_n^\epsilon|_{\H^1}^2.
	\end{align*}
	For the It{\^o} correction, 
	\begin{align*}
		T_4 
		&= \sum_{k=1}^n |\nabla \Delta \la m_n^\epsilon \times g_k + g_k \ra|^2_{\L^2} \\
		&\leq c \sum_{k=1}^n \la |\nabla \Delta m_n^\epsilon \times g_k|_{\L^2}^2 + |\Delta m_n^\epsilon \times \nabla g_k|_{\L^2}^2 + |\nabla m_n^\epsilon \times \Delta g_k|^2_{\L^2} + |m_n^\epsilon \times \nabla \Delta g_k|_{\L^2}^2 + |\nabla \Delta g_k|^2_{\L^2} \ra \\
		&\leq c |\nabla \Delta m_n^\epsilon|^2_{\L^2} \sum_{k=1}^n |g_k|_{\L^\infty}^2 + |\Delta m_n^\epsilon|^2_{\L^4} \sum_{k=1}^n |\nabla g_k|_{\L^2} |\Delta g_k|_{\L^2} 
		+ c |\nabla m_n^\epsilon|^2_{\L^4} \sum_{k=1}^n |\Delta g_k|_{\L^2} |\nabla \Delta g_k|_{\L^2} \\
		&\quad + c \la |m_n^\epsilon|_{\L^\infty}^2 + 1 \ra \sum_{k=1}^n |\nabla \Delta g_k|^2_{\L^2} \\
		&\leq c |\nabla \Delta m_n^\epsilon|^2_{\L^2} + c \la |\Delta^2 m_n^\epsilon|_{\L^2} |\nabla m_n^\epsilon|_{\L^2} + |\nabla m_n^\epsilon|_{\L^2}^2 \ra
		+ c \la |\nabla \Delta m_n^\epsilon|^\frac{1}{2}_{\L^2} |\nabla m_n^\epsilon|_{\L^2}^\frac{3}{2} + |\nabla m_n^\epsilon|_{\L^2}^2 \ra \\
		&\quad + c \la |\nabla \Delta m_n^\epsilon|_{\L^2}^{\frac{2}{3}} |m_n^\epsilon|_{\L^2}^\frac{4}{3} + |m_n^\epsilon|_{\L^2}^2 + 1 \ra \\
		&\leq c |\nabla \Delta m_n^\epsilon|^2_{\L^2} + \delta |\Delta^2 m_n^\epsilon|_{\L^2}^2 + c \la \delta^{-1} |m_n^\epsilon|^2_{\H^1} + 1 \ra.
	\end{align*}
	For the diffusion term, we observe that
	\begin{align*}
		\lb \Delta^2 m_n^\epsilon, \Delta m_n^\epsilon \times g_k \rb_{\L^2} 
		&= -\lb \nabla \Delta m_n^\epsilon, \Delta m_n^\epsilon \times \nabla g_k \rb_{\L^2} \\
		&= \lb \Delta g_k \times \nabla \Delta m_n^\epsilon + \nabla g_k \times \Delta^2 m_n^\epsilon, \nabla m_n^\epsilon \rb_{\L^2}, 
	\end{align*}
	and 
	\begin{align*}
		\lb \Delta^3 m_n^\epsilon, m_n^\epsilon \times g_k + g_k \rb_{\L^2} 
		&= \lb \Delta^2 m_n^\epsilon, \Delta m_n^\epsilon \times g_k \rb_{\L^2} 
		+ \lb \Delta^2 m_n^\epsilon, m_n^\epsilon \times \Delta g_k + 2 \nabla m_n^\epsilon \times \nabla g_k + \Delta g_k \rb_{\L^2} \\
		&= \lb \nabla \Delta m_n^\epsilon, \nabla m_n^\epsilon \times \Delta g_k \rb_{\L^2} 
		+ \lb \Delta^2 m_n^\epsilon, m_n^\epsilon \times \Delta g_k + 3 \nabla m_n^\epsilon \times \nabla g_k + \Delta g_k \rb_{\L^2}. 
	\end{align*}
	Thus, 
	\begin{align*}
		\sum_{k=1}^n \lb \Delta^3 m_n^\epsilon, m_n^\epsilon \times g_k + g_k \rb_{\L^2}^2 
		&\leq  
		|\nabla \Delta m_n^\epsilon|_{\L^2}^2 |\nabla m_n^\epsilon|_{\L^2}^2 \sum_{k=1}^n|\Delta g_k|^2_{\L^\infty} \\
		&\quad + c |\Delta^2 m_n^\epsilon|_{\L^2}^2 \la |m_n^\epsilon|_{\L^2}^2 \sum_{k=1}^n|\Delta g_k|^2_{\L^\infty} + |\nabla m_n^\epsilon|_{\L^2}^2 \sum_{k=1}^n|\nabla g_k|_{\L^\infty}^2 \ra \\
		&\quad + c |\Delta^2 m_n^\epsilon|_{\L^2}^2 \sum_{k=1}^n|\Delta g_k|^2_{\L^2} \\
		&\leq 
		c |\nabla \Delta m_n^\epsilon|_{\L^2}^2 |\nabla m_n^\epsilon|_{\L^2}^2 
		+ c |\Delta^2 m_n^\epsilon|_{\L^2}^2 \la |m_n^\epsilon|_{\H^1}^2 + 1 \ra.
	\end{align*}
	Then for $ p \in [1,\infty) $, by Burkholder-Davis-Gundy inequality, 
	\begin{align*}
		\E \left[ \sup_{t \in[0,T]} \left| T_5(t) \right|^p \right] 
		&\leq b_p \E \left[ \la \sum_{k=1}^n \int_0^{T} \lb \Delta^3 m_n^\epsilon(t), m_n^\epsilon(t) \times g_k + g_k \rb_{\L^2}^2 \ dt \ra^{\frac{p}{2}} \right] \\
		&\leq c b_p \E \left[ \sup_{t \in [0,T]} |\nabla m_n^\epsilon(t)|_{\L^2}^p \la \int_0^T |\nabla \Delta m_n^\epsilon(t)|_{\L^2}^2 \ dt \ra^{\frac{p}{2}} \right] \\
		&\quad + c b_p \E \left[ 2^\frac{p}{2} \delta^{-\frac{p}{2}}\sup_{t \in [0,T]} \la |m_n^\epsilon(t)|^2_{\H^1} +1 \ra^\frac{p}{2} \la \int_0^T \frac{\delta}{2} |\Delta^2 m_n^\epsilon(t)|_{\L^2}^2 \ dt \ra^\frac{p}{2} \right] \\
		&\leq 
		\frac{1}{2} c b_p \E \left[ \sup_{t \in [0,T]} |\nabla m_n^\epsilon(t)|_{\L^2}^{2p} + \int_0^T |\nabla \Delta m_n^\epsilon(t)|_{\L^2}^{2p} \ dt \right] \\
		&\quad + 2^{p-1} c b_p^2 \delta^{-p} \E \left[ \sup_{t \in [0,T]} \la |m_n^\epsilon|^2_{\H^1} + 1 \ra^p \right] \\
		&\quad + 2^{-(p+1)}\delta^p \E \left[ \la \int_0^T |\Delta^2 m_n^\epsilon(t)|_{\L^2}^2 \ dt \ra^p \right].
	\end{align*}	
	We choose a sufficiently small $ \delta \in (0,1) $ such that the sum of coefficients of $ |\Delta^2 m_n^\epsilon|_{\L^2}^2 $ in $ T_1, \ldots, T_4 $ is negative; more precisely, $ -1 + \frac{51}{8} \delta < \frac{1}{2} $. 
	Then combining the estimates, we deduce that for $ p \in [1,\infty) $, 
	\begin{align*}
		&\E \left[ \sup_{t \in [0,T]} |\nabla \Delta m_n^\epsilon(t)|^{2p}_{\L^2} + \epsilon^p \la \int_0^T |\nabla \Delta^2 m_n^\epsilon(t)|_{\L^2}^2 \ dt \ra^p + \frac{1}{2} \la \int_0^T |\Delta^2 m^\epsilon(t)|^2_{\L^2} \ dt \ra^p \right] \\
		&\leq 
		c \E \left[ |\nabla \Delta m_n^\epsilon(0)|_{\L^2}^{2p} +1 \right] \\
		&\quad + c \la \epsilon^{\frac{8}{3}} \delta^{4} \ra^{-p} \E \left[ \sup_{t \in [0,T]} \la |\nabla m_n^\epsilon|_{\L^2}^2 + 1 \ra^p \la \epsilon^{\frac{5}{3}} |\Delta m_n^\epsilon|_{\L^2}^{2} \ra^p \la \epsilon \int_0^T |\nabla \Delta m_n^\epsilon|^2_{\L^2} \ dt \ra^p \right] \\
		&\quad + c \la \epsilon \delta^{4} \ra^{-p} \E \left[ \sup_{t \in [0,T]} \la |\nabla m_n^\epsilon|_{\L^2}^2 + |\nabla m_n^\epsilon|_{\L^2}^4 \ra^p \la \epsilon \int_0^T |\nabla \Delta m_n^\epsilon|^2_{\L^2} \ dt \ra^p \right] \\
		&\quad + c \delta^{-p} \E \left[ \int_0^T \sup_{s \in [0,t]} |\nabla \Delta m_n^\epsilon(s)|_{\L^2}^{2p} \ dt \right] \\
		&\quad + c \delta^{-4} \E \left[ \sup_{t \in [0,T]} \la |m_n^\epsilon(t)|_{\H^1}^{2p} + |m_n^\epsilon(t)|_{\H^1}^{6p} + |m_n^\epsilon(t)|_{\H^1}^{8p} + |m_n^\epsilon(t)|_{\H^1}^{18p}  \ra \right]. 
	\end{align*}
	By Lemmas \ref{Lemma: mn^e est H1} and \ref{Lemma: mn^e est H2}, 
	\begin{align*}
		\E \left[ \sup_{t \in [0,T]} |\nabla \Delta m_n^\epsilon(t)|^{2p}_{\L^2} + \la \epsilon \int_0^T |\nabla \Delta^2 m_n^\epsilon(t)|_{\L^2}^2 \ dt \ra^p \right] 
		&\leq c \epsilon^{-\frac{8}{3}p} + c \E \left[ \int_0^T \sup_{s \in [0,t]} |\nabla \Delta m_n^\epsilon(s)|_{\L^2}^{2p} \ dt \right] + c. 
	\end{align*}
	Finally, by Gronwall's lemma,
	\begin{align*}
		\E \left[ \sup_{t \in [0,T]} |\nabla \Delta m_n^\epsilon(t)|^{2p}_{\L^2} \right] &\leq c (1+\epsilon^{-\frac{8}{3}p}) e^c,
	\end{align*}
	for some constant $ c $ that may depend on $ C_1, C_g, T, p $ but not on $ n $ and $ \epsilon $.

\subsection{Proof of Theorem \ref{Lemma: m 1D regularity}}
	We consider the Faedo-Galerkin approximation as in Section \ref{Section: uniform estimates}, with $ \epsilon = 0 $ and $ m_n :=: m_n^0 $. 
	Since the $ L^{2p}(\Omega;L^\infty(0,T;\H^1)) $-estimate of $ m_n $ was obtained in \cite[Lemma 3.2 and 3.3]{BGL_sLLB}, we only focus on $ \Delta m_n $. 
	By It{\^o}'s lemma, 
	\begin{align*}
		\frac{1}{2}|\Delta m_n(t)|^2_{\L^2} -\frac{1}{2} |\Delta m_n(0)|^2_{\L^2} 
		&= 
		\int_0^{t} \lb \Delta^2 m_n(s), \Delta m_n(s) + m_n(s) \times \Delta m_n(s) \rb_{\L^2} \ ds \\
		&\quad - \int_0^{t} \lb \Delta^2 m_n(s), (1+|m_n(s)|^2) m_n(s) \rb_{\L^2} \ ds \\
		&\quad + \frac{1}{2} \sum_{k=1}^\infty \int_0^{t} \lb \Delta^2 m_n(s), (m_n(s) \times g_k) \times g_k \rb_{\L^2} \ ds \\
		&\quad + \frac{1}{2} \sum_{k=1}^\infty \int_0^{t} |\Delta \la m_n(s) \times g_k + g_k \ra |^2_{\L^2} \ ds \\
		&\quad + \sum_{k=1}^\infty \int_0^{t} \lb \Delta^2 m_n(s), m_n(s) \times g_k + g_k \rb_{\L^2} dW_k(s) \\
		&=: \int_0^{t} H_1(s) \ ds + \int_0^{t} H_2(s) \ ds + \int_0^{t} H_3(s) \ ds + \int_0^{t} H_4(s) \ ds + H_5(t). 
	\end{align*}
	We have
	\begin{align*}
		H_1 
		&= \lb \Delta^2 m_n, \Delta m_n + m_n \times \Delta m_n \rb_{\L^2} \\
		&= -|\nabla \Delta m_n|^2_{\L^2} - \lb \nabla \Delta m_n, \nabla m_n \times \Delta m_n \rb_{\L^2} \\
		&\leq -|\nabla \Delta m_n|^2_{\L^2} + |\nabla \Delta m_n|_{\L^2} |\nabla m_n|_{\L^4} |\Delta m_n|_{\L^4} \\
		&\leq -|\nabla \Delta m_n|^2_{\L^2} + c|\nabla \Delta m_n|_{\L^2} \la |\nabla \Delta m_n|_{\L^2}^\frac{1}{8} |\nabla m_n|_{\L^2}^\frac{7}{8} + |\nabla m_n|_{\L^2} \ra \la |\nabla \Delta m_n|_{\L^2}^\frac{5}{8} |\nabla m_n|_{\L^2}^\frac{3}{8} + |\nabla m_n|_{\L^2} \ra \\
		&\leq -|\nabla \Delta m_n|^2_{\L^2} 
		+ c|\nabla \Delta m_n|_{\L^2}^\frac{7}{4} |\nabla m_n|_{\L^2}^\frac{5}{4} 
		+ c|\nabla \Delta m_n|_{\L^2}^\frac{9}{8} |\nabla m_n|_{\L^2}^\frac{15}{8} 
		+ c|\nabla \Delta m_n|_{\L^2}^\frac{13}{8} |\nabla m_n|_{\L^2}^\frac{11}{8} 
		+ c|\nabla m_n|_{\L^2}^2 \\
		&\leq (-1+\delta) |\nabla \Delta m_n|^2_{\L^2} + c \la |\nabla m_n|_{\L^2}^{10} + |\nabla m_n|_{\L^2}^{\frac{30}{7}} + |\nabla m_n|_{\L^2}^{\frac{22}{3}} + |\nabla m_n|_{\L^2}^2 \ra \\
		&\leq (-1+\delta) |\nabla \Delta m_n|^2_{\L^2} + c_1 \la |\nabla m_n|_{\L^2}^{10} +1 \ra,
	\end{align*}
	where $ c_1 $ depends on $ \delta^{-1} $. 
	The estimates of $ H_2, H_3 $ and $ H_4 $ are similar to those in the proof of Lemma \ref{Lemma: mn^e est H2}, so we omit some details here. 
	Since $ \H^1 \hookrightarrow \L^\infty $ for $ 1 $-dimensional domain, we have 
	\begin{align*}
		H_2
		&= -\lb \Delta^2 m_n, (1+|m_n|^2) m_n \rb_{\L^2} \\
		&\leq -|\Delta m_n|^2_{\L^2} + c |\nabla \Delta m_n|_{\L^2} |m_n|_{\L^\infty}^2 |\nabla m_n|_{\L^2} \\
		&\leq -|\Delta m_n|^2_{\L^2} + c |\nabla \Delta m_n|_{\L^2} |m_n|_{\H^1}^2 |\nabla m_n|_{\L^2} \\
		&\leq -|\Delta m_n|^2_{\L^2} + \delta |\nabla \Delta m_n|_{\L^2}^2 + c_2 |m_n|_{\H^1}^6. 
	\end{align*}
	For the Stratonovich correction and It{\^o} correction,  
	\begin{align*}
		H_3 
		&\leq \delta |\nabla \Delta m_n|^2_{\L^2} + c_3 |m_n|^2_{\H^1}, \\
		H_4 
		&\leq c_4 \la 1 + |m_n|^2_{\H^1} + |\Delta m_n|^2_{\L^2} \ra. 
	\end{align*}
	For $ p \in [1,\infty) $, by Burkholder-Davis-Gundy inequality, 
	\begin{align*}
		\E \left[ \sup_{t \in[0,T]} \left| H_5(t) \right|^p \right] 
		&= \E \left[ \sup_{t \in [0,T]} \left| \sum_{k=1}^\infty \int_0^t \lb \Delta^2 m_n(s), m_n(s) \times g_k + g_k \rb_{\L^2} dW_k(s) \right|^p \right] \\
		&\leq b_p \E \left[ \left| \sum_{k=1}^\infty \int_0^T \lb \Delta^2m_n(s), m_n \times g_k + g_k \rb_{\L^2}^2 \ ds \right|^{\frac{p}{2}} \right] \\
		&= b_p \E \left[ \left| \sum_{k=1}^\infty \int_0^T \lb \Delta m_n(s), m_n(s) \times \Delta g_k + 2\nabla m_n(s) \times \nabla g_k + \Delta g_k \rb_{\L^2}^2 \ ds \right|^{\frac{p}{2}} \right] \\
		&\leq c b_p \E \left[ \left| \int_0^T |\Delta m_n(s)|^2_{\L^2} \la |m_n(s)|_{\L^\infty}^2 + |\nabla m_n(s)|_{\L^2}^2 +1 \ra \ ds \right|^{\frac{p}{2}} \right] \\
		&\leq c b_p \E \left[ \sup_{t \in [0,T]} \la |m_n(s)|_{\H^1}^2 +1 \ra^p \right] 
		+ c b_p \E \left[ \int_0^T |\Delta m_n(s)|^{2p}_{\L^2} \ ds \right].
	\end{align*}	
	Combining the estimates, we deduce that
	\begin{align*}
		&\frac{1}{2} \sup_{t \in [0,T]} |\Delta m_n(t)|^2_{\L^2}
		+(1-3\delta) \int_0^T |\nabla \Delta m_n(s)|^2_{\L^2} \ ds \\
		&\leq \frac{1}{2} |\Delta m_n(0)|_{\L^2}^2 + c \la |m_n|_{\H^1}^{10} + 1 \ra 
		+ c |\Delta m_n|^2_{\L^2} 
		+ \sup_{t \in[0,T]} \left| H_5(t) \right|.
	\end{align*}
	We choose $ \delta < \frac{1}{3} $. Raising both sides to power $ p $ and taking expectation, 
	\begin{align*}
		&\E \left[ \sup_{t \in [0,T]} |\Delta m_n(t)|^{2p}_{\L^2} +(1-3\delta)^p \la \int_0^T |\nabla \Delta m_n(s)|^2_{\L^2} \ ds\ra^p \right] \\
		&\leq c\E \left[|\Delta m_0|_{\L^2}^{2p} \right] + c \E \left[ \int_0^T \sup_{s \in [0,t]} |\Delta m_n(s)|^{2p}_{\L^2} \ dt \right] 
		+ c \E\left[ \sup_{t \in [0,T]} \la |m_n(t)|_{\H^1}^{10} + 1 \ra^p \right] \\
		&\leq c \E \left[ \int_0^T \sup_{s \in [0,t]} |\Delta m_n(s)|^{2p}_{\L^2} \ dt \right] + c,
	\end{align*}
	for some constant $ c $ that may depend on $ C_1, C_g, T, p $ and $ \delta^{-1} $. 
	Hence, 
	\begin{align*}
		\E \left[ \sup_{t \in [0,T]} |\Delta m_n(t)|^{2p}_{\L^2} \right] 
		&\leq ce^c, \quad
		\E \left[ \la \int_0^T |\nabla \Delta m_n(s)|^2_{\L^2} \ ds \ra^p \right]
		\leq c(e^c+1),
	\end{align*}
	by Gronwall's lemma. 
	Then we can repeat the compactness, tightness, and convergence arguments as in Section \ref{Section: tightness and convergence} to show the existence of a weak martingale solution $ m $ of the stochastic LLB equation with estimate in $ L^{2p}(\Omega; L^\infty(0,T;\H^2) \cap L^2(0,T;\H^3)) $ when $ d=1 $.
	
	As a result, using a similar argument to the  $ \mathcal{C}^\alpha([0,T];\H^1) $-estimate of $ m^\epsilon $ as in Theorem \ref{Theorem: m^e}, there exists a positive constant $ c $ such that
	\begin{align*}
		\E \left[ |m|_{\mathcal{C}^\alpha([0,T]; \H^1)}^p \right] < c, 
	\end{align*}
	for $ \alpha \in (0,\frac{1}{2}) $ and $ p \in [1,\infty) $.

\section{Proof of error estimates}\label{Section: Proofs of error estimates}

\subsection{Proof of Theorem \ref{Theorem: eh conve}}
	Taking $ \phi $ be $ P_h e_h^{(n)} $, we have
	\begin{align*}
		&\lb e_h^{(n)}-e_h^{(n-1)}, P_h e_h^{(n)} \rb_{\L^2} \\
		&= - \int_{\tau_{n-1}^{\epsilon,R}}^{\tau_n^{\epsilon,R}} \la \epsilon \lb \Delta m^\epsilon(s) - \Delta m_h^{(n)}, \Delta P_h e_h^{(n)} \rb_{\L^2} + \lb \nabla m^\epsilon(s) - \nabla m_h^{(n)}, \nabla P_h e_h^{(n)} \rb_{\L^2} \ra \ ds \\
		&\quad - \int_{\tau_{n-1}^{\epsilon,R}}^{\tau_n^{\epsilon,R}} \lb m^\epsilon(s) \times \nabla m^\epsilon(s) - m_h^{(n-1)} \times \nabla m_h^{(n)}, \nabla P_h e_h^{(n)} \rb_{\L^2} \ ds \\
		&\quad - \int_{\tau_{n-1}^{\epsilon,R}}^{\tau_n^{\epsilon,R}} \lb (1 + |m^\epsilon(s)|^2) m^\epsilon(s) - (1+|m_h^{(n-1)}|^2) m_h^{(n)}, P_h e_h^{(n)} \rb_{\L^2} \ ds \\
		&\quad + \frac{1}{2} \sum_{k=1}^\infty \int_{\tau_{n-1}^{\epsilon,R}}^{\tau_n^{\epsilon,R}} \lb \la m^\epsilon(s)-m_h^{(n-1)} \ra \times g_k, g_k \times P_h e_h^{(n)} \rb_{\L^2} \ ds \\
		&\quad + \int_{t_{n-1}}^{t_n} \lb P_h e_h^{(n)}, \la m^\epsilon(s) - m_h^{(n-1)} \ra \times dW(s) \rb_{\L^2} \frac{\Delta \tau_n^{\epsilon,R}}{\Delta t} \\
		&= \int_{\tau_{n-1}^{\epsilon,R}}^{\tau_n^{\epsilon,R}} \la I_1^{(n)} + I_2^{(n)} + I_3^{(n)} \ra \ ds + \sum_{k=1}^\infty \int_{\tau_{n-1}^{\epsilon,R}}^{\tau_n^{\epsilon,R}} I_4^{(n)} \ ds + I_5^{(n)} .
	\end{align*}
	The left-hand-side can be re-written as
	\begin{align*}
		\lb e_h^{(n)}-e_h^{(n-1)}, P_h e_h^{(n)} \rb_{\L^2} 
		&= \lb e_h^{(n)}-e_h^{(n-1)}, e_h^{(n)} \rb_{\L^2} + \lb e_h^{(n)}-e_h^{(n-1)}, P_he_h^{(n)}-e_h^{(n)} \rb_{\L^2} \\
		&= \frac{1}{2} \la |e_h^{(n)}|^2_{\L^2} - |e_h^{(n-1)}|^2_{\L^2} + |e_h^{(n)}-e_h^{(n-1)}|^2_{\L^2} \ra - I_0^{(n)}
	\end{align*}
	where $ I_0^{(n)} = \langle e_h^{(n)}-e_h^{(n-1)}, e_h^{(n)}-P_he_h^{(n)} \rangle_{\L^2} $. 
	Hence, taking the sum over $ j=1, \ldots, n $
	\begin{align}
			&\frac{1}{2} \la |e_h^{(n)}|^2_{\L^2} + \sum_{j=1}^n |e_h^{(j)}-e_h^{(j-1)}|^2_{\L^2} \ra \nonumber \\
			&= \frac{1}{2} |e_h^{(0)}|^2_{\L^2} 
			+ \sum_{j=1}^n I_0^{(j)} 
			+ \sum_{j=1}^n \int_{\tau_{j-1}^{\epsilon,R}}^{\tau_j^{\epsilon,R}} \la I_1^{(j)} + I_2^{(j)} + I_3^{(j)} \ra \ ds 
			+ \sum_{j=1}^n \sum_{k=1}^\infty \int_{\tau_{j-1}^{\epsilon,R}}^{\tau_j^{\epsilon,R}} I_4^{(j)} \ ds 
			+ \sum_{j=1}^n I_5^{(j)}, \label{eq: 1/2|e_h^n|^2} 
	\end{align}
	where $ e_h^{(0)} = m_0 - m_h^{(0)} = m_0 - P_h m_0 $. 
	
	\underline{\bf Estimate $ e_h^{(0)} $.} 
	\begin{equation}\label{eq: eh0 est}
		\begin{aligned}
			|e_h^{(0)}|_{\L^2}^2 
			&= |m_0 - P_h m_0 |_{\L^2}^2 
			\leq ch^2 |m_0|_{\H^1}^2, \\
			|\nabla e_h^{(0)}|_{\L^2}^2 
			&= |\nabla (m_0 - P_h m_0) |_{\L^2}^2 
			\leq ch^2 |\nabla m_0|_{\H^1}^2. 	
		\end{aligned}
	\end{equation}
	
	\underline{\bf Estimate $ I_0 $.}
	\begin{align*}
		I_0^{(j)} 
		&= \langle e_h^{(j)}-e_h^{(j-1)}, e_h^{(j)}-P_h e_h^{(j)} \rangle_{\L^2} \\ 
		&= \langle e_h^{(j)}-e_h^{(j-1)}, m^\epsilon(\tau_j^{\epsilon,R}) - P_h m^\epsilon(\tau_j^{\epsilon,R}) \rangle_{\L^2} \\
		&\leq \delta |e_h^{(j)}-e_h^{(j-1)}|^2_{\L^2} + c\delta^{-1}|P_h m^\epsilon(\tau_j^{\epsilon,R})-m^\epsilon(\tau_j^{\epsilon,R})|^2_{\L^2} \\
		&= \delta |e_h^{(j)}-e_h^{(j-1)}|^2_{\L^2} + c\delta^{-1} h^2 |m^\epsilon(\tau_j^{\epsilon,R})|^2_{\H^1}.
	\end{align*}
	Thus, with the assumption $ Nh \leq C_2 $, 
	\begin{equation}\label{eq: I0 est}
		\begin{aligned}
			\sum_{j=1}^n I_0^{(j)} 
			&\leq \delta \sum_{j=1}^n |e_h^{(j)}-e_h^{(j-1)}|^2_{\L^2} + c\delta^{-1} h^2 \sum_{j=1}^n |m^\epsilon(\tau_j^{\epsilon,R})|^2_{\H^1} \\
			&\leq \delta \sum_{j=1}^n |e_h^{(j)}-e_h^{(j-1)}|^2_{\L^2} + c\delta^{-1} h |m^\epsilon|^2_{L^\infty(0,T; \H^1)} \\
			&= \delta \sum_{j=1}^n |e_h^{(j)}-e_h^{(j-1)}|^2_{\L^2} + h \varphi_0^{(n)} ,
		\end{aligned}
	\end{equation}
	for some constant $ c>0 $, where
	\begin{align*}
		\E \left[\max_{n \leq N} \varphi_0^{(n)} \right]
		&= \E \left[ c \delta^{-1} |m^\epsilon|^2_{L^\infty(0,T; \H^1)} \right]
		\leq c_0. 
	\end{align*}
	
	\underline{\bf Estimate $ I_1 $.}
	\begin{align*}
		I_1^{(j)} 
		&= - \epsilon \lb \Delta m^\epsilon(s) - \Delta m_h^{(j)}, \Delta P_h e_h^{(j)} \rb_{\L^2} - \lb \nabla m^\epsilon(s) - \nabla m_h^{(j)}, \nabla P_h e_h^{(j)} \rb_{\L^2} 
		=: I_{1a}^{(j)} + I_{1b}^{(j)}. 
	\end{align*}
	We first estimate $ I_{1a}^{(j)} $. 
	\begin{align*}
		\epsilon^{-1} I_{1a}^{(j)}
		&= - \lb \Delta m^\epsilon(s) - \Delta m_h^{(j)}, \Delta P_h e_h^{(j)} \rb_{\L^2} \\
		&= - \lb \Delta (m^\epsilon(s) - m_h^{(j)}), \Delta (P_h e_h^{(j)} - e_h^{(j)}) \rb_{\L^2} 
		- \lb \Delta (m^\epsilon(s) - m^\epsilon(\tau_j^{\epsilon,R})) + \Delta e_h^{(j)}, \Delta e_h^{(j)} \rb_{\L^2} \\
		&\leq ch |m^\epsilon(\tau_j^{\epsilon,R})|_{\H^3} |\Delta (m^\epsilon(s) - m_h^{(j)})|_{\L^2}
		- \lb \Delta (m^\epsilon(s) - m^\epsilon(\tau_j^{\epsilon,R})), \Delta e_h^{(j)} \rb_{\L^2} - |\Delta e_h^{(j)}|^2_{\L^2} \\
		&\leq ch |m^\epsilon|_{L^\infty(0,\tau^{\epsilon,R};\H^3)} \la |\Delta m^\epsilon(s)|_{\L^2} + |\Delta m_h^{(j)}|_{\L^2} \ra - (1-\delta) |\Delta e_h^{(j)}|^2_{\L^2} \\ 
		&\quad + c\delta^{-1} |\Delta (m^\epsilon(s) - m^\epsilon(\tau_j^{\epsilon,R}))|_{\L^2}^2 \\
		&\leq ch \la \frac{1}{2}|m^\epsilon|^2_{L^\infty(0,\tau^{\epsilon,R};\H^3)} + |\Delta m^\epsilon(s)|^2_{\L^2} + |\Delta m_h^{(j)}|^2_{\L^2} \ra - (1-\delta) |\Delta e_h^{(j)}|^2_{\L^2} \\ 
		&\quad + c\delta^{-1} |\Delta (m^\epsilon(s) - m^\epsilon(\tau_j^{\epsilon,R}))|_{\L^4}^2 \\
		&\leq ch \la \frac{3}{2}|m^\epsilon|^2_{L^\infty(0,\tau^{\epsilon,R};\H^3)} + |\Delta m_h^{(j)}|^2_{\L^2} \ra - (1-\delta) |\Delta e_h^{(j)}|^2_{\L^2} \\ 
		&\quad + c\delta^{-1} |\nabla (m^\epsilon(s) - m^\epsilon(\tau_j^{\epsilon,R}))|_{\H^2}^{\frac{3}{2}} |\nabla(m^\epsilon(s) - m^\epsilon(\tau_j^{\epsilon,R}))|_{\L^2}^\frac{1}{2} \\
		&\leq ch \la \frac{3}{2}|m^\epsilon|^2_{L^\infty(0,\tau^{\epsilon,R};\H^3)} + |\Delta m_h^{(j)}|^2_{\L^2} \ra - (1-\delta) |\Delta e_h^{(j)}|^2_{\L^2} \\ 
		&\quad + c\delta^{-1} |m^\epsilon|_{L^\infty(0,\tau^{\epsilon,R};\H^3)}^{\frac{3}{2}} |m^\epsilon|_{\mathcal{C}^\alpha([0,\tau^{\epsilon,R}]; \H^1)}^\frac{1}{2} |\tau_j^{\epsilon,R} -s|^\frac{\alpha}{2}.
	\end{align*}
	Then,
	\begin{equation}\label{eq: reg-I1a est}
		\begin{aligned}
			&\sum_{j=1}^n \int_{\tau_{j-1}^{\epsilon,R}}^{\tau_j^{\epsilon,R}} I_{1a}^{(j)} \ ds + \epsilon (1-\delta)\sum_{j=1}^n |\Delta e_h^{(j)}|^2_{\L^2} \Delta \tau_j^{\epsilon,R} \\
			&\leq 
			ch \la \frac{3}{2} \epsilon T |m^\epsilon|^2_{L^\infty(0,\tau^{\epsilon,R};\H^3)} + \epsilon \sum_{j=1}^n |\Delta m_h^{(j)}|^2_{\L^2} \Delta \tau_j^{\epsilon,R} \ra \\
			&\quad + c\epsilon\delta^{-1} |m^\epsilon|_{L^\infty(0,\tau^{\epsilon,R};\H^3)}^{\frac{3}{2}} |m^\epsilon|^\frac{1}{2}_{\mathcal{C}^\alpha([0,\tau^{\epsilon,R}]; \H^1)} (\Delta t)^\frac{\alpha}{2} T. 
		\end{aligned}
	\end{equation}
	Then, we estimate $ I_{1b}^{(j)} $. 
	\begin{align*}
		I_{1b}^{(j)}
		&= - \lb \nabla m^\epsilon(s) - \nabla m_h^{(j)}, \nabla P_h e_h^{(j)} \rb_{\L^2} \\
		&= - \lb \nabla m^\epsilon(s) - \nabla m^\epsilon(\tau_j^{\epsilon,R}) + \nabla e_h^{(j)}, \nabla P_h e_h^{(j)} \rb_{\L^2} \\
		&= - \lb \nabla m^\epsilon(s) - \nabla m^\epsilon(\tau_j^{\epsilon,R}), \nabla P_h e_h^{(j)} \rb_{\L^2} - \lb \nabla e_h^{(j)}, \nabla (P_h m^\epsilon(\tau_j^{\epsilon,R}) - m^\epsilon(\tau_j^{\epsilon,R})) \rb_{\L^2} - |\nabla e_h^{(j)}|^2_{\L^2}.
	\end{align*}
	Then,
	\begin{align*} 
		I_{1b}^{(j)} + |\nabla e_h^{(j)}|^2_{\L^2} 
		&\leq c\delta^{-1} \la |\nabla m^\epsilon(\tau_j^{\epsilon,R}) - \nabla m^\epsilon(s)|^2_{\L^2} + |\nabla (P_h m^\epsilon(\tau_j^{\epsilon,R}) - m^\epsilon(\tau_j^{\epsilon,R}))|^2_{\L^2} \ra \\
		&\quad + \delta \la |\nabla P_h e_h^{(j)}|^2_{\L^2} +|\nabla e_h^{(j)}|^2_{\L^2} \ra \\
		&\leq c\delta^{-1} |\nabla m^\epsilon(\tau_j^{\epsilon,R}) - \nabla m^\epsilon(s)|^2_{\L^2} + c\delta^{-1} h^2 |m^\epsilon(\tau_j^{\epsilon,R})|^2_{\H^2}
		+ 2\delta |\nabla e_h^{(j)}|^2_{\L^2}. 
	\end{align*}
	With $ \Delta \tau_j^{\epsilon,R} \leq \Delta t $,
	\begin{equation}\label{eq: I1b est}
		\begin{aligned}
			&\sum_{j=1}^n \int_{\tau_{j-1}^{\epsilon,R}}^{\tau_j^{\epsilon,R}} I_{1b}^{(j)} \ ds + (1-2\delta)\sum_{j=1}^n |\nabla e_h^{(j)}|^2_{\L^2} \Delta \tau_j^{\epsilon,R} \\
			&\leq 
			c\delta^{-1} |\nabla m^\epsilon|_{\mathcal{C}^\alpha([0,\tau^{\epsilon,R}]; \L^2)}^2 (\Delta t)^{2 \alpha} T 
			+ c\delta^{-1} h^2 |m^\epsilon|^2_{L^\infty(0,\tau^{\epsilon,R};\H^2)} T.
		\end{aligned}
	\end{equation}
	Combining $ I_{1a}^{(j)} $ and $ I_{1b}^{(j)} $ estimates, we have
	\begin{equation}\label{eq: I1 est}
		\begin{aligned}
			&\sum_{j=1}^n \int_{\tau_{j-1}^{\epsilon,R}}^{\tau_j^{\epsilon,R}} I_1^{(j)} \ ds 
			+ (1-2\delta)\sum_{j=1}^n |\nabla e_h^{(j)}|^2_{\L^2} \Delta \tau_j^{\epsilon,R} 
			+ \epsilon (1-\delta)\sum_{j=1}^n |\Delta e_h^{(j)}|^2_{\L^2} \Delta \tau_j^{\epsilon,R} \\
			&\leq 
			c \epsilon\delta^{-1} |m^\epsilon|_{L^\infty(0,\tau^{\epsilon,R};\H^3)}^{\frac{3}{2}} |m^\epsilon|^\frac{1}{2}_{\mathcal{C}^\alpha([0,\tau^{\epsilon,R}]; \H^1)} (\Delta t)^\frac{\alpha}{2}
			+ c \delta^{-1} |\nabla m^\epsilon|_{\mathcal{C}^\alpha([0,\tau^{\epsilon,R}]; \L^2)}^2 (\Delta t)^{2\alpha} \\
			&\quad + c h \la \epsilon |m^\epsilon|^2_{L^\infty(0,\tau^{\epsilon,R};\H^3)} + \epsilon \sum_{j=1}^n |\Delta m_h^{(j)}|^2_{\L^2} \Delta \tau_j^{\epsilon,R} + \delta^{-1} h |m^\epsilon|^2_{L^\infty(0,\tau^{\epsilon,R};\H^2)} \ra \\
			&= (\Delta t)^\frac{\alpha}{2}\psi_{1a}^{(n)} + (\Delta t)^{2\alpha}\psi_{1b}^{(n)} + h \varphi_1^{(n)} ,
		\end{aligned}
	\end{equation}
	where 
	\begin{align*}
		\E \left[ \max_{n \leq N} \psi_{1a}^{(n)} \right] \leq c_1 \epsilon^{-\frac{5}{3}}, \quad 
		\E \left[ \max_{n \leq N} \psi_{1b}^{(n)} \right] \leq c_1 \epsilon^{-\frac{8}{3}}, \quad
		\E \left[ \max_{n \leq N} \varphi_1^{(n)} \right] \leq c_1 \epsilon^{-\frac{5}{3}}. 
	\end{align*}
	
	\underline{\bf Estimate $ I_2 $.}
	\begin{align*}
		I_2^{(j)} 
		&= \lb m_h^{(j-1)} \times \nabla m_h^{(j)} - m^\epsilon(s) \times \nabla m^\epsilon(s), \nabla P_h e_h^{(j)} \rb_{\L^2} \\
		&= \lb \la m^\epsilon(\tau_{j-1}^{\epsilon,R}) - m^\epsilon(s) \ra \times \nabla m^\epsilon(s), \nabla P_h e_h^{(j)} \rb_{\L^2} \\
		&\quad + \lb m_h^{(j-1)} \times \la \nabla m^\epsilon(\tau_j^{\epsilon,R}) - \nabla m^\epsilon(s) \ra, \nabla P_h e_h^{(j)} \rb_{\L^2} \\
		&\quad - \lb e_h^{(j-1)} \times \nabla m^\epsilon(s), \nabla P_h e_h^{(j)} \rb_{\L^2} 
		- \lb m_h^{(j-1)} \times \nabla e_h^{(j)}, \nabla P_h e_h^{(j)} \rb_{\L^2} \\
		&= I_{2a}^{(j)} + I_{2b}^{(j)} + I_{2c}^{(j)} + I_{2d}^{(j)}.
	\end{align*}
	Then, with $ \H^1 \hookrightarrow \L^4 $, 
	\begin{align*}
		I_{2a}^{(j)} 
		&= \lb \la m^\epsilon(\tau_{j-1}^{\epsilon,R}) - m^\epsilon(s) \ra \times \nabla m^\epsilon(s), \nabla P_h e_h^{(j)} \rb_{\L^2} \\
		&\leq c\delta^{-1} \left| \la m^\epsilon(\tau_{j-1}^{\epsilon,R}) - m^\epsilon(s) \ra \times \nabla m^\epsilon(s) \right|_{\L^2}^2 
		+ \delta |\nabla P_h e_h^{(j)}|^2_{\L^2} \\
		&\leq c\delta^{-1} \la |m^\epsilon(\tau_{j-1}^{\epsilon,R}) - m^\epsilon(s)|_{\H^1}^2 |\nabla m^\epsilon(s)|_{\L^4}^2 \ra 
		+ \delta|\nabla e_h^{(j)}|^2_{\L^2}, \\
		\sum_{j=1}^n \int_{\tau_{j-1}^{\epsilon,R}}^{\tau_j^{\epsilon,R}} I_{2a}^{(j)} \ ds 
		&\leq c\delta^{-1} \sum_{j=1}^n \int_{\tau_{j-1}^{\epsilon,R}}^{\tau_j^{\epsilon,R}} |m^\epsilon(s)-m^\epsilon(\tau_{j-1}^{\epsilon,R})|_{\H^1}^2 |\nabla m^\epsilon(s)|_{\L^4}^2 \ ds
		+ \delta \sum_{j=1}^n |\nabla e_h^{(j)}|^2_{\L^2} \Delta \tau_j^{\epsilon,R} \\
		&\leq c\delta^{-1} |m^\epsilon|_{\mathcal{C}^\alpha([0,\tau^{\epsilon,R}]; \H^1)}^2 (\Delta t)^{2\alpha} \int_0^{\tau_n^{\epsilon,R}} |m^\epsilon(s)|_{\H^2}^2 \ ds
		+ \delta \sum_{j=1}^n |\nabla e_h^{(j)}|^2_{\L^2} \Delta \tau_j^{\epsilon,R}. 
	\end{align*}	
	Also, using $ \Delta \tau_j^{\epsilon,R} \leq \Delta \tau_{j-1}^R $, we have
	\begin{align*}
		I_{2b}^{(j)} 
		&= \lb m_h^{(j-1)} \times \la \nabla m^\epsilon(\tau_j^{\epsilon,R})- \nabla m^\epsilon(s) \ra, \nabla P_h e_h^{(j)} \rb_{\L^2} \\
		&= \lb m_h^{(j-1)} \times \la \nabla m^\epsilon(\tau_j^{\epsilon,R})- \nabla m^\epsilon(s) \ra, \nabla P_h e_h^{(j)} \rb_{\L^2} \\
		&\leq \delta^{-1} |m_h^{(j-1)}|^2_{\L^4} |\nabla m^\epsilon(\tau_j^{\epsilon,R})- \nabla m^\epsilon(s)|_{\L^4}^2 + \delta |\nabla e_h^{(j)}|^2_{\L^2} \\
		&\leq c\delta^{-1} |m_h^{(j-1)}|_{\H^1}^2 |\nabla m^\epsilon(\tau_j^{\epsilon,R})- \nabla m^\epsilon(s)|_{\H^1} |\nabla m^\epsilon(\tau_j^{\epsilon,R})- \nabla m^\epsilon(s)|_{\L^2} 
		+ \delta |\nabla e_h^{(j)}|^2_{\L^2} , \\
		\sum_{j=1}^n \int_{\tau_{j-1}^{\epsilon,R}}^{\tau_j^{\epsilon,R}} I_{2b}^{(j)} \ ds
		&\leq c\delta^{-1} |\nabla m^\epsilon|_{\mathcal{C}^\alpha([0,\tau^{\epsilon,R}]; \L^2)} (\Delta t)^{\alpha} |m^\epsilon|_{L^\infty(0,\tau^{\epsilon,R};\H^2)} \sum_{j=1}^n |m_h^{(j-1)}|_{\H^1}^2 \Delta \tau_{j-1}^{\epsilon,R} \\
		&\quad + \delta \sum_{j=1}^n |\nabla e_h^{(j)}|^2_{\L^2} \Delta \tau_j^{\epsilon,R}.
	\end{align*}	
	Similarly,
	\begin{align*}
		I_{2c}^{(j)}
		&= - \lb e_h^{(j-1)} \times \nabla m^\epsilon, \nabla P_h e_h^{(j)} \rb_{\L^2} \\
		&= \lb e_h^{(j-1)} \times \nabla P_h e_h^{(j)}, \nabla m^\epsilon \rb_{\L^2} \\
		&= -\lb \nabla e_h^{(j-1)} \times \nabla P_h e_h^{(j)} + e_h^{(j-1)} \times \Delta P_h e_h^{(j)}, m^\epsilon \rb_{\L^2} \\
		&\leq 
		|\nabla e_h^{(j-1)}|_{\L^2} |\nabla P_h e_h^{(j)}|_{\L^4} |m^\epsilon|_{\L^4} 
		+ |\Delta e_h^{(j)}|_{\L^2} |e_h^{(j-1)}|_{\L^4} |m^\epsilon|_{\L^4} \\
		&\leq 
		c|\nabla e_h^{(j-1)}|_{\L^2} \la |\Delta e_h^{(j)}|_{\L^2}^\frac{3}{4} |e_h^{(j)}|_{\L^2}^\frac{1}{4} + |e_h^{(j)}|_{\L^2} \ra |m^\epsilon|_{\H^1} \\
		&\quad + c|\Delta e_h^{(j)}|_{\L^2} \la |\Delta e_h^{(j-1)}|_{\L^2}^\frac{1}{4} |e_h^{(j-1)}|_{\L^2}^\frac{3}{4} + |e_h^{(j-1)}|_{\L^2} \ra |m^\epsilon|_{\H^1} \\
		&\leq 
		\delta |\nabla e_h^{(j-1)}|_{\L^2}^2 + c \delta^{-1} \la |\Delta e_h^{(j)}|_{\L^2}^\frac{3}{2} |e_h^{(j)}|_{\L^2}^\frac{1}{2} + |e_h^{(j)}|_{\L^2}^2 \ra |m^\epsilon|_{\H^1}^2 \\
		&\quad + \delta \epsilon |\Delta e_h^{(j)}|^2_{\L^2} + c (\delta \epsilon)^{-1} \la |\Delta e_h^{(j-1)}|_{\L^2}^\frac{1}{2} |e_h^{(j-1)}|_{\L^2}^\frac{3}{2} + |e_h^{(j-1)}|_{\L^2}^2 \ra |m^\epsilon|_{\H^1}^2 \\
		&\leq 
		\delta |\nabla e_h^{(j-1)}|_{\L^2}^2 
		+ \frac{3}{4} \delta^{\frac{1}{3}} \epsilon |\Delta e_h^{(j)}|_{\L^2}^2 
		+ c \delta^{-5} \epsilon^{-3} |e_h^{(j)}|_{\L^2}^2 |m^\epsilon|_{\H^1}^8 
		+ c \delta^{-1} |e_h^{(j)}|_{\L^2}^2 |m^\epsilon|_{\H^1}^2 \\
		&\quad + \delta \epsilon |\Delta e_h^{(j)}|^2_{\L^2} 		
		+ \frac{1}{4} \delta \epsilon |\Delta e_h^{(j-1)}|_{\L^2}^2 
		+ c \delta^{-\frac{5}{3}} \epsilon^{-\frac{5}{3}} |e_h^{(j-1)}|_{\L^2}^2 |m^\epsilon|_{\H^1}^{\frac{8}{3}} 
		+ c \delta^{-1} \epsilon^{-1} |e_h^{(j-1)}|_{\L^2}^2 |m^\epsilon|_{\H^1}^2 \\
		&\leq 
		\delta |\nabla e_h^{(j-1)}|_{\L^2}^2 + \la \frac{3}{4} \delta^\frac{1}{3} + \delta \ra \epsilon |\Delta e_h^{(j)}|^2_{\L^2} + \frac{1}{4} \delta \epsilon |\Delta e_h^{(j-1)}|_{\L^2}^2 \\
		&\quad + c \delta^{-5} \epsilon^{-3} |e_h^{(j)}|_{\L^2}^2 |m^\epsilon|_{\H^1}^8 + c \delta^{-\frac{5}{3}} \epsilon^{-\frac{5}{3}} |e_h^{(j-1)}|_{\L^2}^2 |m^\epsilon|_{\H^1}^{\frac{8}{3}} \\
		&\quad + c \delta^{-1} \la |e_h^{(j)}|_{\L^2}^2 + \epsilon^{-1}|e_h^{(j-1)}|_{\L^2}^2 \ra |m^\epsilon|_{\H^1}^2,  
	\end{align*} 	
	and then with $ \delta, \epsilon \in (0,1) $,
	\begin{align*}
		\sum_{j=1}^n \int_{\tau_{j-1}^{\epsilon,R}}^{\tau_j^{\epsilon,R}} I_{2c}^{(j)} \ ds
		&\leq 
		\delta \sum_{j=1}^n |\nabla e_h^{(j)}|_{\L^2}^2 \Delta \tau_{j}^{\epsilon,R} 
		+ \la \frac{3}{4} \delta^\frac{1}{3} + \frac{5}{4} \delta \ra \epsilon \sum_{j=1}^n |\Delta e_h^{(j)}|^2_{\L^2} \Delta \tau_j^{\epsilon,R} \\
		&\quad + c \delta^{-5} \epsilon^{-3} \sum_{j=1}^n |e_h^{(j)}|_{\L^2}^2 \int_{\tau_{j-1}^{\epsilon,R}}^{\tau_{j}^{\epsilon,R}} \la |m^\epsilon(s)|_{\H^1}^8 + |m^\epsilon(s)|_{\H^1}^{\frac{8}{3}} + |m^\epsilon(s)|_{\H^1}^2 \ra \ ds \\
		&\quad + \frac{1}{4} \delta \epsilon |\Delta e_h^{(0)}|_{\L^2}^2 \Delta t + \delta |\nabla e_h^{(0)}|_{\L^2}^2 \Delta t \\
		&\quad + c \delta^{-\frac{5}{3}} \epsilon^{-\frac{5}{3}} |e_h^{(0)}|_{\L^2}^2 \la |m^\epsilon|_{L^\infty(0,T;\H^1)}^{\frac{8}{3}} + |m^\epsilon|_{L^\infty(0,T;\H^1)}^2 \ra.
	\end{align*}
	For $ I_{2d}^{(j)} $,  
	\begin{align*}
		I_{2d}^{(j)}
		&= - \lb m_h^{(j-1)} \times \nabla e_h^{(j)}, \nabla P_h e_h^{(j)} \rb_{\L^2} \\
		&= - \lb m_h^{(j-1)} \times \nabla e_h^{(j)}, \nabla \la P_h e_h^{(j)} - e_h^{(j)} \ra \rb_{\L^2} \\
		&= \lb m_h^{(j-1)} \times \nabla \la P_h m^\epsilon(\tau_j^{\epsilon,R}) - m^\epsilon(\tau_j^{\epsilon,R}) \ra, \nabla e_h^{(j)} \rb_{\L^2} \\
		&\leq c\delta^{-1} |m_h^{(j-1)}|^2_{\L^4} \left|\nabla \la P_h m^\epsilon(\tau_j^{\epsilon,R}) - m^\epsilon(\tau_j^{\epsilon,R}) \ra \right|^2_{\L^4} + \delta |\nabla e_h^{(j)}|^2_{\L^2} \\
		&\leq c\delta^{-1} |m_h^{(j-1)}|^2_{\H^1} \left|\nabla \la P_h m^\epsilon(\tau_j^{\epsilon,R}) - m^\epsilon(\tau_j^{\epsilon,R}) \ra \right|_{\L^2} \left|\nabla \la P_h m^\epsilon(\tau_j^{\epsilon,R}) - m^\epsilon(\tau_j^{\epsilon,R}) \ra \right|_{\H^1} 
		+ \delta |\nabla e_h^{(j)}|^2_{\L^2} \\
		&\leq c \delta^{-1} |m_h^{(j-1)}|^2_{\H^1} h |m^\epsilon(\tau_j^{\epsilon,R})|^2_{\H^2} + \delta |\nabla e_h^{(j)}|^2_{\L^2}, 
	\end{align*}
	and
	\begin{align*}
		\sum_{j=1}^n \int_{\tau_{j-1}^{\epsilon,R}}^{\tau_j^{\epsilon,R}} I_{2d}^{(j)} \ ds
		&\leq c \delta^{-1} h \sum_{j=1}^n |m_h^{(j-1)}|^2_{\H^1} |m^\epsilon(\tau_j^{\epsilon,R})|^2_{\H^2} \Delta \tau_j^{\epsilon,R} + \delta \sum_{j=1}^n |\nabla e_h^{(j)}|^2_{\L^2} \Delta \tau_j^{\epsilon,R} \\
		&\leq c \delta^{-1} h |m^\epsilon|^2_{L^\infty(0,\tau^{\epsilon,R};\H^2)} \sum_{j=1}^n |m_h^{(j-1)}|^2_{\H^1} \Delta \tau_{j-1}^{\epsilon,R} + \delta \sum_{j=1}^n |\nabla e_h^{(j)}|^2_{\L^2} \Delta \tau_j^{\epsilon,R}.
	\end{align*}
	Hence, 
	\begin{align}\label{eq: I2 est}
		\begin{aligned}
			&\sum_{j=1}^n \int_{\tau_{j-1}^{\epsilon,R}}^{\tau_j^{\epsilon,R}} I_2^{(j)} \ ds \\
			&\leq 
			c (\delta \epsilon)^{-\frac{5}{3}} \la |m^\epsilon|_{L^\infty(0,T;\H^1)}^{\frac{8}{3}} + |m^\epsilon|_{L^\infty(0,T;\H^1)}^2 \ra |e_h^{(0)}|_{\L^2}^2 
			+ \delta \la |\nabla e_h^{(0)}|_{\L^2}^2 + \epsilon |\Delta e_h^{(0)}|_{\L^2}^2 \ra \Delta t \\
			&\quad + c \delta^{-1} (\Delta t)^{2\alpha}  |m^\epsilon|_{\mathcal{C}^\alpha([0,\tau^{\epsilon,R}]; \H^1)}^2 \int_0^{\tau_n^{\epsilon,R}} |m^\epsilon(s)|_{\H^2}^2 \ ds \\
			&\quad + c \delta^{-1} (\Delta t)^{\alpha} |m^\epsilon|_{\mathcal{C}^\alpha([0,\tau^{\epsilon,R}]; \H^1)} |m^\epsilon|_{L^\infty(0,\tau^{\epsilon,R};\H^2)} \sum_{j=1}^n |m_h^{(j-1)}|_{\H^1}^2 \Delta \tau_{j-1}^{\epsilon,R} \\
			&\quad + c \delta^{-1} h |m^\epsilon|^2_{L^\infty(0,\tau^{\epsilon,R};\H^2)} \sum_{j=1}^n |m_h^{(j-1)}|^2_{\H^1} \Delta \tau_{j-1}^{\epsilon,R} \\
			&\quad + 4\delta \sum_{j=1}^n |\nabla e_h^{(j)}|^2_{\L^2} \Delta \tau_j^{\epsilon,R} + 2\delta^\frac{1}{3} \epsilon \sum_{j=1}^n |\Delta e_h^{(j)}|^2_{\L^2} \Delta \tau_j^{\epsilon,R} \\
			&\quad + c \delta^{-5} \epsilon^{-3} \sum_{j=1}^n |e_h^{(j)}|_{\L^2}^2 \int_{\tau_{j-1}^{\epsilon,R}}^{\tau_{j}^{\epsilon,R}} \la |m^\epsilon(s)|_{\H^1}^8 + |m^\epsilon(s)|_{\H^1}^{\frac{8}{3}} + |m^\epsilon(s)|_{\H^1}^2 \ra \ ds \\
			&= \eta_2 + (\Delta t)^{2\alpha} \psi_{2a}^{(n)} + (\Delta t)^\alpha \psi_{2b}^{(n)} + h \varphi_2^{(n)} \\
			&\quad + 4\delta \sum_{j=1}^n |\nabla e_h^{(j)}|^2_{\L^2} \Delta \tau_j^{\epsilon,R} + 2\delta^\frac{1}{3} \epsilon \sum_{j=1}^n |\Delta e_h^{(j)}|^2_{\L^2} \Delta \tau_j^{\epsilon,R} \\
			&\quad + c \delta^{-5} \epsilon^{-3} \sum_{j=1}^n |e_h^{(j)}|_{\L^2}^2 \int_{\tau_{j-1}^{\epsilon,R}}^{\tau_{j}^{\epsilon,R}} \la |m^\epsilon(s)|_{\H^1}^8 + |m^\epsilon(s)|_{\H^1}^{\frac{8}{3}} + |m^\epsilon(s)|_{\H^1}^2 \ra \ ds ,
		\end{aligned}
	\end{align}	
	where 
	\begin{align*}
		\E \left[\eta_2 \right] 
		&\leq c_2 \epsilon^{-\frac{5}{3}} \E \left[|e_h^{(0)}|_{\L^2}^4 \right]^\frac{1}{2} + c_2 \E \left[ |\nabla e_h^{(0)}|_{\L^2}^2 + \epsilon |\Delta e_h^{(0)}|_{\L^2}^2 \right] \Delta t \\
		&\leq c_2 \epsilon^{-\frac{5}{3}} h^2 |m_0|_{\H^1}^2 + c_2 \Delta t \la |\nabla m_0|_{\H^1}^2 + \epsilon^2 |\Delta m_0|_{\H^1}^2 \ra \\
		&\leq c_2 \epsilon^{-\frac{5}{3}} (h^2 + \Delta t) |m_0|_{\H^3}^2, 
	\end{align*}
	and
	\begin{align*}
		\E \left[ \max_{n \leq N} \psi_{2a}^{(n)} \right] \leq c_2 \epsilon^{-\frac{8}{3}} \quad
		\E \left[ \max_{n \leq N} \psi_{2b}^{(n)} \right] \leq c_2 \epsilon^{-\frac{13}{6}}\quad
		\E \left[ \max_{n \leq N} \varphi_2^{(n)} \right] \leq c_2 \epsilon^{-\frac{5}{3}}.
	\end{align*}

	\underline{\bf Estimate $ I_3 $.}
	\begin{align*}
		I_3^{(j)}
		&= -\lb (1 + |m^\epsilon(s)|^2) m^\epsilon(s) - (1+|m_h^{(j-1)}|^2) m_h^{(j)}, P_h e_h^{(j)} \rb_{\L^2} \\
		&= \lb m_h^{(j)}-m^\epsilon(s), P_h e_h^{(j)} \rb_{\L^2} + \lb |m_h^{(j-1)}|^2 (m_h^{(j)}-m^\epsilon(s)), P_h e_h^{(j)} \rb_{\L^2} \\
		&\quad + \lb (|m_h^{(j-1)}|^2 -|m^\epsilon(s)|^2)m^\epsilon(s), P_h e_h^{(j)} \rb_{\L^2} \\
		&= \lb m^\epsilon(\tau_j^{\epsilon,R}) - m^\epsilon(s) - e_h^{(j)}, P_h e_h^{(j)} \rb_{\L^2} \\
		&\quad + \lb |m_h^{(j-1)}|^2 (m_h^{(j)} - m^\epsilon(s)), P_h e_h^{(j)} \rb_{\L^2} \\
		&\quad + \lb \langle m_h^{(j-1)} + m^\epsilon(s), m_h^{(j-1)} - m^\epsilon(\tau_j^{\epsilon,R}) \rangle m^\epsilon(s), P_he_h^{(j)} - e_h^{(j)} \rb_{\L^2} \\
		&\quad + \lb \langle m_h^{(j-1)} + m^\epsilon(s), (m_h^{(j-1)} - m_h^{(j)}) - e_h^{(j)} \rangle m^\epsilon(s), e_h^{(j)} \rb_{\L^2} \\
		&\quad + \lb \langle m_h^{(j-1)} + m^\epsilon(s), m^\epsilon(\tau_j^{\epsilon,R})-m^\epsilon(s) \rangle m^\epsilon(s), P_h e_h^{(j)} \rb_{\L^2} \\
		&= I_{3a}^{(j)} + I_{3b}^{(j)} + I_{3c}^{(j)} + I_{3d}^{(j)} + I_{3e}^{(j)}.
	\end{align*}
	We first estimate $ I_{3a}^{(j)} $.
	\begin{align*}
		I_{3a}^{(j)} 
		&= \lb m^\epsilon(\tau_j^{\epsilon,R}) - m^\epsilon(s) - e_h^{(j)}, P_h e_h^{(j)} \rb_{\L^2} \\
		&\leq \frac{1}{2} |m^\epsilon(s)-m^\epsilon(\tau_j^{\epsilon,R})|^2_{\L^2} + \frac{3}{2} |e_h^{(j)}|^2_{\L^2} \\
		&\leq \frac{1}{2} |m^\epsilon|_{\mathcal{C}^\alpha([0,\tau^{\epsilon,R}]; \L^2)}^2 |\tau_j^{\epsilon,R}-s|^{2\alpha} + \frac{3}{2} |e_h^{(j)}|^2_{\L^2}, \\
		\sum_{j=1}^n \int_{\tau_{j-1}^{\epsilon,R}}^{\tau_j^{\epsilon,R}} I_{3a}^{(j)} \ ds
		&\leq \frac{1}{2} |m^\epsilon|_{\mathcal{C}^\alpha([0,\tau^{\epsilon,R}]; \L^2)}^2 (\Delta t)^{2\alpha} T + \frac{3}{2} \sum_{j=1}^n |e_h^{(j)}|^2_{\L^2} \Delta \tau_j^{\epsilon,R} .
	\end{align*}
	For $ I_{3b}^{(j)} $, since $ \H^1 \hookrightarrow \L^8 \hookrightarrow \L^4 $,
	\begin{align*}
		I_{3b}^{(j)}
		&= \lb |m_h^{(j-1)}|^2 (m^\epsilon(\tau_j^{\epsilon,R}) - m^\epsilon(s) - e_h^{(j)}), P_h e_h^{(j)} \rb_{\L^2} \\
		&= - \lb |m_h^{(j-1)}|^2 e_h^{(j)}, e_h^{(j)} \rb_{\L^2} + \lb |m_h^{(j-1)}|^2 (m^\epsilon(\tau_j^{\epsilon,R}) - m^\epsilon(s)), e_h^{(j)} \rb_{\L^2} \\
		&\quad + \lb |m_h^{(j-1)}|^2 (m_h^{(j)} - m^\epsilon(s)),  P_he_h^{(j)} - e_h^{(j)} \rb_{\L^2} \\
		&\leq - \int_D |m_h^{(j-1)}|^2 |e_h^{(j)}|^2 \ dx \\
		&\quad + \delta \int_D |m_h^{(j-1)}|^2 |e_h^{(j)}|^2 \ dx + c\delta^{-1} \int_D |m_h^{(j-1)}|^2 |m^\epsilon(\tau_j^{\epsilon,R})-m^\epsilon(s)|^2 \ dx \\
		&\quad + |m_h^{(j-1)}|^2_{\L^8} |m_h^{(j)}-m^\epsilon(s)|_{\L^2} |P_h e_h^{(j)} - e_h^{(j)}|_{\L^4} \\
		&\leq - (1-\delta)\int_D |m_h^{(j-1)}|^2 |e_h^{(j)}|^2 \ dx + c\delta^{-1} |m_h^{(j-1)}|_{\L^4}^2 |m^\epsilon(\tau_j^{\epsilon,R})-m^\epsilon(s)|_{\L^4}^2 \\
		&\quad + |m_h^{(j-1)}|^2_{\H^1} \la |m_h^{(j)}|_{\L^2} + |m^\epsilon(s)|_{\L^2} \ra |P_h e_h^{(j)} - e_h^{(j)}|_{\H^1} \\
		&\leq - (1-\delta)\int_D |m_h^{(j-1)}|^2 |e_h^{(j)}|^2 \ dx + c\delta^{-1} |m_h^{(j-1)}|_{\H^1}^2 |m^\epsilon|_{\mathcal{C}^\alpha([0,\tau^{\epsilon,R}];\H^1)}^2 |\tau_j^{\epsilon,R}-s|^{2\alpha} \\
		&\quad + ch |m_h^{(j-1)}|^2_{\H^1} \la |m_h^{(j)}|_{\L^2} + |m^\epsilon(s)|_{\L^2} \ra |m^\epsilon(\tau_j^{\epsilon,R})|_{\H^2} \\
		\sum_{j=1}^n \int_{\tau_{j-1}^{\epsilon,R}}^{\tau_j^{\epsilon,R}} I_{3b}^{(j)} \ ds
		&\leq 
		- (1-\delta)\sum_{j=1}^n \int_D |m_h^{(j-1)}|^2 |e_h^{(j)}|^2 \ dx \ \Delta \tau_j^{\epsilon,R} \\
		&\quad + c\delta^{-1} |m^\epsilon|^2_{\mathcal{C}^\alpha([0,\tau^{\epsilon,R}]; \H^1)} (\Delta t)^{2\alpha} \la \sum_{j=1}^n |m_h^{(j-1)}|_{\H^1}^2 \Delta \tau_{j-1}^{\epsilon,R} \ra \\
		&\quad + ch \la \max_{j \leq n} |m_h^{(j)}|_{\L^2} + |m^\epsilon|_{L^\infty(0,T;\L^2)} \ra |m^\epsilon|_{L^\infty(0,\tau^{\epsilon,R};\H^2)} \sum_{j=1}^n |m_h^{(j-1)}|^2_{\H^1} \Delta \tau_{j-1}^{\epsilon,R}.   
	\end{align*}
	By the property of $ P_h $ and continuous embeddings $ \H^2 \hookrightarrow \L^\infty $ and $ \H^1 \hookrightarrow \L^4 $,  
	\begin{align*}
		I_{3c}^{(j)}
		&= \lb \langle m_h^{(j-1)} + m^\epsilon(s), m_h^{(j-1)} - m^\epsilon(\tau_j^{\epsilon,R}) \rangle m^\epsilon(s), P_he_h^{(j)} - e_h^{(j)} \rb_{\L^2} \\
		&\leq |m^\epsilon(s)|_{\L^\infty} |m_h^{(j-1)} + m^\epsilon(s)|_{\L^4} |m_h^{(j-1)} - m^\epsilon(\tau_j^{\epsilon,R})|_{\L^4} |P_he_h^{(j)} - e_h^{(j)}|_{\L^2} \\
		&\leq c h |m^\epsilon(s)|_{\H^2} \la |m_h^{(j-1)}|^2_{\H^1} + |m^\epsilon|^2_{L^\infty(0,T; \H^1)} \ra |m^\epsilon(\tau_j^{\epsilon,R})|_{\H^1} \\
		\sum_{j=1}^n \int_{\tau_{j-1}^{\epsilon,R}}^{\tau_j^{\epsilon,R}} I_{3c}^{(j)} \ ds
		&\leq c h |m^\epsilon|_{L^\infty(0,T; \H^1)} |m^\epsilon|_{L^\infty(0,\tau^{\epsilon,R}; \H^2)} \sum_{j=1}^n |m_h^{(j-1)}|_{\H^1}^2 \Delta \tau_{j-1}^{\epsilon,R} \\
		&\quad + c h \la |m^\epsilon|_{L^\infty(0,T; \H^1)}^6 + \int_0^{\tau^{\epsilon,R}} |m^\epsilon(s)|^2_{\H^2} \ ds \ra. 
	\end{align*}
	Next, for $ I_{3d}^{(j)} $, we have
	\begin{align*}
		I_{3d}^{(j)} 
		&= \lb \langle m_h^{(j-1)} + m^\epsilon(s), (m_h^{(j-1)} - m_h^{(j)}) - e_h^{(j)} \rangle m^\epsilon(s), e_h^{(j)} \rb_{\L^2} \\
		&= \lb \langle m_h^{(j-1)}, m_h^{(j-1)} - m_h^{(j)} \rangle m^\epsilon(s), e_h^{(j)} \rb_{\L^2} 
		+ \lb \langle m^\epsilon(s), m_h^{(j-1)} - m_h^{(j)} \rangle m^\epsilon(s), e_h^{(j)} \rb_{\L^2} \\
		&\quad - \lb \langle m_h^{(j-1)} + m^\epsilon(s), e_h^{(j)} \rangle m^\epsilon(s), e_h^{(j)} \rb_{\L^2} \\
		&\leq \delta \int_D |m_h^{(j-1)}|^2 |e_h^{(j)}|^2 \ dx + c\delta^{-1} \int_D |m_h^{(j-1)} - m_h^{(j)}|^2 |m^\epsilon(s)|^2 \ dx \\  
		&\quad + |m_h^{j-1}-m_h^{(j)}|_{\L^2} |m^\epsilon(s)|_{\L^8}^2 |e_h^{(j)}|_{\L^4} \\
		&\quad - \int_D \lb m_h^{(j-1)},e_h^{(j)} \rb \lb m^\epsilon(s), e_h^{(j)} \rb \ dx - \int_D \lb m^\epsilon(s), e_h^{(j)} \rb^2 \ dx \\
		&\leq \delta \int_D |m_h^{(j-1)}|^2 |e_h^{(j)}|^2 \ dx + c\delta^{-1} |m^\epsilon(s)|_{\L^\infty}^2 |m_h^{(j-1)} - m_h^{(j)}|^2_{\L^2} \\ 
		&\quad + |m_h^{j-1}-m_h^{(j)}|_{\L^2} |m^\epsilon(s)|_{\H^1}^2 |e_h^{(j)}|_{\H^1} \\
		&\quad - \int_D \frac{1}{\sqrt{2}} \lb m_h^{(j-1)},e_h^{(j)} \rb \sqrt{2}\lb m^\epsilon(s), e_h^{(j)} \rb \ dx - \int_D \lb m^\epsilon(s), e_h^{(j)} \rb^2 \ dx \\
		&\leq c\delta^{-1} |m^\epsilon(s)|_{\L^\infty}^2 |m_h^{(j-1)} - m_h^{(j)}|^2_{\L^2} \\
		&\quad + \delta |e_h^{(j)}|_{\H^1}^2 + c \delta^{-1} |m_h^{j-1}-m_h^{(j)}|_{\L^2}^2 |m^\epsilon(s)|_{\H^1}^4 \\
		&\quad + \la \frac{1}{4}+\delta \ra \int_D |m_h^{(j-1)}|^2 |e_h^{(j)}|^2 \ dx \\
		&\leq \delta |e_h^{(j)}|_{\L^2}^2 + \delta |\nabla e_h^{(j)}|_{\L^2}^2 + c\delta^{-1} |m_h^{j-1}-m_h^{(j)}|_{\L^2}^2 \la |m^\epsilon(s)|_{\H^1}^4 + |m^\epsilon(s)|_{\H^2}^2 \ra \\
		&\quad + \la \frac{1}{4}+\delta \ra \int_D |m_h^{(j-1)}|^2 |e_h^{(j)}|^2 \ dx \\
		\sum_{j=1}^n \int_{\tau_{j-1}^{\epsilon,R}}^{\tau_j^{\epsilon,R}} I_{3d}^{(j)} \ ds 
		&\leq \delta \sum_{j=1}^n |e_h^{(j)}|^2_{\L^2} \Delta \tau_j^{\epsilon,R} + \delta \sum_{j=1}^n |\nabla e_h^{(j)}|^2_{\L^2} \Delta \tau_j^{\epsilon,R}\\
		&\quad + c\delta^{-1} \la |m^\epsilon|_{L^\infty(0,T;\H^1)}^4 + |m^\epsilon|_{L^\infty(0,\tau^{\epsilon,R};\H^2)}^2\ra \la \sum_{j=1}^n |m_h^{(j)} - m_h^{(j-1)}|_{\L^2}^2 \ra \Delta t \\
		&\quad + \la \frac{1}{4}+\delta \ra \sum_{j=1}^n \int_D |m_h^{(j-1)}|^2 |e_h^{(j)}|^2 \ dx \ \Delta \tau_j^{\epsilon,R},
	\end{align*}
	where the last term on the right-hand side can be absorbed by its negative counterpart in $ I_{3b}^{(j)} $ for $ \delta \leq \frac{1}{4} $. 
	Finally, for $ I_{3e}^{(j)} $, we have
	\begin{align*}
		I_{3e}^{(j)} 
		&= \lb \langle m_h^{(j-1)} + m^\epsilon(s), m^\epsilon(\tau_j^{\epsilon,R})-m^\epsilon(s) \rangle m^\epsilon(s), P_h e_h^{(j)} \rb_{\L^2} \\
		&\leq \frac{1}{2} |m^\epsilon(s)|^2_{\L^\infty} |m_h^{(j-1)} + m^\epsilon(s)|_{\L^4}^2 |m^\epsilon(\tau_j^{\epsilon,R})-m^\epsilon(s)|_{\L^4}^2 + \frac{1}{2}|e_h^{(j)}|^2_{\L^2} \\
		&\leq c|m^\epsilon(s)|^2_{\H^2} \la |m_h^{(j-1)}|_{\H^1}^2 + |m^\epsilon(s)|_{\H^1}^2 \ra |m^\epsilon(\tau_j^{\epsilon,R})-m^\epsilon(s)|_{\H^1}^2 + \frac{1}{2}|e_h^{(j)}|^2_{\L^2} \\
		&\leq c|m^\epsilon(s)|^2_{\H^2} \la |m_h^{(j-1)}|_{\H^1}^2 + |m^\epsilon(s)|_{\H^1}^2 \ra |m^\epsilon|_{\mathcal{C}^\alpha([0,\tau^{\epsilon,R}];\H^1)}^2 |\tau_j^{\epsilon,R}-s|^{2\alpha} + \frac{1}{2}|e_h^{(j)}|^2_{\L^2} \\
		\sum_{j=1}^n \int_{\tau_{j-1}^{\epsilon,R}}^{\tau_j^{\epsilon,R}} I_{3e}^{(j)} \ ds 
		&\leq c|m^\epsilon|_{\mathcal{C}^\alpha([0,\tau^{\epsilon,R}];\H^1)}^2 (\Delta t)^{2\alpha} |m^\epsilon|^2_{L^\infty(0,\tau^{\epsilon,R};\H^2)} \la \sum_{j=1}^n |m_h^{(j-1)}|^2_{\H^1} \Delta \tau_{j-1}^{\epsilon,R} + |m^\epsilon|^2_{\L^\infty(0,T;\H^1)} T \ra \\
		&\quad + \frac{1}{2} \sum_{j=1}^n |e_h^{(j)}|^2_{\L^2} \Delta \tau_j^{\epsilon,R}. 
	\end{align*}	
	Hence, for $ \delta < \frac{1}{4} $ and $ \alpha < \frac{1}{2} $, we have
	\begin{equation}\label{eq: I3 est}
		\begin{aligned}
			\sum_{j=1}^n \int_{\tau_{j-1}^{\epsilon,R}}^{\tau_j^{\epsilon,R}} I_3^{(j)} \ ds 
			&\leq \la 2+\delta \ra \sum_{j=1}^n |e_h^{(j)}|^2_{\L^2} \Delta \tau_j^{\epsilon,R} + \delta \sum_{j=1}^n |\nabla e_h^{(j)}|^2_{\L^2} \Delta \tau_j^{\epsilon,R} + (\Delta t)^{2\alpha} \psi^{(n)}_3
			+ h \varphi^{(n)}_3, 
		\end{aligned}
	\end{equation}
	where  
	\begin{equation}\label{def: psij}
		\begin{aligned}
			\psi^{(n)}_3
			&= c |m^\epsilon|_{\mathcal{C}^\alpha([0,\tau^{\epsilon,R}]; \L^2)}^2 
			+ c \delta^{-1} |m^\epsilon|^2_{\mathcal{C}^\alpha([0,\tau^{\epsilon,R}]; \H^1)} \la \sum_{j=1}^n |m_h^{(j-1)}|_{\H^1}^2 \Delta \tau_{j-1}^{\epsilon,R} \ra \\
			&\quad + c \delta^{-1} \la |m^\epsilon|_{L^\infty(0,T;\H^1)}^4 + |m^\epsilon|_{L^\infty(0,\tau^{\epsilon,R};\H^2)}^2\ra \la \sum_{j=1}^n |m_h^{(j)} - m_h^{(j-1)}|_{\L^2}^2 \ra (\Delta t)^{1-2\alpha} \\
			&\quad + c |m^\epsilon|_{\mathcal{C}^\alpha([0,\tau^{\epsilon,R}];\H^1)}^2 |m^\epsilon|^2_{L^\infty(0,\tau^{\epsilon,R};\H^2)} \la \sum_{j=1}^n |m_h^{(j-1)}|^2_{\H^1} \Delta \tau_{j-1}^{\epsilon,R} + |m^\epsilon|^2_{\L^\infty(0,T;\H^1)} \ra, \\
			\varphi^{(n)}_3 
			&= c\la \max_{j \leq n} |m_h^{(j)}|_{\L^2} + |m^\epsilon|_{L^\infty(0,T;\L^2)} + |m^\epsilon|_{L^\infty(0,T; \H^1)} \ra |m^\epsilon|_{L^\infty(0,\tau^{\epsilon,R};\H^2)} \sum_{j=1}^n |m_h^{(j-1)}|^2_{\H^1} \Delta \tau_{j-1}^{\epsilon,R} \\
			&\quad + c \la |m^\epsilon|_{L^\infty(0,T; \H^1)}^6 + \int_0^{\tau^{\epsilon,R}} |m^\epsilon(s)|^2_{\H^2} \ ds \ra. 
		\end{aligned}
	\end{equation}
	By Theorem \ref{Theorem: m^e}, we have
	\begin{align*}
		\E \left[ \max_{n \leq N} \psi^{(n)}_3 \right] \leq c_3\epsilon^{-\frac{13}{3}} , \quad
		\E \left[ \max_{n \leq N} \varphi^{(n)}_3 \right] \leq c_3\epsilon^{-\frac{5}{6}}. 
	\end{align*}

	\underline{\bf Estimate $ I_4 $.}
	\begin{align*}
		I_4^{(j)}
		&= \frac{1}{2}\lb \la m^\epsilon(s) - m_h^{(j-1)} \ra \times g_k, g_k \times P_h e_h^{(j)} \rb_{\L^2} \\
		&= \frac{1}{2}\lb \la m^\epsilon(s) - m^\epsilon(\tau_j^{\epsilon,R}) \ra \times g_k + \la m_h^{(j)} - m_h^{(j-1)} \ra \times g_k, g_k \times P_h e_h^{(j)} \rb_{\L^2} \\
		&\quad + \frac{1}{2} \lb e_h^{(j)} \times g_k, g_k \times e_h^{(j)} + g_k \times \la P_h e_h^{(j)} - e_h^{(j)} \ra \rb_{\L^2} \\
		&\leq c |g_k|_{\L^\infty}^2 \la |m^\epsilon|_{\mathcal{C}^\alpha([0,\tau^{\epsilon,R}];\L^2)}^2 |s - \tau_j^{\epsilon,R}|^{2\alpha} 
		+ |m_h^{(j)} - m_h^{(j-1)}|_{\L^2}^2 
		+ h |m^\epsilon(\tau_j^{\epsilon,R})|_{\H^1}^2 \ra \\	
		&\quad + c |g_k|_{\L^\infty}^2 |e_h^{(j)}|_{\L^2}^2, \\ 
		\sum_{j=1}^n \sum_{k=1}^\infty \int_{\tau_{j-1}^{\epsilon,R}}^{\tau_j^{\epsilon,R}} I_4^{(j)} \ ds 
		&\leq 
		c\sum_{j=1}^n |e_h^{(j)}|_{\L^2}^2 \Delta \tau_j^{\epsilon,R} 
		+ (\Delta t)^{2\alpha} \psi_4^{(n)} 
		+ h \varphi_4^{(n)}, 
	\end{align*}
	where
	\begin{align*}
		\psi_4^{(n)} &= c \la |m^\epsilon|_{\mathcal{C}^\alpha([0,\tau^{\epsilon,R}];\L^2)}^2 + \sum_{j=1}^n |m_h^{(j)}-m_h^{(j-1)}|_{\L^2}^2 (\Delta t)^{1-2\alpha} \ra, \\
		\varphi_4^{(n)} &= c |m^\epsilon|^2_{L^\infty(0,\tau^{\epsilon,R};\H^1)},
	\end{align*}
	and 
	\begin{align*}
		\E \left[ \max_{n \leq N} \la \psi_4^{(n)} + \varphi_4^{(n)} \ra \right] \leq c_4. 
	\end{align*}

	\underline{\bf Estimate $ I_5 $.}
	\begin{align*}
		I_5^{(j)} 
		&= \int_{\tau_{j-1}^{\epsilon,R}}^{\tau_j^{\epsilon,R}} \lb P_h e_h^{(j)}, (m^\epsilon(s) - m_h^{(j-1)}) \times dW(s) \rb_{\L^2} \\
		&= \int_{\tau_{j-1}^{\epsilon,R}}^{\tau_j^{\epsilon,R}} \lb P_h \la e_h^{(j)} - e_h^{(j-1)} \ra, (m^\epsilon(s) - m_h^{(j-1)}) \times dW(s) \rb_{\L^2} \\
		&\quad + \int_{\tau_{j-1}^{\epsilon,R}}^{\tau_j^{\epsilon,R}} \lb P_h e_h^{(j-1)}, (m^\epsilon(s) - m_h^{(j-1)}) \times dW(s) \rb_{\L^2} \\
		&= I_{5a}^{(j)} + I_{5b}^{(j)}. 
	\end{align*}
	Then, 
	\begin{align*}
		I_{5a}^{(j)} 
		&= \lb P_h \la e_h^{(j)} - e_h^{(j-1)} \ra, \int_{\tau_{j-1}^{\epsilon,R}}^{\tau_j^{\epsilon,R}} (m^\epsilon(s) - m_h^{(j-1)}) \times dW(s) \rb_{\L^2} \\
		&\leq \delta \left|e_h^{(j)} - e_h^{(j-1)} \right|^2_{\L^2} + c\delta^{-1} \left| \int_{\tau_{j-1}^{\epsilon,R}}^{\tau_j^{\epsilon,R}} (m^\epsilon(s) - m_h^{(j-1)}) \times dW(s) \right|^2_{\L^2},
	\end{align*}
	and
	\begin{align*}
		\E \left[ \max_{n \leq N} \left| \sum_{j=1}^n I_{5a}^{(j)} \right|  \right]
		&\leq \E \left[ \max_{n \leq N} \sum_{j=1}^n \la \delta |e_h^{(j)}-e_h^{(j-1)}|_{\L^2}^2 + c\delta^{-1} \left| \int_{\tau_{j-1}^{\epsilon,R}}^{\tau_j^{\epsilon,R}} (m^\epsilon(s) - m_h^{(j-1)}) \times dW(s) \right|^2_{\L^2} \ra \right] \\
		&\leq \delta \E \left[ \sum_{j=1}^N |e_h^{(j)}-e_h^{(j-1)}|^2_{\L^2} \right] + c\delta^{-1} \E \left[ \sum_{j=1}^N \left| \int_{\tau_{j-1}^{\epsilon,R}}^{\tau_j^{\epsilon,R}} (m^\epsilon(s) - m_h^{(j-1)}) \times dW(s) \right|^2_{\L^2} \right],
	\end{align*}	
	where
	\begin{align*}
		&\E \left[ \sum_{j=1}^N \left| \int_{\tau_{j-1}^{\epsilon,R}}^{\tau_j^{\epsilon,R}} (m^\epsilon(s) - m_h^{(j-1)}) \times dW(s) \right|^2_{\L^2} \right] \\
		&= \E \left[ \sum_{j=1}^N \int_{t_{j-1}}^{t_j} |m^\epsilon(s) - m_h^{(j-1)}|^2_{\L^2} \mathbbm{1}(t_j \leq \tau^{\epsilon,R}) \ ds \right] \\
		&\leq 2\E \left[ \sum_{j=1}^N |e_h^{(j-1)}|^2_{\L^2} \mathbbm{1}(t_j \leq \tau^{\epsilon,R}) \Delta t \right] + 2\E \left[ \sum_{j=1}^N \int_{t_{j-1}}^{t_j} |m^\epsilon(s) -m^\epsilon(\tau_{j-1}^{\epsilon,R})|_{\L^2}^2 \mathbbm{1}(t_j \leq \tau^{\epsilon,R}) \ ds \right] \\
		&\leq 2\E \left[ \sum_{j=1}^N |e_h^{(j-1)}|^2_{\L^2} \Delta \tau_j^{\epsilon,R} \right] + 2\E \left[ |m^\epsilon|_{\mathcal{C}^\alpha([0,\tau^{\epsilon,R}]; \L^2)}^2 (\Delta t)^{2\alpha} T \right].
	\end{align*}
	Note that
	\begin{align*}
		\sum_{j=1}^n I_{5b}^{(j)}
		&= \sum_{j=1}^n \int_{\tau_{j-1}^{\epsilon,R}}^{\tau_j^{\epsilon,R}} \lb P_h e_h^{(j-1)}, (m^\epsilon(s) - m_h^{(j-1)}) \times dW(s) \rb_{\L^2} \\
		&= \int_0^{n \Delta t} \sum_{j=1}^n \mathbbm{1}_{s \in [\tau_{j-1}^{\epsilon,R},\tau_j^{\epsilon,R})} \sum_{k=1}^\infty \lb P_h e_h^{(j-1)}, (m^\epsilon(s) - m_h^{(j-1)}) \times g_k \rb_{\L^2} dW_k(s), 
	\end{align*}
	and with the assumption \eqref{A: g bound H2},
	\begin{align*}
		\sum_{k=1}^\infty \lb P_h e_h^{(j-1)}, (m^\epsilon(s) - m_h^{(j-1)}) \times g_k \rb_{\L^2}^2 
		&\leq c |e_h^{(j-1)}|_{\L^2}^2 \ |m^\epsilon(s) - m_h^{(j-1)}|_{\L^2}^2.
	\end{align*}
	By Burkholder-Davis-Gundy inequality, 
	\begin{align*}
		\E \left[ \sup_{n \leq N} \left| \sum_{j=1}^n I_{5b}^{j} \right| \right]
		&\leq c\E \left[ \la \int_0^T \sum_{j=1}^N \mathbbm{1}_{s \in [\tau_{j-1}^{\epsilon,R},\tau_j^{\epsilon,R})} |e_h^{(j-1)}|_{\L^2}^2 \ |m^\epsilon(s) - m_h^{(j-1)}|_{\L^2}^2 \ ds \ra^{\frac{1}{2}} \right] \\
		&\leq c\E \left[ \max_{j \leq N} |e_h^{(j-1)}|_{\L^2} \la \int_0^T \sum_{j=1}^N \mathbbm{1}_{s \in [\tau_{j-1}^{\epsilon,R},\tau_j^{\epsilon,R})} |m^\epsilon(s) - m_h^{(j-1)}|_{\L^2}^2 \ ds \ra^{\frac{1}{2}} \right] \\
		&\leq \delta \E \left[ \max_{j \leq N} |e_h^{(j-1)}|_{\L^2}^2 \right] 
		+ c\delta^{-1} \E \left[ \int_0^T \sum_{j=1}^N \mathbbm{1}_{s \in [\tau_{j-1}^{\epsilon,R},\tau_j^{\epsilon,R})} |m^\epsilon(s) - m_h^{(j-1)}|_{\L^2}^2 \ ds \right],
	\end{align*}
	where
	\begin{align*}
		&\E \left[ \int_0^T \sum_{j=1}^N \mathbbm{1}_{s \in [\tau_{j-1}^{\epsilon,R},\tau_j^{\epsilon,R})} |m^\epsilon(s) - m_h^{(j-1)}|_{\L^2}^2 \ ds \right] \\
		&\leq 2\E \left[ \int_0^T \sum_{j=1}^N \mathbbm{1}_{s \in [\tau_{j-1}^{\epsilon,R},\tau_j^{\epsilon,R})} |m^\epsilon(s) - m^\epsilon(\tau_{j-1}^{\epsilon,R})|_{\L^2}^2 \ ds \right] + 2\E \left[ \int_0^T \sum_{j=1}^N \mathbbm{1}_{s \in [\tau_{j-1}^{\epsilon,R},\tau_j^{\epsilon,R})} |e_h^{(j-1)}|_{\L^2}^2 \ ds \right] \\
		&= 2\E \left[ \sum_{j=1}^N  \int_{\tau_{j-1}^{\epsilon,R}}^{\tau_j^{\epsilon,R}} |m^\epsilon(s) - m^\epsilon(\tau_{j-1}^{\epsilon,R})|_{\L^2}^2 \ ds \right] + 2\E \left[ \sum_{j=1}^N |e_h^{(j-1)}|_{\L^2}^2 \Delta \tau_j^{\epsilon,R} \right] \\
		&\leq 2\E \left[ |m^\epsilon|_{C^\alpha([0,\tau^{\epsilon,R}];\L^2)}^2 (\Delta t)^{2\alpha} T \right] + 2\E \left[ \sum_{j=1}^N |e_h^{(j-1)}|_{\L^2}^2 \Delta \tau_j^{\epsilon,R} \right].
	\end{align*}
	Hence, we have
	\begin{equation}\label{eq: I5 est}
		\begin{aligned}
			\E \left[ \max_{n \leq N} \left| \sum_{j=1}^n I_5^{(j)} \right| \right] 
			&\leq 
			\delta \E \left[ \sum_{j=1}^N |e_h^{(j)}-e_h^{(j-1)}|^2_{\L^2} \right] + \delta \E \left[ \sup_{j \leq N} |e_h^{(j-1)}|_{\L^2}^2 \right] \\
			&\quad + c\delta^{-1} \la \E \left[ \sum_{j=1}^N |e_h^{(j-1)}|^2_{\L^2} \Delta \tau_j^{\epsilon,R} \right] + \E \left[ |m^\epsilon|_{\mathcal{C}^\alpha([0,\tau^{\epsilon,R}]; \L^2)}^2 (\Delta t)^{2\alpha} T \right] \ra \\
			&\leq 
			\delta \E \left[ \sum_{j=1}^N |e_h^{(j)}-e_h^{(j-1)}|^2_{\L^2} \right] + \delta \E \left[ \max_{j \leq N} |e_h^{(j-1)}|_{\L^2}^2 \right] \\
			&\quad + c_5\delta^{-1} \E \left[ \sum_{j=1}^N |e_h^{(j-1)}|^2_{\L^2} \Delta \tau_{j-1}^{\epsilon,R} \right] 
			+ c_5 \epsilon^{-\frac{8}{3}} (\Delta t)^{2\alpha}. 
		\end{aligned}
	\end{equation}
	
	For a sufficiently small $ \delta $ such that $ 1-7\delta > 0 $ and $ 1-\delta-2\delta^\frac{1}{3}>0 $, we deduce from the equation \eqref{eq: 1/2|e_h^n|^2} and the estimates \eqref{eq: I0 est} -- \eqref{eq: I3 est} that 
	\begin{align*}
		&\frac{1}{2} |e_h^{(n)}|^2_{\L^2} + \la \frac{1}{2}-\delta \ra \sum_{j=1}^n |e_h^{(j)}-e_h^{(j-1)}|^2_{\L^2} \\
		&\quad + (1-7\delta)\sum_{j=1}^n |\nabla e_h^{(j)}|^2_{\L^2} \Delta \tau_j^{\epsilon,R} 
		+ \epsilon (1-\delta-2\delta^\frac{1}{3})\sum_{j=1}^n |\Delta e_h^{(j)}|^2_{\L^2} \Delta \tau_j^{\epsilon,R} \\
		&\leq 
		\frac{1}{2} |e_h^{(0)}|^2_{\L^2} + \eta_2 
		+ h \sum_{i=0}^4 \varphi_i^{(n)} 
		+ \la (\Delta t)^\frac{\alpha}{2}\psi_{1a}^{(n)} + (\Delta t)^\alpha \psi_{2b}^{(n)} + (\Delta t)^{2\alpha} \la \psi_{1b}^{(n)} + \psi_{2a}^{(n)} + \psi^{(n)}_3 + \psi^{(n)}_4 \ra \ra \\
		&\quad + c \delta^{-5} \epsilon^{-3} \sum_{j=1}^n |e_h^{(j)}|_{\L^2}^2 \int_{\tau_{j-1}^{\epsilon,R}}^{\tau_{j}^{\epsilon,R}} \la |m^\epsilon(s)|_{\H^1}^8 + |m^\epsilon(s)|_{\H^1}^{\frac{8}{3}} + |m^\epsilon(s)|_{\H^1}^2 \ra \ ds 
		+ \la c +\delta \ra \sum_{j=1}^n |e_h^{(j)}|^2_{\L^2} \Delta \tau_j^{\epsilon,R} \\
		&\quad + \left| \sum_{j=1}^n I_5^{(j)} \right|.
	\end{align*}
	Taking the maximum over $ n $ and the expectation,
	\begin{align*}
		&\E \left[\la \frac{1}{2}-\delta \ra \max_{n \leq N} |e_h^{(n)}|^2_{\L^2} + \la \frac{1}{2}-2\delta \ra \sum_{j=1}^N |e_h^{(j)}-e_h^{(j-1)}|^2_{\L^2} \right] \\
		&\quad + (1-7\delta) \E \left[ \sum_{j=1}^N |\nabla e_h^{(j)}|^2_{\L^2} \Delta \tau_j^{\epsilon,R} \right]
		+ \epsilon (1-\delta-2\delta^\frac{1}{3}) \E \left[ \sum_{j=1}^N |\Delta e_h^{(j)}|^2_{\L^2} \Delta \tau_j^{\epsilon,R} \right] \\
		&\leq 
		\la \frac{1}{2}+\delta + c_5 \delta^{-1} \Delta t \ra \E \left[|e_h^{(0)}|^2_{\L^2}\right] + \E \left[\eta_2 \right] \\
		&\quad + \la h+(\Delta t)^{\frac{\alpha}{2}} \ra \E \left[ \sum_{i=0}^4 \max_{n \leq N} \varphi_i^{(n)} 
		+ \la \psi_{1a}^{(n)} + \psi_{1b}^{(n)} + \psi_{2a}^{(n)} + \psi_{2b}^{(n)} + \psi^{(n)}_3 + \psi^{(n)}_4 + c_4 \epsilon^{-\frac{8}{3}} \ra \right] \\
		&\quad + c \delta^{-5} \epsilon^{-3} \E \left[\sum_{j=1}^N |e_h^{(j)}|_{\L^2}^2 \int_{\tau_{j-1}^{\epsilon,R}}^{\tau_{j}^{\epsilon,R}} \la |m^\epsilon(s)|_{\H^1}^8 + |m^\epsilon(s)|_{\H^1}^{\frac{8}{3}} + |m^\epsilon(s)|_{\H^1}^2 \ra \ ds \right] \\
		&\quad + \la c +\delta + c_5\delta^{-1} \ra \E \left[ \sum_{j=1}^N |e_h^{(j)}|^2_{\L^2} \Delta \tau_j^{\epsilon,R} \right] \\
		&\leq 
		c\la h+(\Delta t)^{\frac{\alpha}{2}} \ra \epsilon^{-\frac{13}{3}} 
		+ c \epsilon^{-3} R^8 \E \left[\sum_{n=1}^N \max_{j\leq n}|e_h^{(j)}|_{\L^2}^2 \right] \Delta t,
	\end{align*}
	where the constant $ c $ may depend on $ C_1, C_g, T, \delta $ but not on $ \epsilon $ and $ R $. 
	Then by Gronwall's lemma, 
	\begin{align*}
		\E \left[\max_{n \leq N} |e_h^{(n)}|^2_{\L^2} \right]
		&\leq c\la h+(\Delta t)^{\frac{\alpha}{2}} \ra \epsilon^{-\frac{13}{3}} 
		+ c \epsilon^{-3} R^8 \E \left[\sum_{n=1}^N \max_{j\leq n}|e_h^{(j)}|_{\L^2}^2 \right] \Delta t \\
		&\leq c\la h+(\Delta t)^{\frac{\alpha}{2}} \ra \epsilon^{-\frac{13}{3}} e^{c \epsilon^{-3} R^8}.
	\end{align*}
	Therefore,
	\begin{align*}
		\E \left[ \max_{n \leq N} |e_h^{(n)}|^2_{\L^2} + \sum_{n=1}^N |\nabla e_h^{(n)}|^2_{\L^2} \Delta \tau_n^{\epsilon,R} \right]
		&\leq 
		c\la h+(\Delta t)^{\frac{\alpha}{2}} \ra \epsilon^{-\frac{13}{3}} 
		\la 1 + e^{c \epsilon^{-3} R^8}+ c \epsilon^{-3} R^8 e^{c \epsilon^{-3} R^8} \ra \\
		&\leq
		c \la h+(\Delta t)^{\frac{\alpha}{2}} \ra \epsilon^{-\frac{22}{3}} R^8 e^{c \epsilon^{-3} R^8},
	\end{align*}
	for a constant $ c $ that is independent of $ h, \Delta t, \epsilon $ and $ R $.

\subsection{Proof of Theorem \ref{Theorem: conver 1D}}
	Fix $ R $. Recall that $ e_h^{(n)} = m(\tau_n^R) - m_h^{(n)} $. 
	The equation of $ e_h^{(n)} $ is similar to \eqref{eq: <ehn,phi>} with $ \epsilon = 0 $ and $ m $ in place of $ m^\epsilon $. We obtain
	\begin{align*}
		&\lb e_h^{(n)}-e_h^{(n-1)}, P_h e_h^{(n)} \rb_{\L^2} \\
		&= - \int_{\tau_{n-1}^R}^{\tau_n^R} \lb \nabla m(s) - \nabla m_h^{(n)}, \nabla P_h e_h^{(n)} \rb_{\L^2} \ ds \\
		&\quad - \int_{\tau_{n-1}^R}^{\tau_n^R} \lb m(s) \times \nabla m(s) - m_h^{(n-1)} \times \nabla m_h^{(n)}, \nabla P_h e_h^{(n)} \rb_{\L^2} \ ds \\
		&\quad - \int_{\tau_{n-1}^R}^{\tau_n^R} \lb (1 + |m(s)|^2) m(s) - (1+|m_h^{(n-1)}|^2) m_h^{(n)}, P_h e_h^{(n)} \rb_{\L^2} \ ds \\
		&\quad + \frac{1}{2} \sum_{k=1}^\infty \int_{\tau_{n-1}^R}^{\tau_n^R} \lb \la m(s)-m_h^{(n-1)} \ra \times g_k, g_k \times P_h e_h^{(n)} \rb_{\L^2} \ ds \\
		&\quad + \int_{t_{n-1}}^{t_n} \lb P_h e_h^{(n)}, \la m(s) - m_h^{(n-1)} \ra \times dW(s) \rb_{\L^2} \frac{\Delta \tau_n^R}{\Delta t} \\
		&= \int_{\tau_{n-1}^R}^{\tau_n^R} \la J_1^{(n)} + J_2^{(n)} + J_3^{(n)} \ra \ ds + \sum_{k=1}^\infty \int_{\tau_{n-1}^R}^{\tau_n^R} J_4^{(n)} \ ds + J_5^{(n)} .
	\end{align*}
	The estimates of $ J_1, J_3, J_4 $ and $ J_5 $ are identical to that of $ I_{1b}, I_3, I_4 $ and $ I_5 $ in the proof of Theorem \ref{Theorem: eh conve}, respectively. 
	The estimate of $ J_2 $ is also similar to that of $ I_2 $:
	\begin{align*}
		J_2^{(j)} 
		&= \lb \la m(\tau_{j-1}^R) - m(s) \ra \times \nabla m(s), \nabla P_h e_h^{(j)} \rb_{\L^2} 
		+ \lb m_h^{(j-1)} \times \la \nabla m(\tau_j^R) - \nabla m(s) \ra, \nabla P_h e_h^{(j)} \rb_{\L^2} \\
		&\quad - \lb e_h^{(j-1)} \times \nabla m(s), \nabla P_h e_h^{(j)} \rb_{\L^2} 
		- \lb m_h^{(j-1)} \times \nabla e_h^{(j)}, \nabla P_h e_h^{(j)} \rb_{\L^2} \\
		&= J_{2a}^{(j)} + J_{2b}^{(j)} + J_{2c}^{(j)} + J_{2d}^{(j)},
	\end{align*}
	where $ J_{2a}, J_{2b} $ and $ J_{2d} $ are estimated in the same way as $ I_{2a}, I_{2b} $ and $ I_{2d} $ in the proof of Theorem \ref{Theorem: eh conve}. 
	For $ J_{2c} $, we only require $ \H^1 $-norm of the error:
	\begin{align*}
		J_{2c}^{(j)}
		&= - \lb e_h^{(j-1)} \times \nabla m, \nabla P_h e_h^{(j)} \rb_{\L^2} \\
		&\leq |e_h^{(j-1)}|_{\L^\infty} |\nabla m|_{\L^2}  |\nabla e_h^{(j)}|_{\L^2} \\
		&\leq \la |e_h^{(j-1)}|_{\L^2}^\frac{1}{2} |\nabla e_h^{(j-1)}|_{\L^2}^\frac{1}{2} + |e_h^{(j-1)}|_{\L^2} \ra |\nabla e_h^{(j)}|_{\L^2} |\nabla m|_{\L^2} \\
		&\leq c |e_h^{(j-1)}|_{\L^2}^2 \la |\nabla m|_{\L^2}^4 + |\nabla m|_{\L^2}^2 \ra + \frac{3}{4}\delta |\nabla e_h^{(j-1)}|_{\L^2}^\frac{2}{3} |\nabla e_h^{(j)}|_{\L^2}^\frac{4}{3} + \frac{1}{2}\delta |\nabla e_h^{(j)}|_{\L^2}^2 \\
		&\leq c |e_h^{(j-1)}|_{\L^2}^2 \la |\nabla m|_{\L^2}^4 + |\nabla m|_{\L^2}^2 \ra 
		+ \delta \la |\nabla e_h^{(j-1)}|_{\L^2}^2 + |\nabla e_h^{(j)}|_{\L^2}^2 \ra, \\
		\sum_{j=1}^n \int_{\tau_{j-1}^R}^{\tau_j^R} J_{2c}^{(j)} \ ds
		&\leq 
		\delta |\nabla e_h^{(0)}|^2_{\L^2} +
		\delta \sum_{j=1}^n |\nabla e_h^{(j)}|^2_{\L^2} \Delta \tau_j^R \\
		&\quad + c \sum_{j=1}^n |e_h^{(j-1)}|_{\L^2}^2 \int_{\tau_{j-1}^R}^{\tau_j^R} \la |\nabla m(s) |_{\L^2}^4 + |\nabla m(s) |_{\L^2}^2 \ra \ ds. 
	\end{align*} 	
	The rest of the arguments follow similarly. 
	Using the $ \H^1 $-stopping time and Gronwall's lemma, we deduce
	\begin{align*}
		\E \left[ \max_{n \leq N} |e_h^{(n)}|^2_{\L^2} + \sum_{n=1}^N |\nabla e_h^{(n)}|^2_{\L^2} \Delta \tau_n^R \right]
		&\leq c^* e^{c^* R^4} \la h + (\Delta t)^\alpha \ra,
	\end{align*}
	for some constant $ c^*>0 $ independent of $ R $. 
	We write $ c^\dagger $ instead of $ c $ for the constant in \eqref{eq: P(H1exit) 1D}. 
	Without loss of generality, assume that $ c^* \geq c^\dagger $. 
	Then fix $ q,\beta \in (0,1) $ and choose
	\begin{align*}
		R = R(h,\Delta t) = \la -\frac{q\beta}{c^*} \ln (h + (\Delta t)^\alpha) \ra^\frac{1}{4}.
	\end{align*}
	We obtain 
	\begin{align*}
		&\P \la \max_{n \leq N} |m(t_n) - m_h^{(n)}|^2_{\L^2} + \sum_{n=1}^N |\nabla m(t_n) - \nabla m_h^{(n)}|^2_{\L^2} \Delta t > \gamma (h + (\Delta t)^\alpha)^{1-\beta} \ra \\
		&\leq \P \la \max_{n \leq N} |m(\tau^{R}_n) - m_h^{(n)}|^2_{\L^2} + \sum_{n=1}^N |\nabla m(\tau^{R}_n) - \nabla m_h^{(n)}|^2_{\L^2} \Delta \tau^{R} > \gamma (h + (\Delta t)^\alpha)^{1-\beta} \ra 
		+ \P(\tau^{R} < T) \\
		&\leq c^* \gamma^{-1} \la h+(\Delta t)^{\alpha} \ra^{\beta(1- q)} + c^* (R(h,\Delta t))^{-1} 
		\to 0,
	\end{align*}
	as $ h,\Delta t \to 0 $.

\appendix
\section{}
We quote the existence, uniqueness and regularity results from \cite[Theorems 2.2 and 2.3]{BGL_sLLB} in the following theorems.
\begin{theorem}\label{Theorem: sLLB existence}
		Let $ d=1,2 $. 
		Assume $ |m_0|_{\H^1} < C_1 $ for a certain $ C_1 > 0 $. 
		Then there exists a weak martingale solution $ (\Omega, \mathcal{F}, \mathbb{F}, \P, W, m) $ of \eqref{eq: sLLB} such that 
		\begin{enumerate}[(i)]
			\item 
			for every $ p \in [1,\infty) $, $ \alpha \in (0, \frac{1}{2}) $, 
			\begin{align*}
				m \in L^p(\Omega; L^\infty(0,T; \H^1) \cap L^2(0,T; \H^2)), \\
				\E \left[ |m|_{W^{\alpha, p}(0,T; \L^2)}^q + |m|_{L^\infty(0,T; \H^1) \cap L^2(0,T; \H^2)}^p \right] < c,
			\end{align*}
			and for every $ q \in [1,\frac{4}{3}) $, 
			\begin{align*}
				\E \left[ \la \int_0^T |m(t) \times \Delta m(t)|_{\L^2}^q \ dt \ra^p \right] < c,
			\end{align*}
			where $ c $ is a positive constant depending on $ p, C_1 $ and $ g $; 
			
			\item 
			the following equality holds in $ \L^2 $: 
			\begin{equation}\label{eq: sLLB integral form}
				\begin{aligned}
					m(t)
					&= m(0) + \kappa_1 \int_0^t \Delta m(s) \ ds + \gamma \int_0^t m(s) \times \Delta m(s) \ ds - \kappa_2 \int_0^t (1 + \mu |m|^2) m(s) \ ds \\
					&\quad + \sum_{k=1}^\infty \int_0^t \la \gamma m(s) \times g_k + \kappa_1 g_k \ra \circ dW_k(s);
				\end{aligned}
			\end{equation}
			
			\item 
			for every $ \bar{\alpha} \in [0,\frac{1}{4}) $ and $ \beta \in [0,\frac{1}{2}] $, 
			\begin{align*}
				m \in \mathcal{C}^{\bar{\alpha}}([0,T]; \L^2) \cap \mathcal{C}^\beta([0,T]; \L^{\frac{3}{2}}) \cap \mathcal{C}([0,T]; \H^1_w), \quad \P\text{-a.s.}
			\end{align*}
			where $ \mathcal{C}([0,T]; \H^1_w) $ is the space of weakly (subscript $ w $) continuous functions $ u: [0,T] \to \H^1 $, endowed with the weakest topology such that for any $ g \in \H^1 $, the mapping $ f: \mathcal{C}([0,T]; \H^1_w) \to \mathcal{C}([0,T]; \R) $ given by $ f(u) = \lb u, g \rb_{\H^1} $ is continuous. 
		\end{enumerate}
	\end{theorem}
 
 	\begin{theorem}\label{Theorem: sLLB path uniqueness}
		Assume that $ (\Omega, \mathcal{F}, \mathbb{F}, \P, W, m_1) $ and $ (\Omega, \mathcal{F}, \mathbb{F}, \P, W, m_2) $ are two weak martingale solutions to \eqref{eq: sLLB integral form} such that for $ i=1,2 $,
		\begin{enumerate}[(i)]
			\setlength{\parskip}{0cm}
			\item 
			$ m_1(0) = m_2(0) = m_0 $; 
			
			\item 
			the paths of $ m_i $ lie in $ L^\infty(0,T; \H^1) \cap L^2(0,T; \H^2) $; 
			
			\item 
			each $ m_i $ satisfies \eqref{eq: sLLB integral form}. 
		\end{enumerate} 
		Then, for $ \P $-a.e. $ \omega \in \Omega $, $ m_1(\cdot, \omega) = m_2(\cdot, \omega) $. 
	\end{theorem}


\begin{thebibliography}{99}	
		\bibitem{BreitProhl}
		D. Breit and A. Prohl. 
		Error analysis for 2D stochastic Navier--Stokes equations in bounded domains. 
		arXiv preprint 
		arXiv:2109.06495 (2021).
		
		\bibitem{BGL_sLLB}
		Z. Brze{\'z}niak, B. Goldys  and K.-N. Le.
		Existence of a unique solution and invariant measures for the stochastic Landau-Lifshitz-Bloch equation.
		Journal of Differential Equations. 
		269 (2020), 9471--9507.
		
		\bibitem{Red_book}
		G. Da Prato and J. Zabczyk. 
		Stochastic Equations in Infinite Dimensions. 
		Cambridge university press, 2014.
		
		\bibitem{dunst}
		T. Dunst. 
		Convergence with rates for a time-discretization of the Stochastic Landau–Lifschitz–Gilbert equation. 
		IMA Journal of Numerical Analysis.
		35 (2015), 615--651.
		
		
		\bibitem{Garanin2004}
		D. A. Garanin and O. Chubykalo-Fesenko 
		Thermal fluctuations and longitudinal relaxation of single-domain magnetic particles at elevated temperatures. 
		Phys. Rev. B, 70 (2004), 212409.
		
		\bibitem{Le2016}
		K.-N. Le.
		Weak solutions of the Landau--Lifshitz--Bloch equation. 
		Journal of Differential Equations.
		261 (2016), 6699 --6717.
		 
		\bibitem{JiangWang2019}
		S. Jiang, Q. Ju and H. Wang.
		Martingale weak solutions of the stochastic Landau--Lifshitz--Bloch equation.
		Journal of Differential Equations. 
		266 (5) (2019), pp. 2542--2574
		 
		\bibitem{Evan2012}
		R.F.L. Evans, D. Hinzke, U. Atxitia, U. Nowak, R.W. Chantrell and O. Chubykalo-Fesenko.
		Stochastic form of the Landau--Lifshitz--Bloch equation.
		Phys. Rev. B, 85 (Jan 2012), Article 014433.
		
		
		%
		
	\end{thebibliography}
\end{document}